%% file: tight.tex
\documentclass[12pt]{amsart}
\usepackage{amscd,amssymb}%,showkeys}
\usepackage[mathscr]{eucal}
\usepackage[english]{babel}
\input pstricks
\input pst-node

\usepackage{pstricks, pst-node}
\input xy
\usepackage[all]{xy}

\usepackage{amsmath}
\usepackage[latin1]{inputenc}%
\usepackage{a4}

\newlength{\itemlaenge}

\begin{document}

%%%%%%%%%% Definitionen von Umgebungen %%%%%%%%%%%%%%%%%%%%%%%%

\newtheoremstyle{mytheorem}% name
  {}%      Space above, empty = `usual value'
  {}%      Space below
  {\slshape}% Body font
  {}%         Indent amount (empty = no indent, \parindent = para indent)
  {\scshape}% Thm head font
  {.}%        Punctuation after head
  { }%     Space after thm head: " " = normal interword space;
        %       \newline = linebreak
  {}% Thm head spec

\newtheoremstyle{mydefinition}% name
  {}%      Space above, empty = `usual value'
  {}%      Space below
  {\upshape}% Body font
  {}%         Indent amount (empty = no indent, \parindent = para indent)
  {\scshape}% Thm head font
  {.}%        Punctuation after thm head
  { }%     Space after thm head: " " = normal interword space;
        %       \newline = linebreak
  {}% Thm head spec

\theoremstyle{mytheorem}
\newtheorem{lemma}{Lemma}[section]
\newtheorem{prop}[lemma]{Proposition}
\newtheorem*{prop*}{Proposition}
\newtheorem{prop_intro}{Proposition}
\newtheorem{cor}[lemma]{Corollary}
\newtheorem{cor_intro}[prop_intro]{Corollary}
\newtheorem{thm}[lemma]{Theorem}
\newtheorem{thm_intro}[prop_intro]{Theorem}
\newtheorem*{thm*}{Theorem}
\theoremstyle{mydefinition}
\newtheorem{rem}[lemma]{Remark}
\newtheorem*{rem*}{Remark}
\newtheorem{rem_intro}[prop_intro]{Remark}
\newtheorem{rems_intro}[prop_intro]{Remarks}
\newtheorem*{notation*}{Notation}
\newtheorem{notation}{Notation}
\newtheorem*{warning*}{Warning}
\newtheorem{rems}[lemma]{Remarks}
\newtheorem{defi}[lemma]{Definition}
\newtheorem*{defi*}{Definition}
\newtheorem{defi_intro}[prop_intro]{Definition}
\newtheorem{defis}[lemma]{Definitions}
\newtheorem{exo}[lemma]{Example}
\newtheorem{exos}[lemma]{Examples}
\newtheorem{exo_intro}[prop_intro]{Example}
\newtheorem{exos_intro}[prop_intro]{Examples}
\newtheorem*{fexo_intro}{Fundamental Example}
\newtheorem{nonexo}[lemma]{Non-Example}

\numberwithin{equation}{section}

\newcommand{\bqn}{\begin{eqnarray*}}
\newcommand{\eqn}{\end{eqnarray*}}
\newcommand{\bq}{\begin{eqnarray}}
\newcommand{\eq}{\end{eqnarray}}
\newcommand{\ba}{\begin{aligned}}
\newcommand{\ea}{\end{aligned}}
\newcommand{\be}{\begin{enumerate}}
\newcommand{\ee}{\end{enumerate}}

\newcommand{\bibURL}[1]{{\unskip\nobreak\hfil\penalty50{\tt#1}}}

\def\ti{-\allowhyphens}
\newcommand{\thismonth}{\ifcase\month % case 0 --- impossible!
  \or January\or February\or March\or April\or May\or June%
  \or July\or August\or September\or October\or November%
  \or December\fi}
\newcommand{\thismonthyear}{{\thismonth} {\number\year}}
\newcommand{\thisdaymonthyear}{{\number\day} {\thismonth} {\number\year}}
%
%
%
%%%%%%%%%%%%%%%%ù Definitionen von Symbolen %%%%%%%%%%%%%%%%%%%%%%ù

\newcommand{\Stab}{\operatorname{Stab}}
\newcommand{\sign}{\operatorname{sign}}
\newcommand{\aut}{\operatorname{Aut}}
\newcommand{\End}{\operatorname{End}}
\newcommand{\Is}{\operatorname{Isom}}
\newcommand{\SU}{\operatorname{SU}}
\newcommand{\SL}{\operatorname{SL}}
\newcommand{\Sp}{\operatorname{Sp}}
\newcommand{\ko}{\operatorname{k}}
\newcommand{\PU}{\operatorname{PU}}
\newcommand{\Image}{\operatorname{Image}}
\newcommand{\Lie}{\operatorname{Lie}}
\newcommand{\codim}{\operatorname{codim}}
\newcommand{\esssup}{\operatorname{ess\,sup}}
\newcommand{\supp}{\operatorname{supp}}
\newcommand{\BB}{{\mathbb B}}
\newcommand{\CC}{{\mathbb C}}
\newcommand{\DD}{{\mathbb D}}
\newcommand{\FF}{{\mathbb F}}
\newcommand{\HH}{{\mathbb H}}
\newcommand{\GG}{{\mathbb G}}
\newcommand{\KK}{{\mathbb K}}
\newcommand{\NN}{{\mathbb N}}
\newcommand{\PP}{{\mathbb P}}
\newcommand{\QQ}{{\mathbb Q}}
\newcommand{\RR}{{\mathbb R}}
\newcommand{\TT}{{\mathbb T}}
\newcommand{\ZZ}{{\mathbb Z}}

\newcommand{\Aa}{{\mathcal A}}
\newcommand{\Bb}{{\mathcal B}}
\newcommand{\Cc}{{\mathcal C}}
\newcommand{\Dd}{{\mathcal D}}
\newcommand{\Ee}{{\mathcal E}}
\newcommand{\Ff}{{\mathcal F}}
\newcommand{\Hh}{{\mathcal H}}
\newcommand{\Jj}{{\mathcal J}}
\newcommand{\Kk}{{\mathcal K}}
\newcommand{\Ll}{{\mathcal L}}
\newcommand{\Mm}{{\mathcal M}}
\newcommand{\Nn}{{\mathcal N}}
\newcommand{\Oo}{{\mathcal O}}
\newcommand{\Qq}{{\mathcal Q}}
\newcommand{\Pp}{{\mathcal P}}
\newcommand{\Rr}{{\mathcal R}}
\newcommand{\Ss}{{\mathcal S}}
\newcommand{\Tt}{{\mathcal T}}
\newcommand{\Vv}{{\mathcal V}}
\newcommand{\Xx}{{\mathcal X}}
\newcommand{\Yy}{{\mathcal Y}}
\newcommand{\Ww}{{\mathcal W}}
\newcommand{\Zz}{{\mathcal Z}}

\newcommand{\frakg}{{\mathfrak g}}
\newcommand{\frakk}{{\mathfrak k}}
\newcommand{\frakh}{{\mathfrak h}}
\newcommand{\fraka}{{\mathfrak a}}
\newcommand{\frake}{{\mathfrak e}}
\newcommand{\frakp}{{\mathfrak p}}
\newcommand{\frako}{{\mathfrak o}}
\newcommand{\fraks}{{\mathfrak s}}
\newcommand{\frakl}{{\mathfrak l}}
\newcommand{\frakm}{{\mathfrak m}}
\newcommand{\frakn}{{\mathfrak n}}
\newcommand{\fraku}{{\mathfrak u}}
\newcommand{\frakc}{{\mathfrak c}}
\newcommand{\frakr}{{\mathfrak r}}
\newcommand{\frakpp}{{\mathfrak p_+}}
\newcommand{\frakpm}{{\mathfrak p_-}}
\newcommand{\frakkc}{{{\mathfrak k}_\CC}}
\newcommand{\frakgc}{{{\mathfrak g}_\CC}}
\newcommand{\frakpc}{{{\mathfrak p}_\CC}}
\newcommand{\frakB}{{\mathfrak B}}

\renewcommand{\a}{\alpha}
\newcommand{\e}{\epsilon}
\newcommand{\eps}{\epsilon}
\renewcommand{\b}{\beta}
\newcommand{\g}{\gamma}
\newcommand{\G}{\Gamma}
\renewcommand{\L}{\Lambda}
\renewcommand{\l}{\lambda}

\newcommand{\dD}{{\mathbf D}}
\newcommand{\gG}{{\mathbf G}}
\newcommand{\hH}{{\mathbf H}}
\newcommand{\pP}{{\mathbf P}}
\newcommand{\lL}{{\mathbf L}}
\newcommand{\qQ}{{\mathbf Q}}
\newcommand{\nN}{{\mathbf N}}
\newcommand{\wW}{{\mathbf W}}
\newcommand{\uU}{{\mathbf U}}

\newcommand{\<}{\langle}
\renewcommand{\>}{\rangle}

\def\ol{\overline}

\def\h{{\rm H}}
\def\hb{{\rm H}_{\rm b}}
\def\ehb{{\rm EH}_{\rm b}}
\def\ha{{\rm H}_{(G,K)}}
\def\hc{{\rm H}_{\rm c}}
\def\hcb{{\rm H}_{\rm cb}}
\def\ehbc{{\rm EH}_{\rm cb}}
\def\linfty{L^\infty}
\def\linftyw{L^\infty_{\rm w*}}
\def\linftya{L^\infty_{\mathrm{w*,alt}}}
\def\la{L^\infty_{\mathrm{alt}}}
\def\cb{{\rm C}_{\rm b}}
\def\binfty{\mathcal B^\infty_{\mathrm alt}}

\def\one{\mathbf{1\kern-1.6mm 1}}
\def\sous#1#2{{\raisebox{-1.5mm}{$#1$}\backslash \raisebox{.5mm}{$#2$}}}
\def\rest#1{{\raisebox{-.95mm}{$\big|$}\raisebox{-2mm}{$#1$}}}
\def\homeo#1{{\sl H\!omeo}^+\!\left(#1\right)}
\def\thomeo#1{\widetilde{{\sl \!H}\!omeo}^+\!\left(#1\right)}
\def\bu{\bullet}
\def\weak{weak-* }
\def\property{\textbf{\rm\textbf A}}
\def\cont{\mathcal{C}}
\def\id{{\it I\! d}}
\def\opposite{^{\rm op}}
\def\oddex#1#2{\left\{#1\right\}_{o}^{#2}}
\def\comp#1{{\rm C}^{(#1)}}
\def\ro{\varrho}
\def\ti{-\allowhyphens}
\def\lra{\longrightarrow}

\def\fix{{\operatorname{Fix}}}
\def\Mat{{\operatorname{Mat}}}
\def\Mod{{\bf Mod}}
\def\zent{{\operatorname{Zent}}}
\def\C{{\operatorname{C}}}
\def\c{{\operatorname{c}}}
\def\sym{{\operatorname{Sym}}}
\def\stab{{\operatorname{Stab}}}
\def\arg{{\operatorname{arg}}}
\def\HTP{{\operatorname{HTP}}}
\def\h2{{\operatorname{H_2}}}
\def\h1{{\operatorname{H_1}}}
\def\pr{{\operatorname{pr}}}
\def\rk{{\operatorname{rank}}}
\def\tr{{\operatorname{tr}}}
\def\codim{{\operatorname{codim}}}
\def\nt{{\operatorname{nt}}}
\def\d{{\operatorname{d}}}
\def\Gr{{\operatorname{Gr}}}
\def\id{{\operatorname{Id}}}
\def\ker{{\operatorname{ker}}}
\def\im{{\operatorname{im}}}
\def\r{\operatorname{r}}

\def\PSL{\operatorname{PSL}}
\def\SL{\operatorname{SL}}
\def\Sp{\operatorname{Sp}}
\def\SU{\operatorname{SU}}
\def\SO{\operatorname{SO}}
\def\Spin{\operatorname{Spin}}
\def\PSU{\operatorname{PSU}}

\def\ad{\operatorname{ad}}
\def\Ad{\operatorname{Ad}}
\def\adg{\operatorname{ad}_\frakg}
\def\adp{\operatorname{ad}_\frakp}
\def\bg{B_\frakg}
\def\creg{C_{\rm reg}}
\def\cs{\check S}
\def\cst{{\check S}^{(3)}}
\def\det{{\operatorname{det}}}
\def\deta{{\operatorname{det}_A}}
\def\diag{{\operatorname{diag}}}
\def\gmodp{\gG(\RR)/\pP(\RR)}
\def\gmodq{\gG(\RR)/\qQ(\RR)}
\def\hom{\operatorname{Hom}}
\def\isp{\operatorname{Is}_{\<\cdot,\cdot\>}}
\def\isptwo{\operatorname{Is}_{\<\cdot,\cdot\>}^{(2)}}
\def\ispth{\operatorname{Is}_{\<\cdot,\cdot\>}^{(3)}}
\def\isf{\operatorname{Is}_F}
\def\isfi{\operatorname{Is}_{F_i}}
\def\isft{\operatorname{Is}_F^{(3)}}
\def\isfit{\operatorname{Is}_{F_i}^{(3)}}
\def\isftwo{\operatorname{Is}_F^{(2)}}
\def\lin{\operatorname{Lin}(L_+,L_-)}
\def\ll{{\Ll_1,\Ll_2}}
\def\kahler{ K\"ahler }
\def\kg{\kappa_G}
\def\kx{\kappa_\Xx}
\def\kxb{\kappa_\Xx^{\rm b}}
\def\kxib{\kappa_{\Xx_i}^{\rm b}}

\def\ko{\kappa_\omega}
\def\kob{\kappa_\omega^{\rm b}}
\def\kgjb{\kappa_{G,J}^{\rm b}}
\def\kgajb{\kappa_{G_a,J}^{\rm b}}
\def\kgojb{\kappa_{G_1,J}^{\rm b}}
\def\kgtjb{\kappa_{G_2,J}^{\rm b}}
\def\kgtb{\kappa_{G_2}^{\rm b}}
\def\khjb{\kappa_{H,J}^{\rm b}}
\def\kljb{\kappa_{L_j}^{\rm b}}
\def\klb{\kappa_{\lambda}^{\rm b}}
\def\kmjb{\kappa_{M,J}^{\rm b}}
\def\kmijb{\kappa_{M_i, J_i}^{\rm b}}
\def\kgijb{\kappa_{G_i, J_i}^{\rm b}}
\def\k{\kappa}
\def\kZ{\kappa_Z}
\def\kZb{\kappa_{Z}^{\rm b}}
\def\kZgob{\kappa_{Z_{\frakg_1}}^{\rm b}}
\def\kZgb{\kappa_{Z_\frakg}^{\rm b}}
\def\kZgtb{\kappa_{Z_{\frakg_2}}^{\rm b}}
\def\kZsb{\kappa_{Z_{\fraks\frakl(2,\RR)}}^{\rm b}}
\def\kZib{\kappa_{Z_i}^{\rm b}}
\def\kgb{\kappa_G^{\rm b}}
\def\kgab{\kappa_{G_a}^{\rm b}}
\def\kgib{\kappa_{G,i}^{\rm b}}
\def\kgbb{\kappa_{G,B}^{\rm b}}
\def\kgibb{\kappa_{G_i,B}^{\rm b}}
\def\khib{\kappa_{H,i}^{\rm b}}
\def\kmb{\kappa_M^{\rm b}}
\def\kib{\kappa_i^{\rm b}}
\def\kibt{\tilde\kappa_i^{\rm b}}
\def\oc{\overline c}
\def\om{\overline m}
\def\pii{\pi_i}
\def\piit{\tilde\pi_i}
\def\psupq{{\rm PSU}(p,q)}
\def\psuvi{{\rm PSU}\big(V,\<\cdot,\cdot\>_i\big)}
\def\rg{r_G}
\def\rx{{{\rm r}_\Xx}}
\def\slv{{\rm SL}(V)}
\def\supq{{\rm SU}(p,q)}
\def\suv{{\rm SU}\big(V,\<\cdot,\cdot\>\big)}
\def\suvi{{\rm SU}\big(V,\<\cdot,\cdot\>_i\big)}
\def\jm{J^{(m)}}
\def\val{\operatorname{Val}}
\def\vc{V_\CC}
\def\vconj{\overline v}
\def\wconj{\overline w}

\def\binfty{\mathcal B^\infty_{\mathrm {alt}}}
\def\hcb{{\rm H}_{\rm cb}}
\def\to{\rightarrow}
\def\bu{\bullet}
\def\la{L^\infty_{\mathrm{alt}}}
\def\hb{{\rm H}_{\rm b}}
\def\hc{{\rm H}_{\rm c}}
\def\h{{\rm H}}
\def\kxb{{\kappa_{\Xx}^{\rm b}}}
\def\kyb{{\kappa_{\Yy}^{\rm b}}}
\def\kx{\kappa_\Xx}
\def\sk{\mathsf k}
\renewcommand{\phi}{\varphi}

\def\eI{{^{I}\!\mathrm{E}}}
\def\eII{{^{I\!I}\!\mathrm{E}}}
\def\dI{{^{I}\!d}}
\def\dII{{^{I\!I}\!d}}
\def\fI{{^{I}\!F}}
\def\fII{{^{I\!I}\!F}}

\def\ind{\mathbf{i}}
\def\Ind{\mathrm{Ind}}
\def\No{N\raise4pt\hbox{\tiny o}\kern+.2em}
\def\no{n\raise4pt\hbox{\tiny o}\kern+.2em}
\def\bsl{\backslash}
%
%
%
%
%%%%%%%%%%%%%% Definitionen von Matrizen, etc. %%%%%%%%%%%%%%%%%%%%%%%

\newcommand{\mat}[4]{
\left(
\begin{array}{cc}
{\scriptstyle #1} & {\scriptstyle #2} \\ {\scriptstyle #3} & {\scriptstyle
#4}
\end{array}\right)}

\newcommand{\matt}[6]{
\left(
\begin{array}{cc}
{\scriptstyle #1} & {\scriptstyle #2} \\ {\scriptstyle #3} & {\scriptstyle
#4}\\{\scriptstyle #5} & {\scriptstyle #6}
\end{array}\right)}

\newcommand{\bet}[3]{{\beta(y_{#1},y_{#2},y_{#3})}}
\newcommand{\db}[4]{0=d\beta(y_{#1},y_{#2},y_{#3},y_{#4})
= \beta(y_{#1},y_{#2},y_{#3})+\beta(y_{#1},y_{#3},y_{#4})-\beta(y_{#1},y_{#2},y_{#4})
-\beta(y_{#2},y_{#3},y_{#4})}
%
%
%
%
%
%%%%%%%%%%%%%%%%%%% Titel, etc. %%%%%%%%%%%%%%%%%%ù

\title[Tight homomorphisms and Hermitian symmetric spaces]{Tight homomorphisms \\ and Hermitian symmetric spaces}
\author[M.~Burger]{Marc Burger}
\email{burger@math.ethz.ch}
\address{FIM, ETH Zentrum, R\"amistrasse 101, CH-8092 Z\"urich, Switzerland}
\author[A.~Iozzi]{Alessandra Iozzi}
\email{iozzi@math.ethz.ch}
\address{Department Mathematik, ETH Zentrum, R\"amistrasse 101,
  CH-8092 Z\"urich, Switzerland}
%{Institut f\"ur Mathematik, Universit\"at Basel, Rheinsprung 21,
%CH-4051 Basel, Switzerland}
%\address{D\'epartment de Math\'ematiques, Universit\'e de Strasbourg, 7, rue Ren\'e Descartes, F-67084 Strasbourg Cedex, France}
\author[A.~Wienhard]{Anna Wienhard}
\email{wienhard@math.princeton.edu}
\address{Department of Mathematics, Princeton University, Fine Hall, Washington Road, Princeton, NJ 08540, USA}
\thanks{A.I. and A.W. were partially supported by FNS grant
  PP002-102765. A.W. was partially supported by the National Science
  Foundation under agreement No. DMS-0111298 and No. DMS-0604665.}

%\keywords{ }
%\subjclass{ }

\date{\today}

\begin{abstract} 
We introduce the notion of tight homomorphism into a locally
compact group with nonvanishing bounded cohomology and study these  
homomorphisms in detail when the target is a Lie group of Hermitian type. 
Tight homomorphisms between Lie groups of Hermitian
type give rise to tight totally geodesic maps of Hermitian symmetric
spaces. We show that tight maps
behave in a functorial way with respect to the Shilov
boundary and use this to prove a general structure theorem for tight
homomorphisms. 
Furthermore we classify all tight embeddings of the Poincar\'e disk.
\end{abstract}
\maketitle
%
%
%
%
%
%%%%%%%%%%%%%%%%%%%%%%%%%%%%%%%%% Aufbau, Dateien %%%%%%%%%%%%%%%%%%%%%%%%%%%

%\input{titel}

\tableofcontents

%\vskip2cm
%\input{thanks}
%\input{tight_sympl}
\vskip1cm
\input{intro}

\vskip1cm
\input{tight_hom}

\vskip1cm
\input{shilov}
\vskip1cm
\input{structure_one}

\vskip1cm
\input{tubetype}

\vskip1cm
\input{fiveandahalf}
\vskip1cm
\input{structure_two}

\vskip1cm
\input{criterion}

\vskip1cm
\input{classification}

\appendix
\input{app}
\vskip1cm
\bibliographystyle{amsplain}
\bibliography{refs}
%\vskip1cm
%\input{abstract}
%\vskip1cm
%\input{leben}
\providecommand{\bysame}{\leavevmode\hbox to3em{\hrulefill}\thinspace}
\providecommand{\MR}{\relax\ifhmode\unskip\space\fi MR }
\providecommand{\MRhref}[2]{%
  \href{http://www.ams.org/mathscinet-getitem?mr=#1}{#2}
}
\providecommand{\href}[2]{#2}

\vskip1cm
\end{document}

%% file: intro.tex
\section{Introduction}
Let $L, G$ be locally compact second countable topological 
groups.
% and $E$ a Banach $G$-module. 
A continuous homomorphism $\rho: L\to G$ 
induces canonical pullback maps $\rho^*$ in continuous cohomology and 
$\rho_{\rm b}^\ast$ in continuous bounded cohomology. 
A special feature of continuous bounded cohomology is that 
it comes equipped with a canonical seminorm $\|\cdot\|$
with respect to which $\rho_{\rm b}^*$ is norm decreasing, that is 
\bqn
\big\|\rho_{\rm b}^*(\a)\big\| \leq \|\a\| \text{ for all } \a\in \hcb^\bullet(G,\RR)\,.
\eqn

Given a class $\alpha \in \hcb^\bullet(G,\RR)$ 
we say that a homomorphism $\rho: L\to G$ is {\em $\alpha$-tight} 
if the pullback $\rho_{\rm b}^*$ preserves the norm of $\a$, 
that is $\|\rho_{\rm b}^*(\a)\| = \|\a\|$.

%We give basic properties of general $\a$-tight homomorphisms in
%\S~2. These are immediate consequences of corresponding properties of
%continuous bounded cohomology established via the functorial
%approach in \cite{Burger_Monod_GAFA, Monod_book}. 

For the main part of the article we specialize to the situation
when the target group $G$ is of Hermitian type, i.e. $G$ is a 
connected semisimple Lie group with finite center and without compact factors
such that its associated symmetric space $\Xx$  is Hermitian symmetric. 
Let $J$ be the $G$-invariant complex structure on $\Xx$; combining it with
the unique $G$-invariant Riemannian metric of minimal holomorphic sectional 
curvature $-1$, gives rise to the K\"ahler form $\omega_{G} \in
\Omega^2(\Xx)^G$. We denote by $\kgb \in \hcb^2(G, \RR)$ 
the bounded continuous cohomology class obtained in the familiar way
(see \S~\ref{subsec:hss})
by integration of $\omega_{G}$ over triangles with geodesic sides. 
 
\begin{defi_intro}
Let $L$ be a locally compact second countable topological group and $G$ a group of Hermitian type. A continuous homomorphism
$\rho:L\to G$ is said to be {\em tight}
if $\big\|\rho_{\rm b}^*(\kgb)\big\|= \|\kgb\|$.
\end{defi_intro}

It is implicit in the definition of a tight homomorphism that it depends on the $G$-invariant complex structure $J$ which is part of the data of the Hermitian symmetric space $\Xx$.

\begin{fexo_intro}
Let $\Gamma<\mathrm{SU}(n,1)=:G$ be a cocompact lattice and 
$M:=\Gamma\backslash G$ the corresponding compact hyperbolic manifold.
The (ordinary) K\"ahler class $\rho^\ast(\kg)$ of a representation $\rho:\Gamma\to G$, 
seen as a de Rham class on $M$,
can be paired with the K\"ahler form $\omega_M$ on $M$ 
to give a characteristic number
\bqn
i_\rho:=\frac{\langle\rho^\ast(\kg),\omega_M\rangle}{\langle\omega_M,\omega_M\rangle}
\eqn
which satisfies a Milnor-Wood inequality \cite{Burger_Iozzi_difff}
\bqn
|i_\rho|\leq\r_\Xx\,.
\eqn
Representations such that $i_\rho = \r_\Xx$ are called {\it maximal}\footnote{A modification of the above construction leads to the definition 
of an analogous invariant even in the case of nonuniform lattices,
\cite[\S~5]{Burger_Iozzi_difff}.}. Maximal representations are tight  \cite[Lemma~5.3]{Burger_Iozzi_difff}, 
and in 
fact, they are the most important examples of such.\footnote{For surface groups  
one can easily construct tight
homomorphisms which are not maximal out of a maximal representation of the fundamental group of a lower genus surface. Note however, that there are tight homomorphism of surface groups which are far from being maximal, \cite{Breuillard_Gelander_Souto_Storm}.}
\end{fexo_intro}

%We define a representation $\rho$ to be {\it maximal} if $i_\rho=\r_\Xx$.
%A modification of the above construction leads to the definition 
%of an analogous invariant even in the case of nonuniform lattices,
%\cite[\S~5]{Burger_Iozzi_difff}.

%It is immediate from \cite[Lemma~5.3]{Burger_Iozzi_difff} that 
%maximal representations are tight.  
%In fact, they are most important examples of such 
%and in the case of compact surface groups, 
%are the ones which led us to the study of tight homomorphisms, 
%\cite{Burger_Iozzi_Wienhard_ann, Burger_Iozzi_Labourie_Wienhard, Burger_Iozzi_Wienhard_tol}. 

%It should be remarked however that although 
The study of the structure of tight homomorphisms
is paramount in the classification of maximal representations of compact surface groups
\cite{Burger_Iozzi_Wienhard_ann, Burger_Iozzi_Labourie_Wienhard, Burger_Iozzi_Wienhard_tol}. It should be remarked however that  
the scope of the notion goes well beyond this, 
as for example every surjection of a finitely generated group onto a lattice in $G$ is tight
(see Corollary~\ref{cor:latticetight}). 
% and even for surface groups there are tight
%homomorphisms which are not maximal, \cite{Breuillard_Gelander_Souto_Storm}.
In particular we have:

\begin{prop_intro}\label{cor:mapgroup}
Let $\mathbf{Mod}_g$ be the mapping class group of a closed surface of genus $g\geq 1$. 
Then the natural homomorphism $\mathbf{Mod}_g \to \mathrm{Sp}(2g,\RR)$ is tight. 
\end{prop_intro} 

%\begin{rem_intro} Kotschick shows that there are no maximal representations 
%of a compact surface group which can factor through $\mathbf{Mod}_g \to \mathrm{Sp}(2g,\RR)$,
%\cite{Kotschick_signatures}.
%\end{rem_intro}

%The most important examples of tight homomorphisms, and the ones which led us to
%the study of tight homomorphisms, are maximal representations of surface groups, 
%\cite{Burger_Iozzi_Wienhard_ann, Burger_Iozzi_Labourie_Wienhard, 
%Burger_Iozzi_Wienhard_tol}. 
%But tight homomorphisms are more general, 
%for example every surjection of a finitely generated group onto a lattice in $G$ is tight,
%or the maximal representations of a lattice $\Gamma<\mathrm{SU}(n,1)$
%into a group of Hermitian type, as considered in \cite{Burger_Iozzi_difff},
%are tight homomorphisms.  
%More explicitly, recall that if $\Gamma<\mathrm{SU}(n,1)=:G$ is a cocompact lattice and 
%$M:=\Gamma\backslash G$ is the corresponding compact hyperbolic manifold, 
%the (ordinary) K\"ahler class $\rho^\ast(\kg)$ of a representation $\rho:\Gamma\to G$ 
%seen as a de Rham class on $M$
%can be paired with the K\"ahler form $\omega_M$ on $M$ 
%to give a characteristic number
%\bqn
%i_\rho:=\frac{\langle\rho^\ast(\kg),\omega_M\rangle}{\langle\omega_M,\omega_M\rangle}
%\eqn
%which satisfies a Milnor-Wood inequality
%\bqn
%|i_\rho|\leq\r_\Xx\,.
%\eqn

\medskip
One of the main points of this paper is the following 
structure theorem for tight homomorphisms.

\begin{thm_intro}\label{thm:intro_one}
Let $L$ be a locally compact second countable group, $\gG$ a connected algebraic group 
defined over $\RR$ such
that $G:= \gG(\RR)^\circ$ is of Hermitian type.
Suppose that $\rho: L\to G$ is a continuous tight homomorphism.
Then: 
\begin{enumerate}
\item The Zariski closure $\hH := \ol{\rho(L)}^Z$ is reductive. 
\item The centralizer $\Zz_G(H)$ of $H:= \ol{\rho(L)}^Z(\RR)^\circ$ is
  compact.
\item The symmetric space $\Yy$ corresponding to $H$  is Hermitian and
  $\Yy$ admits a unique 
  $H$-invariant complex structure such that the inclusion $H\to G$ is
  tight and positive.
\end{enumerate}
\end{thm_intro}
To explain the notion of a positive homomorphism, 
let us recall that the complex structure $J$ defines 
a cone $\hc^2(G,\RR)^{\geq0}$ of {\em positive K\"ahler classes} and,
via the isomorphism
\bqn
\xymatrix@1{
 \hcb^2(G,\RR)\ar[r]^\cong
&\hc^2(G,\RR)\,,
}
\eqn
a cone of {\em bounded positive K\"ahler classes} containing in particular $\kgb$.
%The notion of tightness depends only on the $G$-invariant complex structure on $\Xx$, not 
%on the particular choice of a specific cohomology class $\a\in\hcb^2(G, \RR)^{>0}$. 
A continuous homomorphism $\rho: G_1 \to G_2$ between two groups of
Hermitian type is said to be {\em positive} if $\rho_{\rm b}^*\kgtb \in \hcb^2(G_1,\RR)^{\geq0}$.

\medskip

As an immediate application of Theorem~\ref{thm:intro_one}, we have:

\begin{cor_intro} Let $\gG$ a connected algebraic group 
defined over $\RR$ such that $G:= \gG(\RR)^\circ$ is of Hermitian type
and let $\rho:\Gamma\to G$ be a maximal representation
of a lattice $\Gamma<\mathrm{SU}(n,1)$.  Then\footnote{In the case in which $\Gamma$ is the fundamental group of an oriented compact surface possibly with boundary, 
one can reach much stronger conclusions, 
as for example faithfulness and discreteness of $\rho$, 
\cite{Burger_Iozzi_Wienhard_ann, Burger_Iozzi_Labourie_Wienhard, Burger_Iozzi_Wienhard_tol}.}:
\begin{enumerate}
\item The Zariski closure $\hH := \ol{\rho(\Gamma)}^Z$ is reductive. 
\item The centralizer $\Zz_G(H)$ of $H:= \ol{\rho(\Gamma)}^Z(\RR)^\circ$ is
  compact.
\item The symmetric space $\Yy$ corresponding to $H$  is Hermitian and
  $\Yy$ admits a unique 
  $H$-invariant complex structure such that the inclusion $H\to G$ is
  tight and positive.
\end{enumerate}
\end{cor_intro}

%
%Statement (1) that $\ol{\rho(L)}^Z$ is reductive implies
%by the deep result of Corlette \cite{Corlette_88} the following 
%
%\begin{cor_intro}\label{cor:harmonicmap_intro}
%Let $\G=\pi_1(M)$ be the fundamental group of a closed Riemannian
%manifold and let $\rho:\G\to G$ be a tight homomorphism. Then there exists a
%$\rho$-equivariant harmonic map $f:\widetilde{M} \to \Xx$. 
%\end{cor_intro}

Our study of tight homomorphisms relies on the study of a parallel notion of tightness 
for totally geodesic maps of symmetric spaces. 
Namely, let $\Xx_1, \Xx_2$ be Hermitian symmetric spaces of noncompact type.

\begin{defi_intro} A totally geodesic map $f: \Xx_1 \to \Xx_2$ is said to be {\em tight} if 
\bq\label{eq:sup}
\sup_{\Delta\subset \Xx_1} \int_\Delta f^*\omega_{G_2} =  \sup_{\Delta\subset \Xx_2}
\int_\Delta \omega_{G_2}\,,  
\eq
where the supremum is taken over all smooth oriented triangles $\Delta$ with geodesic sides in
$\Xx_1$ and $\Xx_2$, respectively. 
\end{defi_intro}

A first immediate fact is the following:

\begin{prop_intro}
Let $G_1$, $G_2$ be Lie groups of Hermitian type and let $\Xx_1, \Xx_2$ 
be the corresponding symmetric spaces. 
A homomorphism $\rho: G_1 \to G_2$ is tight if and only if 
the corresponding totally geodesic map $f: \Xx_1 \to \Xx_2$ is tight.
\end{prop_intro}

%Since the pullback of the K\"ahler form does not vanish,
%tight embeddings are never totally real but they are also not necessarily holomorphic. 
%Nevertheless 
Tight embeddings behave nicely at infinity. 
Recall that the Shilov boundary $\cs_\Xx$ of the Hermitian
symmetric space $\Xx$ is the unique closed $G$-orbit in the
(topological) compactification of the bounded symmetric domain
realization $\Dd$ of $\Xx$, and can be identified with $G/Q$, 
where $Q<G$ is an appropriate parabolic subgroup. 
Two points $x,y$ in $\cs_\Xx$ are said to be {\em transverse}
if $(x,y)$ is in the unique open $G$-orbit in $\cs_\Xx \times
\cs_\Xx$.

\begin{thm_intro}\label{thm:intro_2a}
Let $G_1, G_2$ be Lie groups of Hermitian type and $\cs_1$,
$\cs_2$ the Shilov boundaries of the associated symmetric spaces. 
Let $\rho: G_1\to G_2$ be a continuous tight homomorphism 
and $f:\Xx_1\to\Xx_2$ the corresponding totally geodesic tight map.
Then there exists a $\rho$-equivariant continuous map $\check{f}: \cs_1  \to \cs_2$
which extends $f$ and which maps transverse pairs to transverse pairs.
\end{thm_intro}

Using this theorem we can establish a general existence result for boundary maps. 
Let $\G$ be a countable discrete  group and $(B,\nu)$ a Poisson boundary for $\G$. 
Recall that under these conditions, 
the amenability of the $\Gamma$-action on $(B,\nu)$ 
insures the existence of a $\rho$-equivariant measurable map 
from $B$ to the space of probability measures on $G/P$, 
where $P$ is a minimal parabolic in $G$.  
Under some conditions, such as for instance Zariski density 
of the image of the representation $\rho$, one can deduce the existence
of such a map with values in $G/P$ (see \cite{Burger_Iozzi_supq, Burger_Iozzi_Wienhard_kahler}).

For tight homomorphisms we have the general existence result:

\begin{thm_intro}\label{thm:boundarymapintro}
Let $\gG$ be a semisimple algebraic group defined over $\RR$ 
such that $G:=\gG(\RR)$ is of Hermitian type and
let $\rho:\G\to G$ be a tight homomorphism. Then there exists a
$\rho$-equivariant measurable boundary map $\phi:B \to \cs_\Xx$. 
\end{thm_intro}

Recall that Hermitian symmetric spaces fall into two classes,
according to whether or not they admit a genuine generalization analogous to the
upper half plane model of the Poincar\'e disk. Namely, a Hermitian symmetric
space is of {\em tube type} if it is biholomorphically
equivalent to a domain $V\oplus i\Omega$ where $\Omega\subset V$ is a
proper open cone in the real vector space $V$. 
For any Hermitian symmetric space $\Xx$ maximal subdomains of tube type exist, 
they are of the same rank as $\Xx$, holomorphically embedded and pairwise conjugate. 

\begin{thm_intro}\label{thm:intro_three} 
Let $\rho: G_1 \to G_2$ be a tight homomorphism and 
$f: \Xx_1 \to \Xx_2$ the corresponding tight map. 
Then:
\begin{enumerate}
\item If $\Xx_1$ is of tube type, then there exists a unique maximal
  subdomain of tube type $T \subset \Xx_2$ 
such that $f(\Xx_1) \subset T$. Moreover $\rho(G_1)$ preserves $T$. 
\item If $\ker\rho$ is finite and $\Xx_2$ is of tube type, then $\Xx_1$ is of tube type.
\end{enumerate}
\end{thm_intro}

\medskip
Our main tool to study tight embeddings and tight homomorphisms are
{\em diagonal disks}. 
Recall that a {\em maximal polydisk} in $\Xx$ is the image of a holomorphic
and totally geodesic embedding $t: \DD^{\rx} \to \Xx$. Maximal
polydisks arise as complexifications of maximal flats in $\Xx$ and are conjugate. 
It is easy to check that maximal polydisks are tightly embedded. 
A diagonal disk in $\Xx$ is the image of the diagonal 
$\Delta(\DD) \subset  \DD^{\rx}$ under the embedding
$t: \DD^{\rx} \to \Xx$. Diagonal disks are precisely 
tight and holomorphically embedded disks in $\Xx$.

Using diagonal disks we can give 
a simple criterion for tightness of a totally geodesic
embedding $f:\Xx_1\to \Xx_2$ in terms of the corresponding
homomorphism of Lie algebras (see Lemma~\ref{lem:criterion}). 
We apply this criterion to  
classify tight embeddings of the Poincar\'e
disk and obtain:

\begin{thm_intro}\label{thm:hermitianhull_intro}
Suppose that $\Xx$ is a Hermitian symmetric space and
 $f:\DD\to \Xx$ is a tight embedding. 
Then the smallest Hermitian symmetric subspace $\Yy\subset \Xx$ 
containing $f(\DD)$ 
is a product $\Yy= \Pi_{i=1}^k \Yy_i$ of Hermitian symmetric subspaces $\Yy_i$ of $\Xx$, 
where $\Yy_i$ is the Hermitian symmetric space  
associated to the symplectic group ${\rm Sp}(2n_i,\RR)$. 
Moreover, $\sum_{i=1}^k n_i \leq r_\Xx$ and the embedding $f:\DD\to \Yy_i$ 
is equivariant with respect to the irreducible representation 
$\SL(2,\RR)\to \Sp(2n_i, \RR)$.
\end{thm_intro}

Tight embeddings are never totally real but they are also not necessarily holomorphic. 
The irreducible representations $\SL(2,\RR)\to \Sp(2n, \RR)$ provide examples of non-holomorphic tight embeddings of the Poincar\'e disk when $n\geq 2$.

 We are not aware of an example of a non-holomorphic tight embedding of an irreducible Hermitian symmetric space of rank $\r_\Xx \geq 2$. It might be that all tight homomorphisms of higher rank Hermitian symmetric spaces are holomorphic.

We suspect that tight embeddings of Hermitian symmetric spaces that are not of tube type are always holomorphic. 
For tight embeddings of $\CC\HH^n$ into classical Hermitian symmetric spaces of rank $2$ this can be deduced from \cite{Koziarz_Maubon_2}.

%\begin{thm_intro}\label{thm:complexhyp_intro}
%Let $\Xx$ be a Hermitian symmetric space and 
%$f:\HH_\CC^n \to \Xx$, $n\geq 2$, a tight embedding. Then 
%$f$ is holomorphic. 
%Furthermore $f$ factors through the diagonal
%embedding 
%$f:\HH_\CC^n \to (\HH_\CC^n)^{\r_\Xx}\to \Xx$.
%\end{thm_intro}
%
%In \S~2 we first give the basic properties of general $\a$-tight
%homomorphisms, then we specialize to tight homomorphism into Lie
%groups of type (RH) (reductive Hermitian), 
%which is a slight generalization of Lie groups of Hermitian type, and introduce tight embeddings.  
%In \S~3 we recall the definition of the Bergmann cocycle 
%and its relation with the Shilov boundary. 
%In \S~4 we prove the first structure theorem for tight embeddings. 
%In \S~5 we recall some facts about tube type domains 
%and prove the second structure theorem for tight embeddings. 
%In \S~6 we deduce the structure theorem for tight homomorphisms and state some of its consequences. 
%In \S~7 we establish a criterion for totally geodesic embeddings to be tight, 
%which allows us to give several examples and non-examples. 
%In \S~8 we classify tight embeddings of the Poincar\'e disk.

\medskip

{\it Acknowledgments:}
We thank Domingo Toledo for useful discussions about tight embeddings
of complex hyperbolic spaces. 
%%% Local Variables: 
%%% mode: latex
%%% TeX-master: "tight"
%%% End: 

%% file: tight_hom.tex
\section{Tight Homomorphisms}\label{sec:tight_hom}
\subsection{Continuous Bounded Cohomology}
In this section we recall some properties of bounded continuous cohomology which
are used in the sequel. For proofs and a comprehensive account of
continuous bounded cohomology the reader is referred to
\cite{Monod_book, Burger_Monod_GAFA}.

If $G$ is a locally compact second countable group, %and let $\RR$ be a
%trivial $G$-module. 
then 
\bqn
\cb(G^{k+1},\RR):=\big\{f:G^{k+1}\to\RR :\,f\text{ is continuous  and }\|f\|_\infty <\infty\big\}
\eqn
is a $G$-module via the action
\bqn
(hf)(g_0, \dots, g_k) = f(h^{-1} g_0, \dots , h^{-1} g_k)\,.
\eqn

The {\em continuous bounded cohomology} $\hcb^\bullet(G, \RR)$ of $G$ with
coefficients in $\RR$  
is the cohomology of the complex of $G$-invariants 
\bqn
\xymatrix{
0 \ar[r] & \cb(G,\RR)^G \ar[r]^-{d} & \cb(G^{2},\RR)^G \ar[r]^-{d} & \cdots
}
\eqn
where 
$d$ is the usual homogeneous coboundary operator defined for 
$f\in\cb(G^k,\RR)$ by 
\bqn
df(g_0, \dots g_{k}) :=
\sum_{i=0}^{k} f(g_0, \dots, \hat{g_i}, \dots\,,g_{k})\,.
\eqn

The supremum norm gives  $\cb(G^{\bullet+1}, \RR)^G$ the structure of a Banach space 
and induces a canonical seminorm $\|\cdot\|$ on $\hcb^\bullet(G, \RR)$
\bqn
\|\a\|= \inf_{[f]=\a} \|f\|_\infty\,.
\eqn

The inclusion of complexes $\cb(G^{\bullet+1},\RR) \subset \mathrm{C}(G^{\bullet+1},\RR)$,
where $\mathrm{C}(G^{\bullet+1},\RR)$ denotes the space of continuous real valued functions, 
induces a natural comparison map 
\bq\label{eq:comparison}
\c_G:\hcb^\bullet(G, \RR)
\to \hc^\bullet(G, \RR)
\eq
from continuous bounded cohomology to continuous cohomology.
Moreover, any continuous homomorphism  $\rho:L\to G$ of locally compact groups
induces canonical pullbacks both in continuous cohomology and in continuous bounded cohomology,
such that the diagram
\bqn
\xymatrix{
\hcb^\bullet(G, \RR) \ar[d]^{\c_G} \ar[r]^{\rho_{\rm b}^*} & \hcb^\bullet(L,\RR) \ar[d]^{\c_L}\\
\hc^\bullet(G,\RR) \ar[r]^{\rho^*} & \hc^\bullet(L,\RR)
}
\eqn
commutes. In particular, if $L<G$ is a closed subgroup, the pullback
given by the inclusion is the restriction map.

\begin{prop}\label{prop:bounded_properties}
\begin{enumerate}
\item Let $L$ be a locally compact second countable group and $L_0<L$ a closed subgroup.
If $L_0$ is of finite index in $L$, then the restriction map 
\bqn
\ba
\hcb^\bullet (L, \RR) &\longrightarrow\hcb^\bullet (L_0, \RR)\\
\kappa\qquad&\longmapsto\quad\kappa|_{L_0}
\ea
\eqn
is an isometric isomorphism\footnote{The statement is true more generally 
if there exists an $L$-invariant mean on $L/L_0$, but we shall not need this 
here.}, \cite[Proposition~2.4.2]{Burger_Monod_GAFA}.
\item If $R \lhd G$ is a closed amenable normal subgroup, the canonical projection 
$p:G\to G/R$ induces an isometric isomorphism via the pullback
\bqn
\xymatrix@1{
p_\mathrm{b}^\ast:\,\hcb^\bullet(G/R, \RR) \ar[r]& \hcb^\bullet(G, \RR)}
\eqn
in continuous bounded cohomology, \cite[Corollary~8.5.2]{Monod_book}.
\item  The seminorm $\|\,\cdot\,\|$ on $\hcb^2(G,\RR)$ is a norm which turns
it into a Banach space, \cite[Theorem 2]{Burger_Monod_GAFA}.
\item If $G = G_1 \times \cdots \times G_n$ is a direct
product of locally compact second countable groups, then the map
\bq\label{eq:isomsum}
\ba
\hcb^2(G,\RR) \to & \prod_{i=1}^n \hcb^2(G_i,\RR)\\
\k \qquad\mapsto&\quad (\k_{|_{G_i}})_{i=1}^n
\ea
\eq
into the Cartesian product of the continuous bounded cohomology of the factors,
is an isometric isomorphism\footnote{By a slight abuse of notation, 
we indicate by $\kappa|_{G_i}$ the pullback of $\kappa\in\hcb^2(G,\RR)$
by the homomorphism
\bqn
\ba
G_i&\to G=G_1\times\dots\times G_n\\
g_i&\longmapsto(e,\dots,g_i,\dots,e)\,.
\ea
\eqn}, \cite[Corollary 4.4.1.]{Burger_Monod_GAFA}, that is 
\bq\label{eq:norm}
\|\k\| = \sum_{i=1}^n \big\|\k_{|_{G_i}}\big\|\,.
\eq
%\item Let $G$ be a semisimple real Lie group with finite center. Then the
%comparison map 
%$\c_G: \hcb^2(G,\RR) \to \hc^2(G,\RR)$ is an isomorphism, and the dimension of
%  $\hc^2(G, \RR)$ equals the number of simple factors of $G$ which are
%  of Hermitian type, \cite{Burger_Monod_JEMS}. {\bf PUT IT LATER.}
\end{enumerate}
\end{prop}

\begin{rem}\label{rem:norm_decomp}
The fact that the isomorphism $\hcb^2(G,\RR) \cong \prod_{i=1}^n
\hcb^2(G_i,\RR)$ is isometric is not stated as such in 
\cite{Burger_Monod_GAFA}, but it follows from the proof. 
Moreover an explicit inverse to the map 
$\hcb^2(G,\RR) \to  \prod_{i=1}^n \hcb^2(G_i,\RR)$ 
in (\ref{eq:isomsum}) is given by 
\bqn
\ba
\prod_{i=1}^n \hcb^2(G_i,\RR) \to &\hcb^2(G,\RR)\\
(\k_i)_{i=1}^n \quad\mapsto& \sum_{i=1}^n (p_i)_\mathrm{b}^* \k_i\,, 
\ea
\eqn 
where $p_i: G \to G_i$ is the projection onto the $i$-th factor. 
\end{rem}

We record the following fact as a consequence of Proposition~\ref{prop:bounded_properties}.

\begin{cor}\label{cor:2.3} Let $L=H\cdot R$, 
where $L$ is a locally compact second countable group,
$H$ and $R$ are closed subgroups of $L$. 
We assume that $R$ is amenable and normal in $L$.  
Then the restriction map
\bqn
\hcb^\bullet(L,\RR)\to\hcb^\bullet(H,\RR)
\eqn
is an isometric isomorphism.
\end{cor}

\begin{proof} We have the following commutative diagram
\bqn
\xymatrix{
 L\ar[r]^p
&L/R\\
 H\ar@{^{(}->}[u]^i\ar[r]_-{p|_H}
&H/H\cap R\ar[u]_{\ol i}\,,
}
\eqn
where $\ol i$ is the topological isomorphism induced by $i$;
according to Proposition~\ref{prop:bounded_properties}(2),
$p_\mathrm{b}^\ast$ as well as $(p|_H)_\mathrm{b}^\ast$ and $\ol i_\mathrm{b}^\ast$
induce isometric isomorphisms.  This implies the assertion about $i_\mathrm{b}^\ast$.
\end{proof}

\subsection{General Facts about $\a$-Tight Homomorphisms}
The pullback $\rho_{\rm b}^*$ is seminorm decreasing with respect to the
canonical seminorm in continuous bounded cohomology, that is 
$\big\|\rho_{\rm b}^{*}(\a)\big\| \leq \|\a\|$ for all  $k\in \NN$ and all $\a\in \hcb^k(G,\RR)$. 
So, it is natural to give the following

\begin{defi}\label{def:tight_gen}
Let $L,G$ be 
locally compact second countable topological groups and $\a \in
\hcb^\bullet(G, \RR)$ a continuous bounded
cohomology class. 
A continuous homomorphism $\rho: L\to G$ is said to be {\em $\a$-tight}, if 
$\|\rho_{\rm b}^*(\a)\|=\|\a\|$.
\end{defi}

\begin{lemma}\label{lem:composition}
Let $L,G, H$ be locally compact second countable groups. Suppose that $\rho:L\to H$,
$\psi: H\to G$ are continuous homomorphisms. Let $\a \in \hcb^\bullet(G, \RR)$.
Then $\psi$ is $\a$-tight and $\rho$ is $\psi_{\rm b}^*(\a)$-tight 
if and only if $\psi \circ \rho$ is $\a$-tight.
\end{lemma}

\begin{proof}
The statement follows immediately from the chain of inequalities
\bqn
\|(\psi\circ\rho)_{\rm b}^* \a\| = \|\rho_{\rm b}^* \psi_{\rm b}^* \a\| \leq
\|\psi_{\rm b}^* \a\|\leq \|\a\|\,.
\eqn
\end{proof}
The following properties of tight homomorphisms are straightforward consequences
of the properties of continuous bounded cohomology summarized 
in Proposition~\ref{prop:bounded_properties}.

\begin{lemma}\label{lem:tight_properties} 
Let $L,G$ be locally compact second countable topological groups, 
$\a \in \hcb^\bullet(G,\RR)$ and $\rho:L\to G$ an $\a$-tight homomorphism.
\begin{enumerate}
\item Let $H<G$ be a closed subgroup. If the image $\rho(L)$ is contained in $H$ 
then $\rho$ is $\a_{|_H}$-tight and 
        $\|\a_{|_H}\|=\|\a\|$
\item Let $L_0<L$ be a closed subgroup of finite index\footnote{See the footnote in 
Proposition~\ref{prop:bounded_properties}(1).} in $L$. Then
        $\rho_{|_{L_0}}$ is $\a$-tight and 
\bqn
\big\|(\rho_{|_{L_0}})_{\rm b}^*\a\big\|= \|\rho_{\rm b}^*\a\|= \|\a\|\,.
\eqn
\item Let $R \lhd G$ be a closed amenable normal subgroup, $p: G \to G/R$ 
the canonical projection. Then the homomorphism 
        $p \circ\rho : L\to G/R$ is tight
        with respect to the class $(p_{\rm b}^*)^{-1} (\alpha)\in\hcb^2(G/R,\RR)$.
\item Let $\a \in \hcb^2(G,\RR)$ and if  $G:=G_1 \times \cdots G_n$, let  
$p_i: G \to G_i$ be the projection onto the $i$-th factor, $i=1,\dots,n$. 
Then $\rho_i = p_i \circ \rho: L \to G_i$ is $\alpha_{|_{G_i}}$-tight for all $i$.
\end{enumerate}
\end{lemma}

\begin{proof}
(1) Since $\rho(L)$ is contained in $H$ we have that 
  $\rho_{\rm b}^*\a=\rho_{\rm b}^*(\a_{|_H})$. If $\rho$ is $\a$-tight, then 
\bqn
\|\a\| = \|\rho_{\rm b}^*\a\| = \|\rho_{\rm b}^*(\a_{|_H})\| \leq
  \|\a_{|_H}\|\,.
\eqn 
Since $H<G$ is a subgroup we have that $\|\a_{|_H}\|\leq \|\a\|$ and the claim follows.\\

(2) Since $(\rho_{|_{L_0}})_{\rm b}^*\a$ is the restriction to $L_0$ of the class
$(\rho_{\rm b}^*\alpha)\in\hcb^\bullet(L,\RR)$, by Proposition~\ref{prop:bounded_properties}(1)
and tightness of $\rho$ we have that 
\bqn
\|({\rho_{|_{L_0}}})_{\rm b}^*\a \| = \|\rho_{\rm b}^*\a\| = \|\a\|\,.
\eqn

(3) The facts that $p_\mathrm{b}^\ast$ is an isometric isomorphism
(Proposition~\ref{prop:bounded_properties}(2)) and that $\rho$ is $\alpha$-tight
give rise to the following chain of equalities
\bqn
  \|(p\circ\rho)_{\rm b}^*(p_{\rm b}^*)^{-1}\a\| 
 =\|\rho_{\rm b}^* \a\|
 =\|\a\| 
 =\|(p_{\rm b}^*)^{-1}\a\|\,.
\eqn

(4) By Proposition~\ref{prop:bounded_properties}(1) and Remark~\ref{rem:norm_decomp}
we have that 
\bqn
\a=\sum_{i=1}^n(p_i)_\mathrm{b}^\ast(\a|_{G_i})\,.
\eqn
Then 
\bqn
\rho^\ast(\a)=\sum_{i=1}^n(\rho_i)_\mathrm{b}^\ast(\a|_{G_i})\,,
\eqn
so that, by (\ref{eq:norm}) and $\a$-tightness of $\rho$, we have 
\bqn
 \sum_{i=1}^n\|\a|_{G_i}\|
=\|\a\|
=\|\rho^\ast(\a)\|
=\|\sum_{i=1}^n (\rho_i)_{\rm b}^* (\a|_{G_i})\|
\leq\sum_{i=1}^n \|(\rho_i)_{\rm b}^* (\a|_{G_i})\|\,.
\eqn
The assertion now follows since
\bqn
\|(p_i)_\mathrm{b}^\ast(\a|_{G_i})\|\leq\|\a|_{G_i}\|\,.
\eqn
\end{proof}

%\subsection{Tight Homomorphisms}\label{subsec:tight_hom}

\subsection{Groups of Type (RH)}\label{subsec:hss}

Let $G$ be a connected reductive Lie groups with compact center. Then $G = G_c \cdot G_{nc}$, 
where $G_c$ is the largest compact connected normal subgroup of $G$ and $G_{nc}$ 
is the product of all connected noncompact almost simple factors of $G$. 
Then $G_{nc}$ has finite center and $G_c \cap G_{nc}$ is finite.

\begin{defi}\label{defi:RH}\be
\item A Lie group $G$  is of {\em Hermitian type} if it is connected semisimple 
with finite center and no compact factors and 
its associated symmetric space is Hermitian. 

\item A group $G$ is of {\em type (RH)} (reductive  Hermitian)  
if it is a connected reductive Lie group with compact center such that $G_{nc}$ is 
of Hermitian type. 
\ee
\end{defi}

If $G$ is a group of type (RH) and $\Xx$ is the symmetric space associated to $G_{nc}$ 
we have a homomorphism
\bq\label{eq:q}
\xymatrix@1{
 q:G \ar[r] 
&G_\Xx}
\eq
into $G_\Xx:=\Is(\Xx)^\circ$, which is surjective with compact kernel
so that $G$ acts properly on $\Xx$. 

Let $\Jj$ be the $G$-invariant complex structure on $\Xx$ and $g_\Xx$ the 
$G$-invariant Riemannian metric on $\Xx$, 
normalized so that the minimal holomorphic sectional curvature on every irreducible factor equals $-1$. 
We denote by $\omega_\Xx \in \Omega^2(\Xx)^{G}$ the $G$-invariant 
two-form
\bqn
\omega_\Xx (X,Y) := g_\Xx(X, \Jj Y)
\eqn
which is called the {\it K\"ahler form} of $\Xx$.

Choosing a base point $x_0 \in \Xx$ any $G$-invariant two-form
$\omega \in \Omega^2(\Xx)^{G}$ gives rise to a continuous cocycle 
\bq\label{eq:int}
\ba
c_\omega: G\times G \times G &\longrightarrow \qquad\quad\RR\\
(g_0, g_1, g_2) &\mapsto\frac{1}{2\pi} \int_{\Delta(g_0 x_0, g_1 x_0, g_2 x_0)}
\omega\,, 
\ea
\eq
where $\Delta(g_0 x_0, g_1 x_0, g_2 x_0)$ denotes a smooth oriented triangle with
geodesic sides and vertices $g_0 x_0, g_1 x_0, g_2 x_0$.
Let $\kappa_\omega = [c_\omega] \in \hc^2(G,\RR)$ denote 
the corresponding continuous cohomology class; then the maps
$\omega\mapsto\kappa_\omega$ implement the Van Est isomorphism \cite{Van_Est}
\bqn
\Omega^2(\Xx)^{G} \cong \hc^2(G,\RR)\,.
\eqn

It is well known that if $M$ is a connected simple Lie group with finite center, then
\bqn
\dim\hc^2(M,\RR)=0\text{ or } 1;
\eqn
in fact, the dimension is nonzero (hence 1) if and only if the associated
symmetric space $\Mm$ carries 
a $M$-invariant complex structure and hence is Hermitian symmetric.
Then $\Omega^2(\Mm)^M=\RR\,\omega_\Mm$, 
and with the above notations and normalizations we have:

\begin{thm}\cite{Domic_Toledo, Clerc_Orsted_2}\label{thm:dtco} If $\Mm$ is an irreducible
Hermitian symmetric space we have that
\bqn
 \frac{1}{2\pi}\left|\sup_{\Delta\subset\Mm}\int_\Delta\omega_\Mm\right|
=\frac{1}{2\pi}\sup_{\Delta\subset\Mm}\int_\Delta\omega_\Mm
=\frac{\r_\Mm}{2}\,,
\eqn
where $\r_\Mm$ denotes the rank of $\Mm$.
\end{thm}

In particular $c_{\omega_\Mm}$ defines a bounded class $\kappa_M^\mathrm{b}\in\hcb^2(M,\RR)$ 
which corresponds to $\kappa_\omega\in\hc^2(M,\RR)$ under the comparison map in (\ref{eq:comparison}),
and it was shown in \cite{Burger_Monod_JEMS} that the comparison map 
\bq\label{eq:comp2}
c_M:\hcb^2(M,\RR)\to\hc^2(M,\RR)
\eq
is an isomorphism in degree two.  
%In particular $c_{\omega_\Mm}$, and hence $c_\omega$ for any $\omega\in\Omega^2(\Mm)^M$,
%defines a bounded class $\kob\in\hcb^2(M,\RR)$ 
%which corresponds to $\kappa_\omega\in\hc^2(M,\RR)$ under the comparison map in (\ref{eq:comparison}),
The following result for the canonical norm in continuous bounded
cohomology could also in principle be deduced from \cite{Domic_Toledo, Clerc_Orsted_2}. 

\begin{thm}\label{Domic_Toledo, Clerc_Orsted_2}\label{thm:co-norm} With the above notations and if
$\Mm$ is irreducible, then
\bqn
\|\kappa_M^\mathrm{b}\|=\frac{\r_\Mm}{2}\,.
\eqn
\end{thm}

Strictly speaking the concept of bounded continuous classes and their norms does
not occur in \cite{Domic_Toledo, Clerc_Orsted_2}; what the authors show is that 
for a specific -- and hence any -- cocompact torsionfree lattice $\Gamma<M$, 
the singular bounded class in $\hb^2(\Gamma\backslash \Mm)$ defined by integration
of the K\"ahler form on straight simplices has Gromov norm $\pi\r_\Mm$. 
Using this and the isometric isomorphism between bounded singular cohomology 
of $\Gamma\backslash\Mm$ and bounded (group) cohomology of $\Gamma$,
one could deduce the above theorem.  We shall however give in the Appendix
a direct proof which in particular avoids the construction of lattices 
with specific properties in $M$.

\medskip
Let now $G$ be a group of type (RH), 
$\Xx=\Xx_1\times\dots\times\Xx_n$ a decomposition into irreducible factors,
and $\omega_{\Xx,i}:=p_i^\ast(\omega_{\Xx_i})$, where $p_i:\Xx\to\Xx_i$
is the projection onto the $i$-th factor.
Then 
\bq\label{eq:basis-omega}
\big\{\omega_{\Xx,i}\in\Omega^2(\Xx)^{G_\Xx}:\,1\leq i\leq n\big\}
\eq
gives a basis of $\Omega^2(\Xx)^{G_\Xx}$ and, in view of the Van Est isomorphism \cite{Van_Est},
\bq\label{eq:basis-h2}
\big\{\kappa_{\Xx,i}:=\kappa_{\omega_{\Xx,i}}\in\hc^2(G_\Xx,\RR):\,1\leq i\leq n\big\}
\eq
gives a basis of $\hc^2(G_\Xx,\RR)$. 
%changed $G$ to $G_\Xx$. 
Moreover, since it is the group $G_\Xx$ which acts effectively on $\Xx$,
it is obvious that $\Omega^2(\Xx)^G=\Omega^2(\Xx)^{G_\Xx}$,  and
hence the cohomology class defined by the cocycle $\c_\omega$ in (\ref{eq:int}) 
can be thought of as a cohomology class in $\hc^2(G,\RR)$.  
Hence the map $q$ in (\ref{eq:q}) defines an isomorphism
\bq\label{eq:isom-isom}
q^\ast:\hc^2(G_\Xx,\RR)\to\hc^2(G,\RR)
\eq
and we denote by
\bq\label{eq:abuse-not}
\big\{\kappa_{G,i}:=q^\ast(\kappa_{\omega_{\Xx,i}})\in\hc^2(G,\RR):\,1\leq i\leq n\big\}
\eq
the corresponding basis of $\hc^2(G,\RR)$.

If $\pi_i:G_\Xx\to G_{\Xx_i}$ denotes the projection onto the $i$-th factor, 
then we have that the analogous map to (\ref{eq:isom-isom}) in bounded cohomology 
\bq\label{eq:isom-isom1}
q_\mathrm{b}^\ast:\hcb^2(G_\Xx,\RR)\to\hcb^2(G,\RR)
\eq
and
\bq\label{eq:isom-isom2}
\ba
\prod_{i=1}^n\hcb^2(G_{\Xx_i},\RR)\to&\hcb^2(G_\Xx,\RR)\\
(\kappa_i)\qquad\longmapsto&\sum_{i=1}^n(p_i)_\mathrm{b}^\ast(\kappa_i)
\ea
\eq
are now isometric isomorphisms: for \eqref{eq:isom-isom1} this follows from
Proposition~\ref{prop:bounded_properties}(2) and the fact that $\ker q$ is compact,
and for \eqref{eq:isom-isom2} it follows already from Remark~\ref{rem:norm_decomp}.

Let $\kappa_{\Xx,i}^\mathrm{b}$ be the bounded class of $G_\Xx$
defined by $c_{\omega_{\Xx,i}}=c_{\omega_{\Xx_i}}\circ p_i$;
then it follows from the isomorphisms in (\ref{eq:comp2}) and (\ref{eq:isom-isom2})
that 
\bq\label{eq:basis-h2b}
\big\{\kappa_{\Xx,i}^\mathrm{b}\in\hcb^2(G_\Xx,\RR):\,1\leq i\leq n\big\}
\eq
gives a basis of $\hcb^2(G_\Xx,\RR)$ and, analogously to before,
\bq\label{eq:abuse-not-b}
\big\{\kappa_{G,i}^\mathrm{b}:=q_\mathrm{b}^\ast(\kappa_{\Xx,i}^\mathrm{b})\in\hcb^2(G,\RR):
\,1\leq i\leq n\big\}
\eq
a basis of $\hcb^2(G,\RR)$.

Thus if $\omega=\sum_{i=1}^n\lambda_i\omega_{\Xx,i}$ is any element in $\Omega^2(\Xx)^G$
written in the above basis (\ref{eq:basis-omega}), then
\bqn
\kappa_\omega^\mathrm{b}=\sum_{i=1}^n\lambda_i\kappa_{G,i}^\mathrm{b}
\eqn
is the bounded class in $\hcb^2(G,\RR)$ defined by
$c_\omega=\sum_{i=1}^n\lambda_i c_{\omega_{\Xx,i}}$
and corresponding to $\omega$ under the isomorphism $\Omega^2(\Xx)^G\to\hcb^2(G,\RR)$.
Moreover applying the isometric isomorphism in (\ref{eq:isom-isom2}) and Theorem~\ref{thm:co-norm}
we have that 
\bq\label{eq:marcela}
 \|\kappa_\omega^\mathrm{b}\|
=\sum_{i=1}^n|\lambda_i|\,\|\kappa_{G,i}^\mathrm{b}\|
=\sum_{i=1}^n|\lambda_i|\frac{\r_{\Xx_i}}{2}\,,
\eq
and in particular
\bq\label{eq:2.16}
\|\kgb\|=\frac{\r_\Xx}{2}\,.
\eq

With the same notations we have:

\begin{prop}\label{prop:norm_red} For any Hermitian symmetric space $\Xx$ we have that 
\bqn
     \frac{1}{2\pi} \sup_{\Delta\subset \Xx} \left|\int_\Delta\omega\right|
     = \frac{1}{2\pi} \sup_{\Delta\subset \Xx} \int_\Delta\omega
     = \sum_{i=1}^n |\lambda_i| \frac{\r_{\Xx_i}}{2}\,,
\eqn
where $\Delta \subset \Xx$ runs through all smooth triangles with geodesic sides in $\Xx$.
\end{prop}

\begin{proof}
We assume here Theorem~\ref{thm:dtco} and we focus on the nonirreducible case.
If $\omega = \sum_{i=1}^{n} \lambda_i \omega_{\Xx,i}$, then
\bqn
 \int_\Delta  \omega 
=\sum_{i=1}^{n} \lambda_i \int_\Delta \omega_{\Xx,i} 
=\sum_{i=1}^{n} \lambda_i \int_{p_i(\Delta)} \omega_{\Xx_i}
\eqn 
so that 
\bqn
\ba
     \left|\frac{1}{2\pi}\int_\Delta\omega\right|
    =\frac{1}{2\pi}\left|\sum_{i=1}^n\lambda_i\int_{p_i(\Delta)}\omega_{\Xx_i}\right|
 \leq\frac{1}{2\pi}\sum_{i=1}^n|\lambda_i|\,\left|\int_{p_i(\Delta)}\omega_{\Xx_i}\right|
 \leq\sum_{i=1}^n|\lambda_i|\,\frac{\r_{\Xx_i}}{2}\,,
\ea
\eqn
where we used Theorem~\ref{thm:dtco} in the last inequality.
For the opposite inequality, 
let $\eps >0$ and $\Delta_i \subset \Xx_i$ be smooth triangles with geodesic sides such that 
\bqn
\int_{\Delta_i} \omega_{\Xx_i} \geq \pi\r_{\Xx_i}  -\eps\,.
\eqn
More precisely let $\sigma_i^+: \Delta \to \Xx_i$ be a parametrization of $\Delta_i$ 
with geodesic sides, and $\sigma_i^-$ be 
the parametrization of $\Delta_i$ with the opposite orientation. 
Then let $\Delta_\Xx$ be the image of the map 
\bqn
\sigma: \Delta &\longrightarrow& \qquad\quad\Xx\\
t&\mapsto& \big(\sigma_1^{\eta_1}(t) , \cdots , \sigma_n^{\eta_n}(t)\big)\,, 
\eqn
where $\eta_i = \mathrm{sign}(\lambda_i)$. 
Then 
\bqn
\int_{\Delta_\Xx} \omega_{\Xx,i} = \int_{\sigma_i^{\eta_i}} \omega_{\Xx_i}
\eqn
and 
\bqn
     \int_{\Delta_\Xx} \omega 
   = \sum_{i=1}^{n} \lambda_i \int_{\sigma_i^{\eta_i}} \omega_{\Xx_i}
   = \sum_{i=1}^{n} |\lambda_i| \int_{\Delta_i} \omega_{\Xx_i}
\geq \pi\left( \sum_{i=1}^{n} |\lambda_i| \r_{\Xx_i}\right) - \eps\sum_{i=1}^{n} |\lambda_i|\,.
\eqn
Since this holds for any $\eps >0$ the proof is complete.
\end{proof}

\begin{defi}
Let $L$ be a locally compact second countable topological group and
$G$ a group of type (RH). 
A continuous homomorphism $\rho:L\to G$ is said to be {\em tight}, 
if $\rho$ is $\kgb$-tight, that is if 
\bqn
\|\rho^\ast(\kgb)\|=\|\kgb\|\,.
\eqn
\end{defi}

\begin{prop}\label{prop:tight_rank} Let $H,G$ be Lie groups of type (RH), $\Xx$  
the symmetric space associated to $G$ and $\Yy=\Yy_1\times\dots\times\Yy_n$ the symmetric space associated to $H$, 
where the $\Yy_i$'s are irreducible.  Let $\rho:H\to G$ be a homomorphism
and assume that $\rho^\ast(\kgb)=\sum_{i=1}^n\lambda_i\khib$
Then $\rho$ is tight if and only if $\rx = \sum_{i=1}^n |\lambda_i| \r_{\Yy_i}$.
\end{prop}
\begin{proof}
We have by \eqref{eq:marcela}
\bqn
 \|\rho^\ast(\kgb)\|
=\sum_{i=1}^n|\lambda_i|\,\|\khib\|
\eqn
and since $\|\kgb\| = \frac{\rx}{2}$ (by \eqref{eq:2.16})
and $\|\khib\| = \frac{\r_{\Yy_i}}{2}$ (by Theorem~\ref{thm:co-norm}), 
the assertion follows immediately.
\end{proof}

\subsection{Tight Maps}\label{subsec:tight_maps}
Let now $f:\Xx_1\to\Xx_2$ be a totally geodesic map,
where $\Xx_1,\Xx_2$ are Hermitian symmetric spaces of noncompact type..
Recall that this means that given any geodesic $c:\RR\to\Xx_1$,
the path $f\circ c:\RR\to\Xx_2$ is a geodesic, possibly not parametrized by arclength;
in fact, $f\circ c$ might be the constant map.
At any rate, for every triangle $\Delta_1\subset\Xx_1$
with geodesic sides, $f(\Delta_1)\subset\Xx_2$ is so as well, and hence
\bq\label{eq:trivial-ineq}
    \sup_{\Delta_1\subset\Xx_1}\int_{\Delta_1}f^\ast(\omega_{\Xx_2})
\leq\sup_{\Delta_2\subset\Xx_2}\int_{\Delta_2}\omega_{\Xx_2}
\eq
where the supremum in each side is taken over all triangles with geodesic sides.

\begin{defi}\label{defi:tight_emb}
\be
\item A totally geodesic map $f: \Xx_1 \to \Xx_2$ is {\em tight} if 
equality holds in (\ref{eq:trivial-ineq}).
\item We say that a subsymmetric space $\Yy\subset\Xx_2$ is {\em tightly embedded}
if the inclusion map is tight.
\ee
\end{defi}

\begin{prop} Let $f:\Xx_1\to\Xx_2$ be a totally geodesic map and $\Yy\subset f(\Xx_1)$
a subsymmetric space.
\be
\item The map $f$ is tight if and only if the subsymmetric space $f(\Xx_1)$
is tightly embedded in $\Xx_2$.
\item If $\Yy$ is tightly embedded in $\Xx_2$, then $f$ is tight.
\ee
\end{prop}
\begin{proof} For the first assertion it suffices to observe that every triangle
in $f(\Xx_1)$ with geodesic sides is the image of a triangle in $\Xx_1$ with geodesic sides.
This implies that
\bqn
  \sup_{\Delta_1\subset\Xx_1}\int_{\Delta_1}f^\ast(\omega_{\Xx_2})
 =\sup_{\Delta_2\subset f(\Xx_1)}\int_{\Delta_2}\omega_{\Xx_2}
\eqn
which shows the first assertion.

The second assertion follows immediately from the above and the inequalities
\bqn
    \sup_{\Delta\subset\Yy}\int_{\Delta}\omega_{\Xx_2}
\leq\sup_{\Delta_1\subset f(\Xx_1)}\int_{\Delta_1}\omega_{\Xx_2}
\leq\sup_{\Delta_2\subset\Xx_2}\int_{\Delta_2}\omega_{\Xx_2}\,.
\eqn
\end{proof} 
%
%\begin{cor}\label{cor:composition}
%Let $f: \Xx_1 \to \Xx_2$ be a totally geodesic map. If the image
%$f(\Xx_1)$ contains a symmetric subspace $\Yy\subset \Xx_2$ which is
%tightly embedded in $\Xx_2$, then $f$ is a tight embedding.
%\end{cor}
%
%\begin{proof} By hypothesis we have that 
%\bqn
%    \sup_{\Delta\subset\Yy}\int_{\Delta}\omega_{\Xx_2}
%   =\sup_{\Delta_2\subset\Xx_2}\int_{\Delta_2}\omega_{\Xx_2}
%\geq\sup_{\Delta_1\subset\Xx_1}\int_{\Delta_1}f^\ast\omega_{\Xx_2}
%   =\sup_{\Delta_1\subset\Xx_1}\int_{f(\Delta_1)}\omega_{\Xx_2}\,.
%\eqn
%The assertion now follows from the fact that any geodesic simplex in $\Yy$
%is in the image via $f$ of a geodesic simplex in $\Xx_1$.
%\end{proof}

Given a homomorphism $\rho: G_1 \to G_2$ of Lie groups of type (RH), 
let $x_1 \in \Xx_1$ be a base point, 
$K_1 = \stab_{G_1}(x_1)$ the corresponding maximal compact subgroup and 
$x_2 \in \Xx_2$ a point such that $\rho(K_1)(x_2) = x_2$. 
Then $\rho$ gives rise to a map 
\bqn
f: \Xx_1 \to \Xx_2
\eqn
defined by $f(g x_1):=\rho(g) x_2$, which is $\rho$-equivariant and totally geodesic.
 
\begin{lemma}\label{lem:diagram}
The diagram 
\bqn
\xymatrix{
\hcb^2(G_2,\RR) \ar[dd] \ar[rrr]^{\rho^*_{\rm b}} &&& \hcb^2(G_1,\RR) \ar[dd]\\
&\Omega^2(\Xx_2)^{G_2}\ar[r]^{f^*} \ar[ul]_{\cong} \ar[dl]_{\cong} & \Omega^2(\Xx_1)^{G_1} \ar[ur]^{\cong} \ar[dr]^{\cong}&\\
\hc^2(G_2, \RR)\ar[rrr]^{\rho^*} &&& \hc^2(G_1, \RR)
}
\eqn
commutes.
\end{lemma}

\begin{proof}
Let $\omega \in \Omega^2(\Xx_2)^{G_2}$. Using the points $x_1 \in \Xx_1$ 
and $x_2=f(x_1) \in \Xx_2$ in the construction of the cocycles we have
\bqn
c_{f^*(\omega)} = c_\omega \circ \rho\,.
\eqn
\end{proof}

\begin{cor}\label{cor:2.16}
In the above situation the following are equivalent
\begin{enumerate}
\item The homomorphism $\rho: G_1 \to G_2$ is tight.
\item The totally geodesic map $f: \Xx_1 \to \Xx_2$ is tight.
\end{enumerate}
\end{cor}
\begin{proof}
By Lemma~\ref{lem:diagram} we have that 
$\rho^\ast(\kappa^\mathrm{b}_{G_2})=\kappa^\mathrm{b}_{f^\ast(\omega_{\Xx_2})}$.
Thus, applying Proposition~\ref{prop:norm_red}
we have that 
\bqn
 \|\rho^\ast(\kappa^\mathrm{b}_{G_2})\|
=\|\kappa^\mathrm{b}_{f^\ast(\omega_{\Xx_2})}\|
=\frac{1}{2\pi}\sup_{\Delta_1\subset\Xx_1}\int_{\Delta_1}f^\ast\omega_{\Xx_2}
\eqn
and
\bqn
 \|\kappa^\mathrm{b}_{G_2}\|
=\frac{1}{2\pi}\sup_{\Delta_2\subset\Xx_2}\int_{\Delta_2}\omega_{\Xx_2}
\eqn
from which the equivalence follows readily.
\end{proof}

From Corollary~\ref{cor:2.16} and Proposition~\ref{prop:tight_rank} we immediately deduce
the following:

\begin{cor}\label{cor:tight_rank}
Let $H, G$ be Lie groups of type (RH)
with associated Hermitian symmetric spaces $\Yy$ and $\Xx$, 
$\rho: H\to G$ a homomorphism and $f:\Yy\to\Xx$ 
a $\rho$-equivariant totally geodesic map.
Let $\Yy= \Yy_1 \times \cdots \times \Yy_n$ be the decomposition 
into irreducible factors and suppose that 
$f^*(\omega_\Xx) = \sum_{i=1}^n \lambda_i \omega_{\Yy,i}$.
The totally geodesic map $f:\Yy\to\Xx$ is tight  
if and only if $\rx = \sum_{i=1}^n |\lambda_i| \r_{\Yy_i}$.
\end{cor}

Before stating the next corollary let us recall the following 

\begin{defi}\label{defi:maxi-polydisk}
\be
\item A {\em maximal polydisk} in $\Xx$ is the image of a totally geodesic and holomorphic 
embedding $t:\DD^{\r_\Xx} \to \Xx$  of a product of $r_\Xx$ Poincar\'e disks.
\item A {\em diagonal disk} in $\Xx$ is the image of the diagonal 
$\Delta(\DD) \subset \DD^{r_\Xx}$ under an embedding 
$t:\DD^{\r_\Xx} \to \Xx$ of $\DD^{\r_\Xx}$ as a maximal polydisk. 
In particular $d:= t\circ \Delta: \DD \to \Xx$ is a totally geodesic and holomorphic 
embedding. 
\ee
\end{defi}

Maximal polydisks arise as complexifications of maximal flats in $\Xx$, 
and hence are conjugate under $G_\Xx$. Moreover, with the normalization chosen 
in \S~\ref{subsec:hss} the embedding $t:\DD^{r_\Xx} \to \Xx$ is isometric.
In fact one can say more, as we have:

\begin{lemma}\label{lem:2.19} A metric on an Hermitian symmetric space $\Xx$ is normalized
if and only if every maximal polydisk $f:\DD^\rx\to\Xx$ is isometrically 
embedded.
\end{lemma}

\begin{proof} If $\Xx$ is irreducible, it follows from the computation 
in \cite[p.~273-274]{Clerc_Orsted_2}, that the holomorphic sectional curvature
is minimal at $u\in T_x\Xx$ if and only if the complex geodesic
obtained by $u$ is the image of a factor of a maximal polydisk
$t:\DD^\rx\to\Xx$. The general case for $\Xx$ not necessarily irreducible
follows immediately.
\end{proof}

Lemma~\ref{lem:2.19} has the following useful consequence:

\begin{cor}\label{cor:2.20} Let $\Yy\subset\Xx$ be a Hermitian subsymmetric space 
of the same rank as $\Xx$.
Then the restriction to $\Yy$ of the normalized metric on $\Xx$ is the normalized metric on $\Yy$.
\end{cor}

\begin{proof} Indeed, every maximal polydisk $t:\DD^{\r_\Yy}\to\Yy$ in $\Yy$ 
is a maximal polydisk in $\Xx$.
\end{proof}

\begin{defi}
A Hermitian symmetric space $\Xx$ is said to be of {\em tube type} 
if $\Xx$ is biholomorphically equivalent to a tube domain of the form  
\bqn
\{ v+ iu \, |\, v\in V \, , \, u \in \Omega\} \subset V\oplus iV\,, 
\eqn
where $V$ is a real vector space and $\Omega\subset V$ is a proper open cone. 
\end{defi}

Every Hermitian symmetric space $\Xx$ 
contains maximal (with respect to the inclusion) subdomains $T$ of
tube type (equal to $\Xx$ itself if $\Xx$ is of tube type)  
which are of the same rank as $\Xx$, holomorphically embedded and conjugate under $G_\Xx$. 
Moreover, the embedding $T\subset\Xx$ is always isometric (see Corollary~\ref{cor:2.20}).  

\begin{cor}\label{cor:first_examples}
\begin{enumerate}
\item Let $f: \Yy \to \Xx$ be a holomorphic and isometric embedding. 
Then $f$ is tight if and only if $\r_\Yy = \rx$.
\item Maximal polydisks $t: \DD^{\rx} \to \Xx$ are tight.
\item Diagonal disks  $d: \DD \to \Xx$ are tight.
\item Maximal tube type subdomains $T \subset \Xx$ are tight.
\end{enumerate}
\end{cor}
\begin{proof}
If $f: \Yy \to \Xx$ is holomorphic and isometric, 
then $f^*(\omega_\Xx) = \omega_\Yy = \sum_{i=1}^n \omega_{\Yy,i}$, 
so (1) follows from the fact that $\r_\Yy = \sum_{i=1}^n \r_{\Yy_i}$ 
and Proposition~\ref{prop:tight_rank}.

Then (2) and (4) follow at once from (1) since the embeddings
$t: \DD^{\rx} \to \Xx$  and $T\subset\Xx$ are holomorphic and isometric.

To see (3), observe that since $t$ is a holomorphic isometry, 
then $t^\ast(\omega_\Xx)=\omega_{\DD^{\r_\Xx}}$;
moreover an easy verification shows that 
$\Delta^\ast(\omega_{\DD^{\r_\Xx}})=\rx\omega_\DD$.
It follows then that $d^\ast(\omega_\Xx)=\rx\omega_\DD$,
so that the assertion follows from Proposition~\ref{prop:tight_rank}.
\end{proof}

Further examples of tight maps and tight homomorphisms
will be discussed in \S~\ref{sec:examples}.

%%% Local Variables: 
%%% mode: latex
%%% TeX-master: "tight"
%%% End: 

%% file: shilov.tex
\section{K\"ahler Classes and the Shilov boundary}\label{sec:shilov}
In this section we collect the facts from the geometry of Hermitian symmetric spaces,
some of which are of independent interest, needed for our purpose.  
Those concerning the geometry of triangles are due,
in the context of irreducible domains, to Clerc and \O rsted \cite{Clerc_Orsted_2};
we present also here the necessary -- easy -- extensions to general domains.

\medskip
Let $\Xx$ be an Hermitian symmetric space of noncompact type with a fixed 
$G_\Xx$-invariant complex structure $\Jj$. 
Fix a maximal compact subgroup $K = \Stab_{G_\Xx}(x_0)$, where $x_0 \in \Xx$ is some base point. 
Let $\frakg = \frakk \oplus \frakp$ be the corresponding Cartan decomposition, 
where $\frakg = {\rm Lie}(G_\Xx)$ and $\frakk = {\rm Lie}(K)$.
There exists a unique element $Z_\Jj$ in the center $\Zz(\frakk)$ of $\frakk$ such that 
$\ad(Z_\Jj)|_{\frakp}$ induces the complex structure $\Jj$ 
under the identification $\frakp \cong T_{x_0}\Xx$. 
The complexification $\frakg_\CC$ of $\frakg$ splits into eigenspaces of $\ad(Z_\Jj)$ as 
\bq\label{eq:gc_dec}
\frakg_\CC = \frakk_\CC \oplus \frakp_+ \oplus \frakp_-\,.
\eq

The Hermitian symmetric space $\Xx$ can be realized as a bounded
symmetric domain 
\bqn 
\Dd \subset \frakp_+ \cong \CC^N\,. 
\eqn 
Let us describe the structure of $\Dd$ more explicitly.  Let us fix
$\frakh \subset \frakk$ a maximal Abelian subalgebra.  Then
$\Zz(\frakk) \subset \frakh$ and $\frakh_\CC$ is a maximal Abelian
subalgebra, indeed a Cartan subalgebra of $\frakg_\CC$. The set $\Psi
= \Psi(\frakg_\CC, \frakh_\CC)$ of roots of $\frakh_\CC$ in
$\frakg_\CC$ decomposes corresponding to the decomposition of
$\frakg_\CC$ in (\ref{eq:gc_dec}) as 
\bqn 
 \Psi
=\Psi^{\frakk_\CC} \cup \Psi^{\frakp_+} \cup \Psi^{\frakp_-}\,, 
\eqn 
here $\Psi^{\frakl} :=
\left\{ \alpha \in \Psi\, |\, \text{ the root space }\frakg_\alpha
  \subset \frakl\right\}$. One can choose an ordering $\Psi = \Psi_+
\cup \Psi_-$ such that $\Psi^{\frakp_\pm} \subset \Psi_{\pm}$.

To every root $\alpha \in \Psi$ we associate a three dimensional
simple subalgebra 
\bq\label{eq:3.2} 
 \frakg_{[\alpha]} 
=\CC H_\alpha \oplus \CC E_\alpha \oplus \CC E_{-\alpha}\,,
\eq
where $H_\alpha \in \frakh_\CC$ is the unique element determined by
$\alpha(H) =2 \frac{\frakB(H,H_\alpha)}{\frakB(H_\alpha, H_\alpha)}$ 
for all $H\in \frakh_\CC$ and $\frakB$ is the Killing form on $\frakg_\CC$. 
The elements $E_\alpha, E_{-\alpha}$ are the elements of $\frakg_{\pm \alpha}$ 
satisfying the relation $[E_\alpha, E_{-\alpha}]= H_\alpha$ and 
$\tau(E_\alpha) = -E_{-\alpha}$, where $\tau$ is the complex conjugation of $\frakg_\CC$ 
with respect to the compact real form $g_U = \frakk \oplus i\frakp$.  
Then 
$\frakp_+=\sum_{\alpha \in \Psi^{\frakp_+}} \CC E_\alpha$ 
and the vectors 
$X_\alpha = E_\alpha+ E_{-\alpha}$, $Y_\alpha = 
i(E_\alpha - E_{-\alpha})$, $\alpha\in\Psi^{\frakp_+}$, 
form a real basis of $\frakp$.

Two roots $\alpha, \beta \in \Psi$ are called strongly orthogonal 
if neither $\alpha + \beta$ nor $\alpha - \beta$ is a root. 
By a theorem of Harish-Chandra 
there exists a maximal set $\Lambda = \{ \g_1, \cdots, \g_r\} \subset \Psi^{\frakp_+}$ 
of $r = \r_\Dd$ strongly orthogonal roots. 
The associated vectors $X_{\g_j} \in \frakp$ span 
a maximal Abelian subspace $\fraka$ of $\frakp$ over $\RR$. 
The bounded symmetric domain $\Dd$ admits the following description 
\bq\label{eq:bdd}
        \Dd
      = \left\{ \Ad(k) \sum_{j=1}^r \tanh(t_j) E_{\g_j} \, :\, k \in K, t_j \in \RR\right\} 
\subset \frakp_+\,;
\eq
moreover we call  
%\bq\label{eq:poly}
%        \Pp_0
%      = \left\{ \Ad(k) \sum_{j=1}^r \tanh(t_j) E_{\g_j} \, :\, k \in \Zz(K), t_j \in \RR\right\} 
%\subset \Dd
%\eq
\bq\label{eq:poly}
        \Pp_0
      = \left\{ \Ad(k) \sum_{j=1}^r \tanh(t_j) E_{\g_j} \, :\, k \in \exp(\frakh), t_j \in \RR\right\} 
\subset \Dd
\eq
the {\em standard maximal polydisk} and

\bq\label{eq:diag}
        \Delta_0
      = \left\{ \Ad(k) \sum_{j=1}^r \tanh(t) E_{\g_j} \, :\, k \in \Zz(K), t \in \RR\right\} 
\subset \Pp_0
\eq
the {\em standard diagonal disk}.  
With the explicit description of $\Dd$ we define the (normalized) Bergmann kernel 
\bqn
k_\Dd: \Dd \times \Dd \to \CC^\times\,, 
\eqn
by 
\bq\label{eq:bergmann}
k_\Dd(z,w) = h_\Dd(z,w)^{-2}\,,
\eq 
where $h_\Dd(z,w)$ is the polarization 
of the unique $K$-invariant polynomial $h$ on $\frakp_+$ such that 
\bqn
h\left(\sum_{j=1}^r s_j E_{\g_j} \right) = \prod_{j=1}^r (1-s_j^2)\,.
\eqn

The (normalized) Bergmann kernel is continuous on $\Dd^2$ 
and gives rise to a Riemannian metric $g_\Dd$, 
called the (normalized) Bergmann metric on $\Dd$, 
which has minimal holomorphic sectional curvature $-1$:
this holds in the irreducible case by \cite[(1.2)]{Clerc_Orsted_2}
and follows for the general case by the naturality under product
of the normalized metric. 
Let us observe that the Bergmann metric and the normalized Bergmann metric 
are equivalent: indeed for an irreducible domain they are proportional -- 
the proportionality factor however depends on the domain, 
see \cite[(1.2)]{Clerc_Orsted_2} for the precise value --
and the Bergmann and normalized Bergmann metrics behave functorially 
with respect to taking products.

The  K\"ahler form  given by 
\bqn
\omega_\Dd = i \partial \ol{\partial} \log k_\Dd(z,z)
\eqn 
corresponds to $\omega_\Xx$ under the isomorphism $\Xx\to\Dd$.

\begin{lemma}\label{lem:distance}
Let $\Dd \subset \CC^N$ be a bounded symmetric domain with Riemannian distance 
$d_\Dd(\cdot, \cdot)$. 
Then there exists a constant $c = c(\Dd)$ such that for all $x,y \in \Dd$
\bqn
d_{\Dd}(x,y) \geq c \|x-y\|_{eucl}\,, 
\eqn
where $\|\cdot\|_{eucl}$ denotes the Euclidean norm on $\CC^N$.
\end{lemma}
\begin{proof} Using the observation above,
it suffices to show the lemma for the distance coming from the Bergmann metric.

Let $b_\Dd$ be this metric, then at every $z\in\Dd$ we have
\bqn
(b_\Dd)_z(\,\cdot\,,\,\cdot\,)=\big\<\,\cdot\,,\Ad \Kk(z,z)^{-1}\,\cdot\,\big\>\,,
\eqn
where $\<\,\cdot\,,\,\cdot\,\>$ is the Hermitian form on $\frakp_+$ coming
from the Killing form and $\Kk$ is the kernel function defined on an open
subset of $\frakp_+\times\frakp_+$ with values in the complexification $K_\CC$
of the maximal compact subgroup
(for definition and details see \cite[\S~5.6 and Proposition~6.2]{Satake_book}).

We need now to estimate the eigenvalues of $\Ad \Kk(z,z)$.
Writing $z=\Ad(k)z_1$ and observing that 
\bqn
\Ad \Kk(z,z)=\Ad(k)\,\Ad \Kk(z_1,z_1)\Ad(k)^{-1}\,,
\eqn
we may assume that 
\bqn
z_1=\sum_{j=1}^{r_\Dd}\xi_j E_{\gamma_j}\,.
\eqn

An explicit calculation (see e.g. \cite[page 71]{Satake_book}) 
shows that the eigenvalues of $\Ad \Kk(z_1, z_1)$ on the 
root space $\frakg_\alpha$ for $\alpha\in\Psi^{\frakp_+}$ are given by 
\bqn
\ba
&(1-|\xi_j|^2)^2,\\
&(1-|\xi_j|^2)(1-|\xi_k|^2), \\
&(1-|\xi_j|^2), \text{ or } \\
&(1-|\xi_j|^2)(1+|\xi_k|^2) + |\xi_k|^4 \,, 
\ea
\eqn
where $1\leq j\neq k\leq \r_\Xx$.
In particular, since $0\leq |\xi_j| <1$ every eigenvalue of $\Ad \Kk(z,z)^{-1}$ is greater 
than $\frac{1}{3}$ and the claim follows.
\end{proof}
%\begin{rem}
%The reader might find the proof of Lemma~\ref{lem:distance} a bit technical. The problem is that one cannot easily control how the 
%Euclidean metric changes under the action of the automorphism group of $\Dd$.
%\end{rem}

%
\subsection{Shilov Boundary}\label{subsec:shilov}
%The Shilov boundary $\cs_\Dd$ of a bounded symmetric domain in $\CC^N$ is 
%the unique minimal closed subset of the topological boundary 
%$\partial \Dd$ of $\Dd$ with the property that all functions $f$, 
%which are continuous on $\ol{\Dd}$ and holomorphic on $\Dd$ satisfy
%$|f(x)| \leq \max_{y\in \cs_\Dd} |f(y)|$ for all $x\in \Dd$.
We shall denote by $G_\Dd$ the connected component of the group $\aut(\Dd)$
of holomorphic automorphisms of $\Dd$.  When $\Dd'\subset\Dd$ is a Hermitian
symmetric subspace, we shall denote by $A_{\Dd'}$ the subgroup of $G_\Dd$
of Hermitian type associated to $\Dd'$; in fact, $A_{\Dd'}$ is the product
of the noncompact connected almost simple factors of the reductive subgroup
$\Nn_{G_\Dd}(\Dd')$, where if $E\subset\ol\Dd$ is any subset, we define
\bqn
\Nn_{G_\Dd}(E):=\big\{g\in G_\Dd:\,g(E)=E\big\}\,.
\eqn

The closure $\ol{\Dd}$ contains a unique closed $G_\Dd$-orbit 
which is the Shilov boundary $\cs_\Dd$ of $\Dd$;
more precisely, the Shilov boundary $\cs_\Dd$ is the $G_\Dd$ orbit 
of the point $\sum_{j=1}^r  E_{\g_j} \subset \frakp_+$,
where $E_{\g_j}$ are the root vectors associated to strongly orthogonal roots 
$\gamma_j\in\Lambda$
(see \eqref{eq:3.2}),
and can hence be realized as $G_\Dd/Q$, 
where $Q$ is the stabilizer in $G_\Dd$ of $\sum_{j=1}^r  E_{\g_j}$.
In particular, if $\Dd$ is irreducible, 
then $Q$ is a maximal parabolic subgroup in $G_\Dd$. 
%Namely 
%let $\Dd = \Dd_1 \times \cdots \times \Dd_n$ be a decomposition into irreducible factors. Then the 
%Shilov boundary $\cs_\Dd = \cs_{\Dd_1} \times \cdots \times \cs_{\Dd_n}$. 

\begin{lemma}\label{lem:3.2}
\begin{enumerate}
\item Let $\Dd = \Dd_1 \times \cdots \times \Dd_n$ be a decomposition into irreducible factors. 
Then the Shilov boundary $\cs_\Dd$ of $\Dd$ is the product 
$\cs_{\Dd_1} \times \cdots \times \cs_{\Dd_n}$
of the Shilov boundaries of the irreducible factors.
\item If $\Pp$ is any maximal polydisk and $\Delta\subset\Pp$ is any diagonal disk,  
then $\cs_{\Delta} \subset \cs_\Pp \subset \cs_{\Dd}$.
\item If $\Dd' \subset \Dd$ is a Hermitian symmetric subspace with $\r_{\Dd'} = \r_{\Dd}$,
then $\cs_{\Dd'} \subset \cs_\Dd$.
\end{enumerate}
\end{lemma}
\begin{proof} (1) This first assertion follows from the characterization of $\cs_\Dd$ 
as the unique closed $G_\Dd$-orbit in $\ol\Dd$.

(2) To see the second assertion observe that since all maximal polydisks
(and their diagonal disks) are conjugate by $G_\Dd$,
it suffices to show the assertion for the standard maximal polydisk $\Pp_0$.
First it is obvious that $\cs_{\Delta_0}\subset\cs_{\Pp_0}$. 
Then let $A_{\Pp_0}$ be the subgroup of Hermitian type of $G_\Dd$
associated to $\Pp_0$; clearly the vector $\sum_{j=1}^r  E_{\g_j}$ is contained in $\cs_{\Pp_0}$
and hence its $A_{\Pp_0}$-orbit is contained in its $G_\Dd$-orbit, 
which implies that $\cs_{\Pp_0}\subset\cs_{\Dd}$ and hence the second assertion.

(3) Finally, let $\Pp\subset\Dd'$ be a maximal polydisk and 
let $A_{\Dd'}$ be the subgroup of Hermitian type of $G_\Dd$ associated to $\Dd'$.
Then $\Pp$ is maximal in $\Dd$ as well 
and hence, by (2), $\cs_{\Pp}\subset\cs_{\Dd'}$.
This, together with the obvious inclusion 
$A_{\Dd'}\big(\cs_{\Pp}\big)\subset G_{\Dd}\big(\cs_{\Pp}\big)$
implies that $\cs_{\Dd'}\subset\cs_\Dd$.
\end{proof}
%\begin{proof}
%Choose a maximal set $\Lambda_i$ of strongly orthogonal roots for every $\Dd_i$,
%then the union $\Lambda = \bigcup_{i=1}^n \Lambda_i$ is a maximal set
%of strongly orthogonal roots for $\Dd$. The Shilov boundary
%$\cs_{\Dd_i}$ is the $G_{\Dd_i}$ orbit of $\sum_{j=1}^{\r_i}
%E_{\g^i_j}$, where $\r_i$ is the rank of $\Dd_i$ and $E_{\g^i_j}$ are
%the root vectors associated to the strongly orthogonal roots $\g^i_j
%\in \Lambda_i$. 
%The Shilov boundary $\cs_{\Dd}$ is the $G_\Dd$ orbit of  $\sum_{i=}^n \sum_{j=1}^{\r_i}
%E_{\g^i_j}$, so we have 
%\bqn
%\cs_{\Dd} = G_\Dd \left(\sum_{i=}^n \sum_{j=1}^{\r_i}
%E_{\g^i_j}\right) = \prod_{i=1}^n G_{\Dd_i}\left(\sum_{j=1}^{\r_i}
%E_{\g^i_j}\right) =  \cs_{\Dd_1} \times \cdots \times \cs_{\Dd_n}.
%\eqn
%
%For the second statement note that since $\r_{\Dd'} = \r_{\Dd}$ we can
%choose a maximal set of strongly orthogonal roots for $\Dd$ which at
%the same time is a maximal set of strongly orthogonal roots for
%$\Dd'$. 
%Then 
%\bqn
%\cs_{\Dd'} = G_{\Dd'} \left(\sum_{j=1}^{\r_\Dd'}E_{\g_j}\right) \subset  G_{\Dd}
%\left(\sum_{j=1}^{\r_\Dd}E_{\g_j}\right) = \cs_{\Dd}, 
%\eqn
%where $E_{\g_j}$ are the root vectors associated to the strongly
%orthogonal roots. 
%\end{proof}

The relationship between the geodesic ray compactification $\Dd(\infty)$ of $\Xx$
and the boundary $\partial\Dd$ of the domain $\Dd$ is far from being simple.  For example,
a point in $\Dd(\infty)$ does not uniquely determine one in $\partial\Dd$;
this is however true if the endpoint of a geodesic ray lies in the Shilov boundary.
In fact we have:

\begin{lemma}\cite[Theorem~9.11]{Loos_notes_UCI}\label{lem:geodesics}
Let $z\in \cs_\Dd$ be a point in the Shilov boundary of $\Dd$
and let $\eta_1,\eta_2:[0,\infty)\to\Dd$ geodesic rays such that 
\be
\item $\lim_{t\to \infty} \eta_1(t) = z$, and
\item $\sup_{t\geq0} d_\Dd\big(\eta_1(t), \eta_2(t)\big) <\infty$. 
\ee
Then $\lim_{t\to \infty} \eta_2(t) = z$.
\end{lemma}
\begin{proof}[Sketch of the proof]
We can assume that the bounded symmetric domain is irreducible and that $z = eQ \in \cs = G_\Dd/Q$.  
The geodesic $\eta_1$ converges to $z\in \cs$ 
if and only if the stabilizer $\Stab_{G_\Dd}(\eta_1)$ is a parabolic subgroup $P$ which is contained in the maximal parabolic subgroup 
$Q<G_\Dd$. 
That $\sup_{t\geq0} d_\Dd\big(\eta_1(t), \eta_2(t)\big) <\infty$ implies that 
$\Stab_{G_\Dd}(\eta_1) = \Stab_{G_\Dd}(\eta_2) = P$. Now $P<Q$ and $P$ cannot be contained in any other conjugate of $Q$, hence $\lim_{t\to \infty} \eta_2(t) = z$.
\end{proof}

We shall as usual say that a geodesic ray $r:[0,\infty)\to\Dd$ is {\em of type $P$},
where $P$ is a parabolic subgroup of $G_\Dd$, if the stabilizer of the point 
in $\Dd(\infty)$ defined by $r$ is $P$ or, what amounts to the same, if
\bqn
P=\left\{g\in G_\Dd:\, \sup_{t\geq0} d_\Dd\big(g\,r(t),r(t)\big)<+\infty\right\}\,.
\eqn
By way of example, we notice that the geodesic
\bq\label{eq:geodesic}
\ba
r_0:[0,\infty)&\longrightarrow\qquad\Dd\\
t\qquad&\mapsto\sum_{j=1}^r\tanh(t)E_{\gamma_j}
\ea
\eq
is of type $Q$; this is the geodesic contained in $\Delta_0$  
connecting $0$ to $\sum_{j=1}^r\tanh(t)E_{\gamma_j}$ in $\ol\Dd$.
We should observe here that there are many geodesics connecting 
$0$ to $\sum_{j=1}^r E_{\gamma_j}$ in $\ol\Dd$, 
and they need not be at finite distance from $r_0$.  
A typical example is given by 
\bqn
t\mapsto\sum_{j=1}^r\tanh(a_jt)E_{\gamma_j}
\eqn
where $0<a_1<\dots<a_r$. However we have the following:

\begin{prop}\label{prop:3.4} For any $x\in\Dd$ and $z\in\cs_\Dd$, 
there is a unique diagonal disk $\Delta_{x,z}\subset\Dd$
with $\{x,z\}\subset\ol{\Delta_{x,z}}$.
Moreover, if $r_{x,z}$ denotes the unique geodesic ray in $\Delta_{x,z}$
joining $x$ to $z$, then $r_{x,z}$ is of type $Q_z:=\stab_{G_\Dd}(z)$.
Furthermore, for every $x_1,x_2\in\Dd$ and $z\in\cs_\Dd$, 
we have that 
\bqn
\sup_{t\geq0}d_\Dd\big(r_{x_1,z}(t),r_{x_2,z}(t)\big)<+\infty\,.
\eqn
\end{prop}
\begin{proof} Concerning the existence of such a disk,
observe that the diagonal $G_\Dd$-action on $\Dd\times\cs_\Dd$ is transitive;
indeed $Q$ acts transitively on $\Dd$.  Thus we may assume that 
$x=0$ and $z=\sum_{j=1}^r\tanh(t)E_{\gamma_j}$.  But then
$\Delta_0$ and $\r_0$ (see (\ref{eq:geodesic}) are the sought for objects.

Let now for the moment $r^\Delta_{x,z}$ denote the geodesic 
joining $x$ to $z$ inside $\ol\Delta$, where $\Delta$ is a diagonal disk.
Let $z\in\cs_\Dd$, and consider $x\in\Delta$, $x'\in\Delta'$ both 
diagonal disks with $z\in\partial\Delta\cap\partial\Delta'$.
Then there is $g\in G_\Dd$ with 
\bqn
g(\Delta)=\Delta'\,,\qquad gx=x'\,,\qquad gz=z\,,
\eqn
that is $g\in Q_z$.  In particular
\bq\label{eq:geodesics}
g\big(r^\Delta_{x,z}\big)=r^{\Delta'}_{x',z}\,.
\eq
Let $\Dd=\Dd_1\times\dots\times\Dd_n$ be a decomposition into irreducible components
and, accordingly, $G_\Dd=G_{\Dd_1}\times\dots\times G_{\Dd_n}$,
$x=(x_1,\dots,x_n)$, $x'=(x'_1,\dots,x'_n)$ and $z=(z_1,\dots,z_n)$.
Now if $p_i:\Dd\to\Dd_i$ denotes the projection on to the $i$-th factor, 
we observe that $p_i\big(r^\Delta_{x,z}\big)$ is a ray 
with parametrization proportional to the arclength and of type $Q_{z_i}$.
Since now $Q_{z_i}$ is maximal parabolic, 
there is a unique geodesic ray
$r_i^{x_i}:[0,\infty)\to\Dd_i$ starting at $x_i$ of type $Q_{z_i}$
and hence
\bqn
r^\Delta_{x,z}(t)=\big(r_i^{x_i}(a_it)\big)
\eqn
for some $a_i>0$.  Similarly,
\bqn
r^{\Delta'}_{x',z}(t)=\big(r_i^{x'_i}(b_it)\big)
\eqn
for some $b_i>0$.  If now $g=(g_1,\dots,g_n)$,
according to (\ref{eq:geodesics}) we have
\bqn
g_i\big(r_i^{x_i}(a_it)\big)=r_i^{x'_i}(b_it)\,,\text{ for all }t\geq0\,,
\eqn
which implies, since $g$ is an isometry, that $a_i=b_i$ for $1\leq i\leq n$.
Finally, since $g_i\in Q_{z_i}$, we have that 
\bqn
\sup_{s\geq0}d_{\Dd_i}\big(r_i^{x_i}(s),r_i^{x'_i}(s)\big)<+\infty
\eqn
and hence
\bq\label{eq:supd}
\sup_{t\geq0}d_\Dd\big(r^{\Delta}_{x,z}(t),r^{\Delta'}_{x',z}(t)\big)<+\infty\,.
\eq

It remains to show the uniqueness assertion.  For that,
let $\Delta$, $\Delta'$ be diagonal disks with $x\in\Delta\cap\Delta'$
and $z\in\partial\Delta\cap\partial\Delta'$.
Because of (\ref{eq:supd}) we have that
\bqn
r^{\Delta}_{x,z}(t)=r^{\Delta'}_{x',z}(t)
\eqn 
for all $t\geq0$.  Thus the holomorphic disks $\Delta,\Delta'$
contain a half line in common and hence coincide.
\end{proof}

\subsection{The Bergmann Cocycle and Maximal Triples}
Let us denote by $\ol{\Dd}$ the closure of $\Dd$ in $\frakp_+$. 
\begin{defi}
We define
\bqn
\ba
\ol{\Dd}^{(2)}:=&\big\{(z,w)\in\ol{\Dd}\times\ol{\Dd}:\,h_\Dd(z,w)\neq0\big\}
\\
\ol{\Dd}^{[2]} :=& \left\{ (z,w) \in \ol{\Dd}^2 \, :\, 
\text{ there exists some geodesic in } \Dd \text{ connecting $z$ to $w$}\right\}\,.
\ea
\eqn
\end{defi}
Then $\ol{\Dd}^{(2)}$ is a star shaped domain and is the maximal subset of $\ol{\Dd}^2$ 
to which the Bergmann kernel $k_\Dd$ extends continuously.  Moreover, expanding
\cite[Proposition~4.1]{Clerc_Orsted_2} to the nonirreducible case, we have that 
\bqn
\ol{\Dd}^{[2]} \subset \ol{\Dd}^{2}\,. 
\eqn

If 
\bqn
\arg k_\Dd: \ol{\Dd}^{(2)} \to \RR
\eqn
is the continuous determination of the argument of $k_\Dd$ vanishing on the diagonal of $\Dd^2$,
then the integral $\int_{\Delta \subset \Xx} \omega_\Xx$ can be expressed in terms of the 
$\arg k_\Dd$.

\begin{prop}\cite{Clerc_Orsted_2}
Let $\Delta(x,y,z) \subset \Dd$ be a smooth oriented triangle 
with geodesic sides and vertices $x,y,z \in \Dd$.
Then 
\bqn
\int_{\Delta(x,y,z)} \omega_\Dd = - [\arg k_\Dd(x,y) + \arg k_\Dd(y,z) + \arg k_\Dd(z,x)] 
\eqn
\end{prop}
 
Define the following subsets of $\ol{\Dd}^3$:
\bqn
\ol{\Dd}^{(3)} &:= &\big\{ (z_1, z_2, z_3) \in \ol{\Dd}^3 \, :\, (z_i, z_j) \in \ol{\Dd}^{(2)} \text{ for all } i\neq j\big\}\\
\ol{\Dd}^{[3]} &:= &\big\{ (z_1, z_2, z_3) \in \ol{\Dd}^3 \, :\, (z_i, z_j) \in \ol{\Dd}^{[2]} \text{ for all } i\neq j\big\}\\
\eqn

\begin{defi}
The Bergmann cocycle 
\bqn
\beta_\Dd: \ol{\Dd}^{(3)} \to \RR 
\eqn
is defined by 
\bqn
\beta_\Dd(x,y,z) := - \frac{1}{2\pi}\big[\arg k_\Dd(x,y) + \arg k_\Dd(y,z) + \arg k_\Dd(z,x)\big]\,.
\eqn
\end{defi}

It is a $G_\Dd$-invariant alternating continuous function, which satisfies the cocycle identity $d\beta(z_1, z_2, z_3, z_4) = 0$ whenever $(z_i, z_j, z_k) \in \ol{\Dd}^{(3)}$ for all $1\leq i,j,k\leq 4$ pairwise distinct. 

In terms of the decomposition $\Dd=\Dd_1\times\dots\times\Dd_n$ into irreducible factors,
we have the following formulas for the above mentioned objects:

\begin{align}
&h_\Dd\big((z_1,\dots,z_n),(w_1,\dots,w_n)\big)=\prod_{i=1}^nh_{\Dd_i}(z_i,w_i)\label{eq:1}\\
&\ol\Dd^{(2)}=\left\{(z,w)\in\ol\Dd^2:\,(z_i,w_i)\in\ol\Dd_i^{(2)},\,1\leq i\leq n\right\}\label{eq:2}\\
&\ol\Dd^{[2]}=\left\{(z,w)\in\ol\Dd^2:\,(z_i,w_i)\in\ol\Dd_i^{[2]},\,1\leq i\leq n\right\}\label{eq:3}\\
&\text{analogous formulas for }\ol\Dd^{(3)}\text{ and }\ol\Dd^{[3]}\label{eq:4}\\
&\beta_\Dd(x,y,z)=\sum_{i=1}^n\beta_{\Dd_i}(x_i,y_i,z_i)
   \text{ whenever }(x,y,z)\in\ol\Dd^{(3)}\label{eq:5}\,.
\end{align}

From \eqref{eq:5} and Theorem~\ref{thm:dtco} we deduce that $|\beta_\Dd|\leq\frac{\r_\Dd}{2}$
and clearly 
\bq\label{eq:sum}\quad\qquad
 \beta_\Dd(x,y,z)=\frac{\r_\Dd}{2}
\text{ if and only if }
 \beta_{\Dd_i}(x_i,y_i,z_i)=\frac{\r_{\Dd_i}}{2}\text{ for all }1\leq i\leq n\,.
\eq

\begin{thm}\label{results_c_o}
Suppose that $(x,y,z)\in\ol\Dd^{[3]}$ is such that $\beta_\Dd(x,y,z)=\frac{\r_\Dd}{2}$, then
\be
\item the points $x,y,z$ lie on the Shilov boundary $\cs_\Dd$, and
\item there exists a unique diagonal disk $d: \ol{\DD} \to \ol{\Dd}$ such that 
$d(1) = x$, $d(i) = y$, $d(-1) = z$. Moreover
\bqn
d(\DD)\subset\big\{p\in\Dd:\,p\text{ is fixed by }\stab_{G_\Dd}(x,y,z)\big\}
\eqn 
with equality if $\Dd$ is irreducible.
\item The group $G_\Dd$ acts transitively on the set
\bqn
\left\{(x,y,z)\in\ol\Dd^{[3]}:\,\beta_\Dd(x,y,z)=\frac{\r_\Dd}{2}\right\}
\eqn
of maximal triples.
\ee
\end{thm}

\begin{proof}
All the above assertions are due to Clerc and \O rsted in the irreducible case,
\cite{Clerc_Orsted_2}.  In the general case, the first assertion follows from
\eqref{eq:sum} and Lemma~\ref{lem:3.2}(1).  

In the second assertion, 
only the uniqueness needs to be verified, but this follows easily from the fact
that a totally geodesic map $\DD\to \DD$ is necessarily isometric.

The last assertion follows immediately from \eqref{eq:sum}, \eqref{eq:4}
and the irreducible case.
\end{proof}

\subsection{On Subdomains of Maximal Rank}
The main goal of this section is to show that if $\Dd',\Dd''$ are subdomains
of $\Dd$ of maximal rank, that is $\r_\Dd=\r_{\Dd'}=\r_{\Dd''}$,
whose Shilov boundaries coincide, then $\Dd'=\Dd''$.
We begin with the following

\begin{lemma}\label{lem:3.9} Let $\Dd'\subset\Dd$ be a subdomain of maximal rank.  Then:
\be
\item $k_{\Dd'}=k_\Dd|_{\Dd'}$;
\item $\ol{\Dd'}^{[3]}=\big(\ol{\Dd'}\big)^3\cap\ol\Dd^{[3]}$;
\item $\beta_{\Dd'}(x,y,z)=\beta_\Dd|_{\ol{\Dd'}^{[3]}}(x,y,z)$
for all $(x,y,z)\in\ol{\Dd'}^{[3]}$.
\ee
\end{lemma} 

\begin{proof} (1) We observed already that the normalized Bergmann metric on $\Dd'$ 
is the restriction of the normalized Bergmann metric on $\Dd$
(see Corollary~\ref{cor:2.20}), thus the first assertion follows readily.

(2) The second assertion is obvious.

(3) To see the third assertion, observe first of all that if $(x,y,z)\in\Dd'^3$,
then the equality follows from (1).  Furthermore, continuity of 
the Bergmann cocycle and the fact that $\Dd'^3$ is dense in $\ol{\Dd'}^{[3]}$
complete the proof.
\end{proof}

\begin{lemma}\label{lem:x+1} Let $\Dd'\subset\Dd$ be a subdomain of maximal rank,
$A_{\Dd'}$ the associated subgroup of Hermitian type and 
$f:\DD\to\Dd'$ a diagonal embedding.  Then
\bqn
\Nn_{G_\Dd}(\cs_{\Dd'}) = \Nn_{G_\Dd}(\Dd') \subset A_{\Dd'} \cdot \Zz_{G_\Dd}(f)\,
\eqn
where $\Zz_{G_\Dd}(f)=\left\{g\in G_\Dd:\, gx=x\text{ for all }x\in f(\DD)\right\}$.
\end{lemma}

\begin{proof} Let $x:=f(1)$, $y:=f(i)$ and $z:=f(-1)$, 
and pick $g\in\Nn_{G_\Dd}(\cs_{\Dd'})$.
We have that since $(x,y,z)\in\cs_{\Dd'}^{[3]}$ then $(gx,gy,gz)\in\cs_{\Dd'}^{[3]}$
and, using Lemma~\ref{lem:3.9}(3),
\bqn
 \frac{\r_{\Dd}}{2}
=\beta_{\Dd'}(x,y,z)
=\beta_\Dd(x,y,z)
=\beta_\Dd(gx,gy,gz)
=\beta_{\Dd'}(gx,gy,gz)\,,
\eqn
which implies by Theorem~\ref{results_c_o}(3) that there exists $h\in A_{\Dd'}$ 
with $hx=gx$, $hy=gy$, $hz=gz$ and thus $g\in A_{\Dd'}\cdot\stab_{G_\Dd}(x,y,z)$.
Since $\stab_{G_\Dd}(x,y,z)\subset\Zz_{G_\Dd}(f)$ (see Theorem~\ref{results_c_o}(2)), 
we obtain that $\Nn_{G_\Dd}(\cs_{\Dd'})\subset A_{\Dd'}\cdot\Zz_{G_\Dd}(f)$.

Let now $p\in f(\DD)\subset\Dd'$. Then
$\Nn_{G_\Dd}(\cs_{\Dd'})\cdot p\subset A_{\Dd'}(p)=\Dd'$,
and since $\Nn_{G_\Dd}(\cs_{\Dd'})=\Nn_{G_\Dd}(\cs_{\Dd'})\cdot A_{\Dd'}$,
we obtain that $\Nn_{G_\Dd}(\cs_{\Dd'})\Dd'=\Dd'$ and hence
$\Nn_{G_\Dd}(\cs_{\Dd'})\subset\Nn_{G_\Dd}(\Dd')$.  
The opposite inclusion is clear.
\end{proof}

\begin{prop}\label{prop:3.11} Let $\Dd',\Dd''$ be subdomains of $\Dd$ of maximal rank,
and assume that $\cs_{\Dd'}=\cs_{\Dd''}$.  Then $\Dd'=\Dd''$.
\end{prop}

\begin{proof}  Let $x,y,z\in\cs_{\Dd'}=\cs_{\Dd''}$ with 
$\beta_\Dd(x,y,z)=\frac{\r_\Dd}{2}$.  Then 
\bqn
 \beta_{\Dd'}(x,y,z)=\frac{\r_{\Dd'}}{2}
\quad\text{ and }\quad
 \beta_{\Dd''}(x,y,z)=\frac{\r_{\Dd''}}{2}\,,
\eqn
and there are diagonal disks 
\bqn
 f_{\Dd'}:\ol\DD\to\ol{\Dd'}
\quad\text{ and }\quad
 f_{\Dd''}:\ol\DD\to\ol{\Dd''}
\eqn
with 
\bqn
\ba
 f_{\Dd'}(1)=&x=f_{\Dd''}(1)\\
 f_{\Dd'}(i)=&y=f_{\Dd''}(i)\\
 f_{\Dd'}(-1)=&z=f_{\Dd''}(-1)\,.
\ea
\eqn
Those are also diagonal disks in $\Dd$ and hence by uniqueness 
we have $f_{\Dd'}=f_{\Dd''}$, and in particular $\Dd'\cap\Dd''\neq\emptyset$.
Pick now $p\in\Dd'\cap\Dd''$ and apply Lemma~\ref{lem:x+1}
to obtain
\bqn
\Dd'=\Nn_{G_\Dd}(\cs_{\Dd'})\cdot p=\Nn_{G_\Dd}(\cs_{\Dd''})\cdot p=\Dd''\,.
\eqn
\end{proof}

%% file: structure_one.tex
\section{Structure Theorem for Tight Embeddings, I}
The main objective in this section is to prove the following structure
theorem.  
\begin{thm}\label{thm:tight_embedding_one} 
Let $H, G$ be Lie groups of Hermitian type, 
$\rho: H \to G$ a continuous tight homomorphism 
and $f: \Dd' \to \Dd$ the corresponding $\rho$-equivariant tight
totally geodesic map. 
Then $f$ extends continuously to a $\rho$-equivariant map  
\bqn
\check{f}: \cs_{\Dd'} \to \cs_{\Dd}\,.
\eqn 
Moreover the centralizer $\Zz_G\big(\rho(H)\big)$ is compact.
\end{thm}
\subsection{The Case of the Poincar\'e Disk}
We will first prove Theorem~\ref{thm:tight_embedding_one} in the case when 
$\Dd' = \DD$ is the Poincar\'e disk. 

\begin{prop}\label{prop:baby}
Let $L$ be a finite covering of $\PU(1,1)$ and $G$ a group of Hermitian type. 
Let $\rho: L\to G$ a continuous tight homomorphism with finite kernel 
and $f: \DD \to \Dd$ the corresponding totally geodesic tight embedding.
Then 
\begin{enumerate}
\item $f$ extends continuously to a map 
\bqn
\check{f}: \partial \DD \to \partial \Dd\,, 
\eqn
which is $\rho$-equivariant and has image $\check{f} (\partial \DD) \subset \cs_\Dd$.
\item If $x\neq y$ in $\partial\DD$, then $\check f(x)$ and $\check f(y)$ are transverse.
\item The centralizer $\Zz_{G}\big(\rho(L)\big)$ is compact.
\end{enumerate}
\end{prop}

\begin{proof} Observe that since $f$ is equivariant, there exists a constant $c>0$ such that 
\bqn
d_\Dd\big(f(x),f(y)\big)=c\,d_\DD(x,y)\,.
\eqn
Next, let $r_1,r_2:\RR^+\to\DD$ be geodesic rays representing a given point
$\xi\in\partial\DD$ and $a\geq0$ such that $\lim_{t\to\infty}d_\DD\big(r_1(t),r_2(t+a)\big)=0$.
Then
\bqn
\lim_{t\to\infty}d_\Dd\big(f\big(r_1(t)\big),f\big(r_2(t+a)\big)\big)=0
\eqn
and, by Lemma~\ref{lem:distance},
\bqn
\lim_{t\to\infty}\big\|f\big(r_1(t)\big)-f\big(r_2(t+a)\big)\big\|_{eucl}=0\,,
\eqn
which shows that the geodesics $t\mapsto f\big(r_1(t)\big)$ and
$t\mapsto f\big(r_2(t)\big)$ have the same endpoints in $\partial\Dd$.  
This produces a well defined and equivariant (continuous) extension 
$\check f:\partial\DD\to\partial\Dd$ of $f$. 
Observe that for all $x\neq y$ in $\partial \DD$, $(\check{f}(x), \check{f}(y)) \in \ol{\Dd}^{[2]}$ 
since $x,y$ and thus $\check{f}(x), \check{f}(y)$ are joined by a geodesic. 

Now, for the (normalized) K\"ahler forms $\omega_\Dd \in \Omega^2(\Dd)^G$ and 
$\omega_\DD\in \Omega^2(\DD)^L$ we have since $f$ is tight
\bqn
f^*(\omega_\Dd) = \eps \omega_\DD\,, 
\eqn
where $|\eps| = \r_\Dd$. 
Composing if necessary with an orientation reversing isometry of $\DD$ 
we may assume $\eps = \r_\Dd$. 

This implies by integration over simplices with geodesic sides in $\DD$ and 
continuity of $\check{f} $, 
as well as of $\beta_\Dd$ on $\ol{\Dd}^{[3]}$ and of $\beta_\DD$ on $\ol{\DD}^{[3]}$ that 
\bqn
\beta\big(\check{f}(x), \check{f}(y), \check{f}(z)\big) = \r_\Dd \beta_\DD(x,y,z) 
\eqn
for all $(x,y,z) \in \ol{\DD}^{[3]}$. 

Applying this to a positively oriented triple $(x,y,z)$ we get 
\bqn
\beta\big(\check{f}(x), \check{f}(y), \check{f}(z)\big) = \frac{\r_\Dd}{2}
\eqn
and hence, by Theorem ~\ref{results_c_o},  $\check{f}(x) \in \cs_\Dd$.
This shows (1).  The second assertion follows from the fact we already 
remarked that $\big(\check f(x),\check f(y)\big)\in\ol\Dd^{[2]}$ if $x\neq y$. 
For the third assertion, 
let $\eta: \RR^+ \to \DD$ be a geodesic ray with $\lim_{t\to\infty}\eta(t) = x$. 
Then $f(\eta)$ is a geodesic ray in $\Dd$ converging to $\check{f}(x)$.  
For $g \in \Zz_G(\rho(L))$ the geodesic ray $g \cdot f(\eta)$ is at bounded distance from $f(\eta)$,
hence Lemma~\ref{lem:geodesics} implies that 
\bqn 
 g\check{f}(x)
=\check{f}(x)\,. 
\eqn 
In particular 
$\Zz_G\big(\rho(L)\big) \subset\Stab_G\big(\check{f}(x), \check{f}(y), \check{f}(z)\big)$ 
which by Theorem~\ref{results_c_o}~(2) is compact.
\end{proof}

\begin{cor}\label{lem:munonzero}
Let $H,G$ be groups of Hermitian type and let $\{\khib\}_{i=1}^n\in\hcb^2(H,\RR)$
be the basis of $\hcb^2(H,\RR)$ corresponding to the decomposition $\Yy=\Yy_1\times\dots\times\Yy_n$
into irreducible factors of the symmetric space $\Yy$
associated to $H$. Let $\rho: H\to G$ be a tight homomorphism and 
assume that $\rho^*_{\rm b} \kgb = \sum_{i=1}^n \lambda_i\khib$.
If $H=H_1\cdots H_n$ is the decomposition of $H$ into connected almost simple groups
where $H_i$ corresponds to $\Yy_i$, then $\lambda_i = 0$ 
if and only if $H_i$ is in the kernel of $\rho$. 
\end{cor}

\begin{proof} 
If $f:\Yy\to\Xx$ is a tight $\rho$-equivariant map, then
\bqn
f^\ast(\omega_\Xx)=\sum_{i=1}^n\lambda_i\omega_{\Yy,i}\,,
\eqn
where $\omega_{\Yy,i}=p_i^\ast(\omega_{\Yy_i})$. 
Assume that $\lambda_i\neq0$ for $1\leq i\leq\ell$ and
$\lambda_{\ell+1}=\cdots=\lambda_n=0$.  For $1\leq i\leq \ell$ define 
$t_i:\DD\to\Yy_i$ to be the embedding as a diagonal disk 
(Definition~\ref{defi:maxi-polydisk})
composed with an isometry of $\DD$ reversing the orientation in the case
in which $\lambda_i<0$, and let $\rho_i:\mathrm{SU}(1,1)\to\Is(\Yy_i)^\circ$
be the associated homomorphism. Let $b_i\in\Yy_i$ be a basepoint and define
\bqn
\ba
t:\DD&\longrightarrow\qquad\Yy_1\times\dots\times\Yy_n\\
z&\mapsto\big(t_1(z),\dots,t_\ell(z),b_{\ell+1},\dots,b_n\big)\,.
\ea
\eqn
and
\bqn
\ba
\pi:\mathrm{SU}(1,1)&\longrightarrow\qquad\Is(\Yy)^\circ\\
g&\mapsto\big(\rho_1(g),\dots,\rho_\ell(g),e,\dots,e\big)\,.
\ea
\eqn
Taking into account that $H$ is a finite extension of $\Is(\Yy)^\circ$, 
let 
\bqn
\widetilde\pi:L\to H
\eqn 
be the lift of $\pi$ to a finite extension $L$ of $\mathrm{SU}(1,1)$.
Then 
\bqn
 t^\ast\left(\sum_{i=1}^n\lambda_i\omega_{\Yy,i}\right)
=\left(\sum_{i=1}^n|\lambda_i|\r_{\Yy_i}\right)\omega_\DD
=\rx\omega_\DD\,,
\eqn
where the last equality follows from the fact that $f$ is tight (Corollary~\ref{cor:tight_rank}).
Thus
\bqn
f\circ t:\DD\to\Xx
\eqn
is tight and equivariant with respect to the homomorphism 
\bqn
\rho\circ\widetilde\pi:L\to G\,.
\eqn
Let now $H=H_1\cdot\dots\cdot H_n$ be the decomposition of $H$ into connected almost simple
groups, where $H_i$ is a finite extension of $\Is(\Yy_i)^\circ$.  
In particular, for $\ell+1\leq j\leq n$, $H_j$ commutes with $\widetilde\pi(L)$
and hence $\rho(H_j)$ commutes with $(\rho\circ\widetilde\pi)(L)$
which implies that, for $\ell+1\leq j\leq n$, $\rho(H_j)$ is contained in 
$\Zz_G\big(\rho\widetilde\pi(L)\big)$
which is compact in virtue of Proposition~\ref{prop:baby}, 
and hence $\rho(H_j)=e$.

The converse, namely that $\rho(H_j)=e$ implies that $\lambda_j=0$ is clear.
\end{proof}

\subsection{Positivity}
Let $G$ be a group of type (RH).  We shall use freely the notation from \S~\ref{subsec:hss}.
In this section we prove that the notion of tightness does not depend on the choice 
of the specific K\"ahler class $\kgb \in \hcb^2(G,\RR)$ which we used to define it
but indeed it depends only on the choice of a $G_\Xx$-invariant complex structure on $\Xx$. 
In the case when $\Xx$ is irreducible this is immediate from $\hcb^2(G, \RR) = \RR\, \kgb$. 
In the general case however,
one could have some ``cancellations'' coming from different factors,
but we are going to set up conditions which will allow us some freedom to
choose the K\"ahler classes according to the context. 

Let again $\Xx=\Xx_1\times\dots\times\Xx_n$ be the decomposition
into irreducible factors. Then any choice of $G_\Xx$-invariant complex structure $\Jj_\Xx$
determines a $G_{\Xx_i}$-invariant complex structure $\Jj_{\Xx_i}$ on $\Xx_i$
and hence an orientation on $\hcb^2(G_{\Xx_i},\RR)$.
Conversely, any choice of orientation on each $\hcb^2(G_{\Xx_i},\RR)$ 
determines a complex structure on $\Xx$.
\begin{defi}
A bounded cohomology class $\alpha \in \hcb^2(G,\RR)$ is {\em positive} if 
\bqn
\alpha = \sum_{i=1}^n \mu_i\kgib
\eqn
with $\mu_i \geq0$ for all $ i= 1, \cdots, n$ and {\em strictly positive} 
if the $\mu_i>0$, for all $i=1,\dots,n$.

The cone of {\em positive K\"ahler classes} in $\hcb^2(G,\RR)$ 
is denoted by $\hcb^2(G,\RR)^{\geq0}$ and the cone of {\em strictly positive K\"ahler
classes} by $\hcb^2(G,\RR)^{>0}$.
\end{defi}

Note that the cone  $\hcb^2(G,\RR)^{\geq0}$ depends only on the complex structure $\Jj$. 
In fact $\hcb^2(G,\RR)^{>0}$ coincides with the set of bounded K\"ahler classes
associated to any $G$-invariant Hermitian metric on $\Xx$
compatible with the complex structure $\Jj$;
in particular we have that $\kgb\in\hcb^2(G,\RR)^{>0}$.

\begin{prop}\label{prop:positive}
Let $\rho:H\to G$ be a homomorphism of a locally compact group $H$ into 
a group $G$ of type (RH).  Then the following are equivalent:
\be
\item $\rho$ is tight;
\item $\rho$ is $\alpha$-tight for some $\alpha\in\hcb^2(G,\RR)^{>0}$;
\item $\rho$ is $\alpha$-tight for all $\alpha\in\hcb^2(G,\RR)^{>0}$;
\item $\rho$ is $\alpha$-tight for all $\alpha\in\hcb^2(G,\RR)^{\geq0}$.
\ee
\end{prop} 

This is a consequence of the special Banach space structure of $\hcb^2$. 
\begin{lemma}\label{lem:Banach_norm}
Let $V$ be a Banach space. Let $v_i\in V$, $i=1,\dots,k$, be vectors 
such that 
\bqn
\left\|\sum_{i=1}^nv_i\right\|=\sum_{i=1}^n\|v_i\|\,.
\eqn
Then for every real numbers $\mu_1,\dots,\mu_n\geq0$, we have that 
\bqn
\left\|\sum_{i=1}^n\mu_i\,v_i\right\|=\sum_{i=1}^n\mu_i\,\|v_i\|\,.
\eqn
\end{lemma}
\begin{proof}
In virtue of the Hahn--Banach theorem the norm of a vector $w\in V$ is 
given by 
\bqn
\|w\|=\sup\big\{|\lambda(w)|:\,
   \lambda:V\to\RR\hbox{ is a linear form of norm }1\big\}\,.
\eqn 
By hypothesis, if we fix $\epsilon>0$, there exists $\lambda:V\to\RR$ a linear form 
of norm 1 such that 
\bqn
 \lambda\left(\sum_{i=1}^nv_i\right)
>\left\|\sum_{i=1}^nv_i\right\|-\eps
=\left(\sum_{i=1}^n\|v_i\|\right)-\eps\eqn
From this and the fact that 
$\lambda(v_i)\leq\|v_i\|$,
we must have that  
\bqn
\lambda(v_i)\geq \|v_i\|-\epsilon
\eqn
for all $i=1,\dots,n$, and hence, if $\mu_i\geq0$,
\bqn
\mu_i\lambda(v_i)\geq\mu_i\|v_i\|-\mu_i\eps\,.
\eqn
But then
\bqn
    \lambda\left(\sum_{i=1}^n\mu_iv_i\right)
\geq\left(\sum_{i=1}^n\mu_i\|v_i\|\right)-\eps\sum_{i=1}^n\mu_i\,,
\eqn
which, since $\eps$ is arbitrary, shows the assertion.
\end{proof}

\begin{proof}[Proof of Proposition~\ref{prop:positive}]
We start by showing that (1)$\Rightarrow$(4).  We first verify 
that the vectors $v_i:=\rho_\mathrm{b}^\ast(\kgib)$
satisfy the hypotheses of Lemma~\ref{lem:Banach_norm}.
We have that
\bqn
 \left\|\sum_{i=1}^nv_i\right\|
=\left\|\sum_{i=1}^n\rho_{\rm b}^\ast\kgib\right\|
=\|\rho_\mathrm{b}^\ast\kgb\|
=\|\kgb\|\,,
\eqn
where the last equality follows from the fact that $\rho$ is tight. 
Moreover, Lemma~\ref{lem:tight_properties}(4) implies that 
\bq\label{eq:4.5}
 \|v_i\|
=\|\rho_{\rm b}^\ast(\kgib)\|
=\|\kgib\|\,,
\eq
and hence
\bqn
 \sum_{i=1}^n\|v_i\|
=\sum_{i=1}^n\|\kgib\|
=\|\kgb\|\,.
\eqn
Thus $\left\|\sum_{i=1}^nv_i\right\|=\sum_{i=1}^n\|v_i\|$,
and applying Lemma~\ref{lem:Banach_norm} we get
\bqn
\ba
  \|\rho_{\rm b}^\ast\alpha\|
&=\left\|\sum_{i=1}^n\mu_i\rho_{\rm b}^\ast\kgib\right\|
 =\sum_{i=1}^n \|\mu_i\rho_{\rm b}^\ast\kgib\|\\
&=\sum_{i=1}^n\mu_i\|\kgib\|
 =\left\|\sum_{i=1}^n \mu_i\kgib\right\|
 =\|\alpha\|\,.
\ea
\eqn
Thus $\rho$ is $\alpha$-tight.

The implications (4)$\Rightarrow$(3)$\Rightarrow$(2) are obvious.

Finally, to see that (2)$\Rightarrow$(1), 
let $\alpha=\sum_{i=1}^n\lambda_i\kgib$ be strictly positive.
Then setting $v_i:=\lambda_i\kgib$ and $\mu_i:=\frac{1}{\lambda_i}$,
the argument above implies that if $\rho$ is $\alpha$-tight
then it is $\kgb$-tight.
\end{proof}

%When $G$ is simple the statement follows readily from the fact that
%$\hcb^2(G,\RR)= \RR\, \kgb$. 
%For the general case, let $\a\in \hcb^2(G,\RR)^{\geq 0}$ be a positive class, 
%\bqn
%\alpha = \sum_{i=1}^n \kgib\,,
%\eqn
%where $\mu_i>0$ and $\kgib\in\hcb^2(G,\RR)^{\geq 0}$, $i= 1, \cdots, n$.
%Tightness of $\rho$, Proposition~\ref{prop:bounded_properties}(3) 
%and Lemma~\ref{lem:tight_properties}(4) imply that  
%\bqn
%  \left\|\sum_{i=1}^n(\rho)_{\rm b}^\ast\kgib\right\|
% =\|\rho_{\rm b}^\ast\kgb\|
% =\|\kgb\|
% =\left\|\sum_{i=1}^n\kgib\right\|
% =\sum_{i=1}^n \big\|\kgib\big\|\,.
%\eqn
%%\bqn
%%\ba
%%  \left\|\sum_{i=1}^n(\rho_i)_{\rm b}^\ast\kxib\right\|
%%&=\left\|\rho_{\rm b}^\ast\left(\sum_{i=1}^n
%%     (p_i\circ q)_{\rm b}^\ast\kxib\right)\right\|
%% =\|\rho_{\rm b}^\ast\kgb\|
%% =\|\kgb\|\\
%%&=\left\|\left(\sum_{i=1}^n (p_i\circ q)_{\rm b}^\ast\kxib\right)\right\|
%% =\sum_{i=1}^n \big\|(p_i\circ q)_{\rm b}^\ast\kxib\big\|\\
%%&=\sum_{i=1}^n \big\|\rho_{\rm b}^\ast(p_i\circ q)_{\rm b}^\ast\kxib\big\|
%% =\sum_{i=1}^n \|(\rho_i)_{\rm b}^\ast\kxib\|\,.
%%\ea
%%\eqn
%Finally, because of Lemma~\ref{lem:Banach_norm}, tightness of $\rho_i$ and positivity of 
%$\alpha$, we conclude that
%\bqn
%\ba
%  \|\rho_{\rm b}^\ast\alpha\|
%&=\left\|\sum_{i=1}^n\mu_i(\rho)_{\rm b}^\ast\kgib\right\|
% =\sum_{i=1}^n \|\mu_i(\rho)_{\rm b}^\ast\kgib\|\\
%&=\sum_{i=1}^n\mu_i\|\kgib\|
% =\left\|\sum_{i=1}^n \mu_i\kgib\right\|
% =\|\alpha\|\,.
%\ea
%\eqn
%\end{proof}

\begin{defi}
A homomorphism $\rho: H \to G$ of groups of type (RH) is said to be {\em positive} if
$\rho_{\rm b}^*\kgb \in \hcb^2(H, \RR)^{\geq0}$ and {\em strictly positive}
if $\rho_{\rm b}^*\kgb \in \hcb^2(H, \RR)^{>0}$.
\end{defi}

The point of the next lemma is to provide a converse to Lemma~\ref{lem:tight_properties}(4),
for which we need the hypothesis of positivity. Remark that it will be essential 
that, with the norm on the continuous bounded cohomology,
we have that if $v,w$ are positive classes then $\|v+w\|=\|v\|+\|w\|$.

\begin{lemma}\label{lem:positive_projection}
Let $H,G$ be of type (RH) and let $\rho:H\to G$ be a continuous homomorphism.
With the notation in \S~\ref{subsec:hss}, if $\rho_i:= q\circ \rho: H\to G_{\Xx_i}$ 
is tight and positive for all $i= 1, \dots,n$ then $\rho$ is tight and positive.
\end{lemma} 

\begin{proof}
Since we have  
\bqn
\kgb = \sum_{i=1}^n \kgib \in
\hcb^2(G,\RR)\,,
\eqn 
then 
\bqn
  \rho_{\rm b}^* \kgb 
= \sum_{i=1}^n (\rho_i)_{\rm b}^* \kappa_{\Xx_i}^{\rm b}\,.
\eqn 
 
Since $(\rho_i)_{\rm b}^* \kappa_{\Xx_i}^{\rm b}$ are positive for all $i = 1, \dots, n$, 
$\rho_{\rm b}^*\kgb$ is positive; this, and the hypothesis that 
\bqn
\|(\rho_i)_{\rm b}^*\kxib\| = \|\kxib\|\,,
\eqn
allow us to deduce that 
\bqn
  \|\rho_{\rm b}^* \kgb\| 
= \left\|\sum_{i=1}^n (\rho_i)_{\rm b}^* \kxib\right\|
= \sum_{i=1}^n \big\|(\rho_i)_{\rm b}^* \kxib\big\|
= \sum_{i=1}^n \|\kxib\| = \|\kxb\| 
= \|\kgb\|\,.
\eqn
\end{proof}

\begin{lemma}\label{lem:composition_pos}
Let $H, G$ be Lie groups of type (RH), $L$ a locally compact group,
$\rho: L\to H$ a tight homomorphism, and $\psi: H\to G$ a positive tight homomorphism.
Then $\psi\circ\rho:L\to G$ is a tight homomorphism.
%Then $\psi\circ\rho: L \to G$ is positive if and
%only if $\psi: H \to G$ is positive.
%$f: \Yy \to \Xx$ and $g:\Xx \to \Mm$ be tight positive embeddings. Then $f\circ g: \Yy \to \Mm$ is tight and positive.
\end{lemma}
\begin{proof}
If $\psi$ is positive, $\psi^*_{\rm b} \kgb \in \hcb^2(H, \RR)^{\geq0}$. 
By Proposition~\ref{prop:positive}, 
if the homomorphism $\rho$ is tight, it is also  $\psi^*_{\rm  b}(\kgb)$-tight and 
Lemma~\ref{lem:composition} concludes the proof.
\end{proof}
%
%Let $\ol{H} = H/\Zz(H)$, $\ol{H} = H_1 \times \cdots \times H_n$ a decomposition into simple factors, 
%$q: \ol{H} \to H$ the quotient map 
%and $p_i: \ol{H} \to H_i$, $i= 1, \cdots, n$ the projection onto the $i$-th factor. We choose the same notation for $L$. 
%Then, since $\rho: L \to H$ is tight and positive, 
%the homomorphisms $p_i \circ q \circ \rho: L \to H_i$ are also tight and positive 
%for all $i= 1, \cdots, n$ 
%(Lemma~\ref{lem:positive_projection}). In particular 
%\bqn
%(p_i \circ q \circ \rho)^*_{\rm b} \khib = 
%\sum_{j=1}^k \lambda_j (p_j \circ q)^*_{\rm b}\kljb \in \hcb^2(L,\RR)^+
%\eqn
%with $\lambda_j>0$ for all $j= 1,\cdots,k$.
%Assume that $\psi: H \to G$ is not positive. Then $\psi^*_{\rm b} \kgb = \sum_{i=1}^n \mu_i(p_i \circ q)^*_{\rm b} \khib$ and there exists at least one $i$ such that $\mu_i <0$. 
%Thus 
%\bqn
%(\psi\circ\rho)^*_{\rm b} \kgb &= &\sum_{i=1}^n \mu_i(p_i \circ q \circ \rho)^*_{\rm b} \khib\\
%&=& \sum_{i=1}^n \mu_i \sum_{j=1}^k \lambda_j\kljb
%\eqn 
%with $\lambda_j>0$ for all $1\leq j\leq k$ and $\mu_i < 0$ for at least one $i$. In particular 
%\bqn
%\sum_{i=1}^n \mu_i \sum_{j=1}^k \lambda_j\kljb \notin \hcb^2(L,\RR)^+
%\eqn 
%and $\psi\circ\rho$ is not positive. 

\begin{lemma}\label{lem:makepositive}
Let $H,G$ be Lie groups of Hermitian type with associated symmetric spaces $\Yy$ and $\Xx$ 
with complex structures $\Jj_\Yy$ and $\Jj_\Xx$. 
Suppose that $\rho: H \to G$ is a tight homomorphism 
and $f: \Yy \to \Xx$ is the corresponding tight map.
Then there exists a complex structure $\Jj'$ on $\Yy$ 
such that $\rho$ is tight and positive with respect to $\Jj'$. 
If moreover $\ker\rho$ is finite, then this structure is unique.
\end{lemma}
\begin{proof}
Since $\rho$ is tight, we have 
$0\neq \rho_{\rm b}^* \kgb \in \hcb^2(H,\RR)$.
So, if $\Yy$ is irreducible, then $\rho: H \to G$ is either positive 
with respect to $\Jj_\Yy$ or with respect to $-\Jj_\Yy$.

In the case when $\Yy$ is not irreducible, let $\Yy= \Yy_1 \times \cdots \times \Yy_n$ be 
the decomposition into irreducible factors, and $\Jj_i$ the complex structure on $\Yy_i$ induced by $\Jj_\Yy$. 
%Using the notation as in Notation~\ref{notation_1} we have {\bf Change notation} 
We have 
\bqn
\rho_{\rm b}^* \kgb= \sum_{i=1}^n \mu_i \kappa_{H,i}^{\rm b}\,.
\eqn
%Assume that $\mu_i <0$ for all  $i \in I \subset \{1, \cdots, n\}$. 
%Then s
Set $\Jj'_i = \eps_i \Jj_i$, where $\eps_i =\mathrm{sign}(\mu_i)$ and 
let $\Jj'$ be the complex structure on $\Yy$ which induces the complex structure $\Jj'_i$ on $\Yy_i$. 
Let ${\kappa'}_{H,i}^{\rm b} = \mathrm{sign}(\mu_i) \khib \in
\hcb^2(H, \RR)$ 
be the basis vectors of $\hcb^2(H,\RR)$ corresponding to $\Jj'$. 
Then 
\bqn
 \rho_{\rm b}^* \kgb
=\sum_{i=1}^n \mathrm{sign}(\mu_i)\mu_i\kappa_{H,i}^{\rm b}\,, 
\eqn
so $\rho$ is positive with respect to $\Jj'$ and tight. 
In case $\ker\rho$ is finite we have that $\mu_i\neq0$ for all $1\leq i\leq n$
(Corollary~\ref{lem:munonzero}) and hence $\Jj'$ is unique.
\end{proof}

\begin{proof}[Proof of Theorem~\ref{thm:tight_embedding_one}]
Let $\Dd'$ and $\Dd$ be the bounded domain realizations respectively of $\Yy$ and $\Xx$. 
Let $f:\Dd' \to \Dd$ be a $\rho$-equivariant totally geodesic tight map.
Because of Lemma~\ref{lem:makepositive} we may assume that $f$ is positive.  
For every $x\in\Dd'$ and $z\in\cs_{\Dd'}$, let $\Delta_{x,z}\subset\Dd'$
be the unique diagonal disk given by Proposition~\ref{prop:3.4} and
\bqn
d_{x,z}:\ol\DD\to\ol{\Delta_{x,z}}
\eqn
the unique totally geodesic map with 
\bqn
d_{x,z}(0)=x\qquad\text{ and }\qquad d_{x,z}(1)=z\,.
\eqn
Then $f\circ d_{x,z}:\DD\to\Dd$ is tight (Lemma~\ref{lem:composition_pos}) and hence,
by Proposition~\ref{prop:baby}  extends to 
\bqn
\ol{f\circ d_{x,z}}:\;\partial\DD\to\cs_{\Dd}\,.
\eqn
We set 
\bqn
\check f_x(z):=\big(\ol{f\circ d_{x,z}}\big)(1)=\lim_{t\to\infty}f\big(r_{x,z}(t)\big)\,.
\eqn
If now $x'$ is another point in $\Dd$, we have that
\bqn
\sup_{t\geq0} d_{\Dd'}\big(r_{x,z}(t),r_{x',z}(t)\big)<+\infty
\eqn
and, since $f$ is totally geodesic, also 
\bqn
\sup_{t\geq0} d_{\Dd}\left(f\big(r_{x,z}(t)\big),f\big(r_{x',z}(t)\big)\right)<+\infty\,.
\eqn
Since $\check f_x(z)\in\cs_{\Dd}$, we deduce, by \eqref{eq:geodesics}, 
that $\check f_{x'}(z)=\check f_x(z)$; thus the extension
$\check f:\cs_{\Dd'}\to\cs_{\Dd}$ is independent of $x$ 
and hence $\rho$-equivariant.
\end{proof}

%% file: tubetype.tex
\section{Tight Embeddings and Tube Type Domains}
\label{sec:tubetype}
Let $\Xx$ be a Hermitian symmetric space and $\Dd$ its bounded
symmetric domain realization. We will use the concepts and notations
from \S~\ref{sec:shilov}. 

The real vectors $X_{\g_j} \in \frakp$ associated to the 
strongly orthogonal roots $\g_j \in \Lambda$, $ j= 1 , \dots , r$, (see
\S~\ref{subsec:shilov}) give rise to the 
Cayley element 
\bqn
c = \exp\left(\frac{\pi}{4} i \sum_{j=1}^r X_{\g_j} \right) \in G_\CC= \exp(\frakg_\CC)\,.
\eqn
\begin{rem}
The Cayley element defines the Cayley transformation $\frakp_+
\supset \Dd \to \Hh \subset \frakp_+$, which sends $\Dd$ to a Siegel
domain $\Hh$, which, if $\Xx$ is of tube type, (see Definition~\ref{defi:maxi-polydisk}) is a tube domain of the form 
$V \oplus i \Omega$.
\end{rem}
The automorphism $\Ad(c)$ of $\frakg_\CC$
is of order $4$ if $\Xx$ is of tube type and of order $8$ if $\Xx$ is not of tube type. 
When $\Xx$ is not of tube type 
$\Ad(c)^4$ is an involution of $\frakg_\CC$ which preserves
$\frakg$ and commutes with the Cartan involution of 
$\frakg = \frakk \oplus \frakp$ (see e.g. \cite[Theorem~4.9]{Koranyi_Wolf_65_annals}). 

We denote by $\frakg_T \subset \frakg$ the fix points of $\Ad(c)^4$ in
$\frakg$ and let $\frakg_T = \frakk_T \oplus \frakp_T$ be its Cartan
decomposition. Then the corresponding Hermitian symmetric space
$\Xx_T$ is of tube type. Furthermore $\Xx_T$ is isometrically and
holomorphically embedded into $\Xx$, the rank of $\Xx_T$ equals the
rank of $\Xx$ and as a bounded symmetric domain $\Xx_T$ is realized as 
\bqn
\Dd_T = \Dd \cap \frakp_T^+\,, 
\eqn
where $\frakp_T^+$ are the fixed points of $\Ad(c)^4$ in $\frakp_+$.

Note that the maximal standard polydisk $\Pp_0$ is contained in  $\Dd_T$ 
(see (\ref{eq:poly})), hence also
$\sum_{j=1}^r E_{\g_j} \in \cs_{\Dd_T} \subset \cs_\Dd$.  
Moreover for the polynomial $h_\Dd$ which is related to the Bergmann
kernel by Equation (\ref{eq:bergmann}) we have 
\bqn
h_{\Dd_T} = {h_\Dd}_{|_{\frakp_T^+}}\,.
\eqn

This implies in particular that $\Dd_T^{(3)} = \Dd^{(3)} \cap
(\frakp_T^+)^3$ and $\beta_{\Dd_T} =
{\beta_\Dd}_{|_{\overline{\Dd_T}^{(3)}}}$.

\begin{lemma}
$\Dd_T$ is a maximal (with respect to inclusion) subdomain of tube
type in $\Dd$.
\end{lemma}
%\begin{proof}
%Assume that $\Dd'$ where another subdomain of tube type with 
%$\Dd_T \subset \Dd' \subset \Dd$. Since $\Pp \subset \Dd_T \subset
%\Dd'$ is a maximal polydisk in both symmetric domains we have that 
%the Cayley element of $\Dd'$ equals the 
%Cayley element of $\Dd_T$, which is the Cayley element of $\Dd$. Since
%$\Dd'$ is of tube type $\Ad(c)_{|_{\Dd'}}$ is of order $4$ and $\Dd'
%\subset \Dd \cap \frakp_T^+ = \Dd_T$.
%\end{proof}

\subsection{The Shilov boundary and Tube Type Domains}
It is well known that the structure of the Shilov boundary $\cs_\Dd$ detects
whether $\Dd$ is of tube type or not, see for example
\cite[Theorem~4.9]{Koranyi_Wolf_65_annals}. Similarly the behavior of
the restriction of the Bergmann cocycle to the Shilov boundary detects
whether $\Dd$ is of tube type or not when $\Dd$ is irreducible.
In the general case we have:

\begin{prop}\cite[Corollary~3.10]{Burger_Iozzi_Wienhard_kahler}\label{prop:bergmann_tubetype}
Let 
\bqn
{\cs}^{(3)}:=\left\{ (z_1, z_2, z_3) \in {\cs}^3 \, :\, (z_i, z_j)
\in {\cs}^{(2)} \text{ for all } i\neq j\right\}
\eqn 
the space of triples of pairwise transverse points in $\cs$. Then $\cs^{(3)} \subset
\ol{\Dd}^{[3]}$ and the Bergmann cocycle $\beta_\Dd$ is well defined
and continuous on $\cs^{(3)}$. Furthermore,
\be
\item If $\Dd$ is of tube type, then 
\bqn
\b_\Dd(\cs^{(3)}) = \left\{-\frac{\r_\Dd}{2},
    -\frac{\r_\Dd}{2} +1 , \cdots , \frac{\r_\Dd}{2}-1, \frac{\r_\Dd}{2}\right\}\,.
\eqn
\item If $\Dd$ is irreducible and not of tube type, then
\bqn
\b_\Dd(\cs^{(3)}) 
=  \left[-\frac{\r_\Dd}{2}, \frac{\r_\Dd}{2}\right]\,.
\eqn
\ee
\end{prop}

Important for our later considerations is the relation 
between transverse pairs in the Shilov boundary 
and maximal subdomains of tube type. 

If $(x,y) \in \cs_\Dd^{(2)}$ we define 
\bqn
(\cs_\Dd)_{x,y} := \left\{ z \in \cs_\Dd \, :\, (z,x) \in \cs^{(2)}, \,
(z,y) \in \cs^{(2)}\right\} 
\eqn
to be the set of points $z \in \cs$ which are transverse to $x$
and to $y$. This is an open and dense set in $\cs_\Dd$.
In the following we shall denote $\sum_{j=1}^rE_{\gamma_j}$ by $E_\Lambda$,
where, as in \S~\ref{sec:shilov}, $\Lambda$ refers 
to the set $\{\gamma_1,\dots,\gamma_r\}$ of strongly orthogonal roots.
The following lemma is crucial and follows immediately from the case
in which $\Dd$ is irreducible, which was proven by Clerc and \O rsted
as the first step in the proof of Theorem~4.7 in \cite{Clerc_Orsted_2}.

\begin{lemma}\cite{Clerc_Orsted_2}\label{lem:co}
Let $z \in \cs_\Dd$ be transverse to $E_\Lambda$ and $-E_\Lambda$ with
\bqn
\big|\beta_\Dd(z,E_\Lambda,-E_\Lambda) \big|=\frac{\r_\Dd}{2}\,.
\eqn 
Then $z\in \cs_{\Dd_T}$.
\end{lemma}

If $(x,y)\in\cs_\Dd^{(2)}$ and $\Dd'\subseteq\Dd$ is a subdomain of Hermitian type, we define
\bq\label{eq:mxy}
  M_{x,y}(\Dd')
:=\left\{z\in\big(\cs_{\Dd'}\big)_{x,y}:\,
    \big|\beta_{\Dd'}(z,E_\lambda,-E_\Lambda)\big|=\frac{\r_\Dd}{2}\right\}\,.
\eq
From the above lemma we now deduce:

\begin{prop}\label{prop:maximalsubset} With the above notation we have that:
\bq\label{eq:sets}
M_{E_\Lambda,-E_\Lambda}(\Dd)=M_{E_\Lambda,-E_\Lambda}(\Dd_T)
\eq
and $\cs_{\Dd_T}$ is the real Zariski closure in $\cs_\Dd$ of 
$M_{E_\Lambda,-E_\Lambda}(\Dd)$.
\end{prop}

\begin{proof} Equation~\eqref{eq:sets} follows from Lemma~\ref{lem:co}
and the fact that $\beta_\Dd|_{\big(\cs_{\Dd_T}\big)^{(3)}}=\beta_{\Dd_T}$.

Since $\Dd_T$ is of tube type, then $M_{E_\Lambda,-E_\Lambda}(\Dd_T)$
is a nonempty open subset of $\cs_{\Dd_T}$,
and hence \eqref{eq:sets} implies that the Zariski closure of
$M_{E_\Lambda,-E_\Lambda}(\Dd)$ is $\cs_{\Dd_T}$
\end{proof}

%\begin{cor}\label{cor:maximalsubset}
%Let $(x,y)=\left(\sum_{j=1}^r E_{\g_j} ,-\sum_{j=1}^r E_{\g_j}\right) \in \cs_{\Dd_T}^{(2)}\subset \cs_\Dd^{(2)}$. 
%then
%\bqn
%\left \{ z \in (\cs_\Dd)_{xy} \, :\, 
%\left|\beta_\Dd\left( z,x ,y\right) \right|=\frac{\r_\Dd}{2}\right\} = 
%\left\{ z \in (\cs_{\Dd_T})_{xy} \, :\, 
%\left|\beta_{\Dd_T}\left( z,x ,y\right) \right|=\frac{\r_\Dd}{2}\right\}\,.
%\eqn 
%And since $\Dd_T$ is of tube type, the latter is a nonempty open
%subset of $(\cs_{\Dd_T})_{xy}$. 
%\end{cor}
%\begin{proof}
%The first is an immediate consequence of the previous lemma. The
%second statement follows from Proposition~\ref{prop:bergmann_tubetype}.
%\end{proof}
%
%\begin{cor}
%$\cs_{\Dd_T}$ is the real Zariski closure of 
%\bqn
%\left\{ z \in
%(\cs_{\Dd})_{xy} \, |\, 
%\left|\beta_{\Dd_T}\left( z,x ,y\right) \right|=\frac{\r_\Dd}{2}\right\}\,, 
%\eqn
%where  $(x,y)=\left(\sum_{j=1}^r E_{\g_j} ,-\sum_{j=1}^r E_{\g_j}\right)$.
%\end{cor}

\begin{prop}\label{prop:transversepoints_tubetype}
Let $\Dd$ be a Hermitian symmetric space and $(x,y) \in \cs_\Dd^{(2)}$ be
a pair of transverse points in its Shilov boundary. Then there
exists a unique maximal subdomain $T_{xy} \subset \Xx$ of tube type
with $x,y \in \cs_{T_{xy}}$.
Moreover $\cs_{T_{x,y}}$ is the real Zariski closure in $\cs_\Dd$ of 
$M_{x,y}(\Dd)$.
\end{prop}

\begin{proof} Observe that $(E_\Lambda,-E_\Lambda)$ are in $\cs_{\Dd_T}$
and $\Dd_T$ is a maximal subdomain of tube type.  Moreover,
since $G_\Dd$ acts transitively on $\cs_\Dd^{(2)}$, we obtain
the existence statement for every pair $(x,y)\in\cs_\Dd^{(2)}$.

Concerning uniqueness, we may assume,
again by transitivity of the $G_\Dd$-action on $\cs_\Dd^{(2)}$,
that $(x,y)=(E_\Lambda,-E_\Lambda)$.  
Let thus $\Dd'\subset\Dd$ be a maximal subdomain of tube type
with $(E_\Lambda,-E_\Lambda)\in\cs_{\Dd'}^{(2)}$.
Since $\r_\Dd=\r_{\Dd'}$, 
we have that $\beta_\Dd|_{\big(\cs_{\Dd'}\big)^{(3)}}=\beta_{\Dd'}$
and hence
\bqn
        M_{E_\Lambda,-E_\Lambda}(\Dd')
\subset M_{E_\Lambda,-E_\Lambda}(\Dd)
      = M_{E_\Lambda,-E_\Lambda}(\Dd_T)
\eqn
which implies, upon taking the real Zariski closure and using 
Proposition~\ref{prop:maximalsubset}, that $\cs_{\Dd'}\subset\cs_{\Dd_T}$.
On the other hand $\dim\Dd'=\dim\Dd_T$, which,
since $\Dd'$ and $\Dd_T$ are of tube type, implies that
$\dim\cs_{\Dd'}=\dim\cs_{\Dd_T}$ and,
together with the previously established inclusion,
that $\cs_{\Dd'}=\cs_{\Dd_T}$; this then implies
by Proposition~\ref{prop:3.11} that $\Dd'=\Dd_T$.
\end{proof}

\begin{rem}
One could prove the uniqueness in
Proposition~\ref{prop:transversepoints_tubetype}  
also by considering the Lie algebra of the stabilizer of 
$(E_\Lambda,-E_\Lambda)\in\cs_\Dd^{(2)}$, 
but for us the characterization of $\cs_{T_{xy}}$ 
obtained as a byproduct of the proof is essential.
\end{rem}

Let $\Tt_\Xx$ be the space of maximal tube type subdomains in $\Xx$. 
Then, since all maximal subdomains of tube type are conjugate, 
$\Tt_\Xx$ is a homogeneous space under $G_\Xx$. The map 
\bqn
\cs_\Xx^{(2)} \to \Tt_\Xx
\eqn 
provided by Proposition~\ref{prop:transversepoints_tubetype} is a $G_\Xx$-equivariant map 
between $G_\Xx$-homogeneous spaces and hence is real analytic. 

\subsection{Structure Theorem for Tight Embeddings, II}
\begin{thm}\label{thm:tight_embedding_tubetype} 
Let $H, G$ be Lie groups of Hermitian type with associated symmetric spaces $\Yy$ and $\Xx$. Let 
$\rho: H\to G$ be a continuous tight homomorphism 
and $f: \Yy \to \Xx$ the induced $\rho$-equivariant tight map.
Then:
\begin{enumerate}
\item If $\Yy$ is of tube type, then there exists a unique maximal tube type subdomain $T\subset \Xx$ 
such that $f(\Yy) \subset T$ and $\rho(H)$ preserves $T$. 
\item If $\rho$ has finite kernel and $\Xx$ is of tube type, then $\Yy$ is of tube type.
\end{enumerate}
\end{thm}

We shall need the following

\begin{lemma}\label{lem:5.11} 
Under the assumption of Theorem~\ref{thm:tight_embedding_tubetype},
let $\check f:\cs_\Yy\to\cs_\Xx$ be the (continuous) equivariant map
given by Theorem~\ref{thm:tight_embedding_one}, and let
\bqn
f^\ast(\omega_\Xx)=\sum_{i=1}^n\lambda_i\omega_{\Yy,i}\,.
\eqn
Then for all $(x,y,z)\in\cs^{(3)}$, we have
\bqn
 \beta_\Xx\big(\check f(x),\check f(y),\check f(z)\big)
=\sum_{i=1}^n\lambda_i\beta_{\Yy_i}(x_i,y_i,z_i)\,.
\eqn
In particular, if $f$ is moreover positive and
$\beta_\Yy(x,y,z)=\frac{\r_\Yy}{2}$, then
\bqn
\beta_\Xx\big(f(x),f(y),f(z)\big)=\frac{\r_\Xx}{2}\,.
\eqn
\end{lemma}

\begin{proof} Let $0\in\Dd_\Yy$ and $r_{0,x},r_{0,y}$ and $r_{0,z}$
be the geodesic rays given by Proposition~\ref{prop:3.4};
then we know that 
\bq\label{eq:5.3}
\check f(x)=\lim_{t\to\infty}f\big(r_{0,x}(t)\big),
\eq
and analogously for $y$ and $z$. Writing $x=(x_1,\dots,x_n)$ and $y$ and $z$
in coordinates, we have
\bqn
\ba
  \beta_\Xx\big(f\big(r_{0,x}(t)\big),f\big(r_{0,y}(t)\big),f\big(r_{0,z}(t)\big)\big)
&=\int_{\Delta\big(r_{0,x}(t),r_{0,y}(t),r_{0,z}(t)\big)}
   f^\ast(\omega_\Xx)\\
&=\sum_{i=1}^n\lambda_i\beta_{\Yy_i}\big(r_{0,x_i}(t),r_{0,y_i}(t),r_{0,z_i}(t)\big)
\ea
\eqn
and by using \eqref{eq:5.3} and the fact that the normalized Bergmann cocycles
extend continuously to $\cs_{\Yy_i}^{(3)}$ and $\cs_{\Xx}^{(3)}$,
we conclude the first claim.

Assume now that $f$ is positive;
then if $\beta_\Yy(x,y,z)=\frac{\r_\Yy}{2}$, then
$\beta_{\Yy_i}(x_i,y_i,z_i)=\frac{\r_{\Yy_i}}{2}$ and hence
\bqn
 \beta_\Xx\big(\check f(x),\check f(y),\check f(z)\big)
=\sum_{i=1}^n\lambda_i\frac{\r_{\Yy_i}}{2}
=\frac{\r_\Xx}{2}\,,
\eqn
where the last equality follows from Corollary~\ref{cor:tight_rank} 
and the fact that $f$ is positive, that is $\lambda_i\geq0$.
\end{proof}

\begin{proof} By changing complex structure on $\Yy$ we may assume
that $f$ is positive (Lemma~\ref{lem:makepositive}).
Let $(\cs_\Yy)_x$ be the set of points in $\Yy$ transverse to $x$
so that $(\cs_\Yy)_{x,y}=(\cs_\Yy)_x\cap(\cs_\Yy)_y$, 
and let us consider the set $M_{x,y}(\Yy)$ defined in \eqref{eq:mxy}.
Let $\check f:\cs_\Yy\to\cs_\Xx$ be the equivariant extension of $f$
given by Theorem~\ref{thm:tight_embedding_one}.
Since $f$ is tight and positive, we have that for every $z\in M_{x,y}(\Yy)$
\bqn
 \big|\beta_\Xx\big(\check f(x),\check f(y),\check f(z)\big)\big|
=\frac{\r_\Yy}{2}
\eqn
(see Lemma~\ref{lem:5.11}), and hence
\bqn
\check f(z)\in\cs_{T_{\check f(x),\check f(y)}}
\eqn
by Proposition~\ref{prop:maximalsubset}, and thus,
\bqn
T_{\check f(x),\check f(y)}=T_{\check f(x),\check f(z)}
\eqn
by the uniqueness statement in Proposition~\ref{prop:transversepoints_tubetype}.

Let $\Tt_\Xx$ be the conjugacy class of maximal tube type domains
seen as a G-homogeneous space and hence as a real analytic variety.
The map
\bq\label{eq:5.12}
\ba
(\cs_\Yy)_x&\longrightarrow\quad\Tt_\Xx\\
z\quad &\mapsto T_{\check f(x),\check f(z)}
\ea
\eq
is real analytic and constant on the subset $M_{x,y}(\Yy)\subset(\cs_\Yy)_x$;
but since $\Yy$ is of tube type, $M_{x,y}(\Yy)$ is open,
$(\cs_\Yy)_x$ is connected and hence the map \eqref{eq:5.12}
is constant on $(\cs_\Yy)_x$.

Let now $(x_1,y_1)$ and $(x_2,y_2)$ be arbitrary elements in $(\cs_\Yy)^{(2)}$
and choose $z\in(\cs_\Yy)_{x_1,x_2}$.  Then we have
\bqn
 T_{\check f(x_1),\check f(y_1)}
=T_{\check f(x_1),\check f(z)}
=T_{\check f(x_2),\check f(z)}
=T_{\check f(x_2),\check f(y_2)}
\eqn
which shows that the map
\bqn
\ba
(\cs_\Yy)^{(2)}&\longrightarrow\quad\Tt_\Xx\\
(x,y)\,\,\,&\mapsto T_{\check f(x),\check f(y)}
\ea
\eqn
is constant and hence its constant value $T\subset\Xx$ is $\rho(H)$-invariant.
We deduce also, since $\check f(x)\in\cs_{T_{\check f(x),\check f(y)}}$,
that $\check f(\cs_\Yy)\subset\cs_T$.

Now, by Theorem~\ref{thm:tight_embedding_one} we know that 
the centralizer of $\rho(H)<G$ in $G$ is compact;
this implies that, given any maximal compact subgroup $K<H$, 
there is a unique point $x_K\in\Xx$ which is $\rho(K)$-fixed.
Since $\rho(H)$ leaves $T$ invariant, 
this implies that $x_K\in T$ and hence that $f(\Yy)\subset T$.

For the second statement, observe that if $\Xx$ is of tube type, 
then $\beta_\Xx$ takes on finitely many values; 
since $\rho$ has finite kernel, this implies that $\lambda_i\neq0$
for all $1\leq i\leq n$ and hence each $\beta_{\Yy_i}$
takes on finitely many values on $\cs_{\Yy_i}^{(3)}$,
and this, together with Proposition~\ref{prop:bergmann_tubetype}, 
implies that $\Yy_i$ is of tube type.
\end{proof}

%% file: fiveandahalf.tex
\section{Extensions to Groups of Type (RH)}\label{sec:fiveandahalf}
Here we indicate the argument extending Theorems~\ref{thm:tight_embedding_one} and 
\ref{thm:tight_embedding_tubetype} to Lie groups of type (RH).
The study of tight homomorphisms of groups of type (RH) can be reduced
to the study of homomorphisms of groups of Hermitian type. In fact,
let $G_1$, $G_2$ be groups of type (RH) and
let $\rho: G_1 \to G_2$ be a continuous homomorphism.
%so that $\|\rho_{\rm b}^*(\kappa_{G_2}^{\rm b})\| = \|\kappa_{G_2}^{\rm b}\|$.
We have the inclusion $\rho(G_{1,nc})\subset G_{2,nc}$
and hence the commutative diagram
\bqn
\xymatrix{
G_1 \ar[r]^{\rho} & G_2 \\
G_{1,nc} \ar[u] \ar[r]^{\rho} &G_{2,nc} \ar[u]
}
\eqn
hence 
 $\rho^*_{\rm b}(\kappa_{G_2}^{\rm b})_{|_{G_{1,nc}}} 
= (\rho)^*_{\rm b} ({\kappa_{G_2}}_{|_{G_{2,nc}}})$.
From this and Corollary ~\ref{cor:2.3} we deduce the equalities
\bqn
 \|\rho^*_{\rm b}(\kappa_{G_2}^{\rm b})\| 
=\left\|\rho^*_{\rm b}\big({\kappa_{G_2}^{\rm b}}_{|_{G_{2,nc}}}\big)\right\| 
    \quad \text{and} \quad
\|\kappa_{G_2}^{\rm b}\| = \left\|{\kappa_{G_2}}_{|_{G_{2,nc}}}\right\|\,,
\eqn
from which it follows that $\rho$ is tight if and only of $\rho_{|_{G_{1,nc}}}$ is tight.

From this and Theorems~\ref{thm:tight_embedding_one} and \ref{thm:tight_embedding_tubetype}
we readily deduce the following:

\begin{thm}\label{thm:tight_embedding_one.2} 
Let $H, G$ be Lie groups of type (RH), 
$\rho: H \to G$ a continuous tight homomorphism 
and $f: \Dd' \to \Dd$ the corresponding $\rho$-equivariant tight
totally geodesic map. 
Then $f$ extends continuously to a $\rho$-equivariant map  
\bqn
\check{f}: \cs_{\Dd'} \to \cs_{\Dd}\,.
\eqn 
Moreover the centralizer $\Zz_G\big(\rho(H)\big)$ is compact.
\end{thm}

\begin{thm}\label{thm:tight_embedding_tubetype.2} 
Let $H, G$ be Lie groups of type (RH) with associated symmetric spaces $\Yy$ and $\Xx$. 
Let $\rho: H\to G$ be a continuous tight homomorphism 
and $f: \Yy \to \Xx$ the induced $\rho$-equivariant tight map.
Then:
\begin{enumerate}
\item If $\Yy$ is of tube type, 
then there exists a unique maximal tube type subdomain $T\subset \Xx$ 
such that $f(\Yy) \subset T$ and $\rho(H)$ preserves $T$. 
\item If $\rho$ has compact kernel and $\Xx$ is of tube type, then $\Yy$ is of tube type.
\end{enumerate}
\end{thm}

%%% Local Variables: 
%%% mode: latex
%%% TeX-master: "tight"
%%% End: 

%% file: structure_two.tex
\section{Structure Theorem for Tight Homomorphisms}\label{sec:structure_two}
In this section we prove the main structure theorem for tight
homomorphisms.
\begin{thm}\label{thm:structure_hom}
Let $L$ be a locally compact group, $\gG$ a connected algebraic group
defined over $\RR$ such
that $G= \gG(\RR)^\circ$ is of type (RH), and $
\rho: L \to G$ a continuous tight homomorphism.
Then: 
\begin{enumerate}
\item The Zariski closure $\hH := \ol{\rho(L)}^Z$ is reductive.
\item The centralizer $\Zz_G(H)$ of $H:= \ol{\rho(L)}^Z(\RR)^\circ$ is
  compact.
\item The group $H$ is of type (RH) and the symmetric space $\Yy$ corresponding to $H$ is Hermitian.
\item There is a unique complex structure on $\Yy$ such that the
  embedding $H\to G$ is tight and positive. 
\end{enumerate}
\end{thm}
\begin{proof}
Set $H:= \ol{\rho(L)}^Z(\RR)^\circ$. Then the inclusion $H \to G$ is tight. Let 
\bqn
H= H^{ss}\cdot R
\eqn 
be the decomposition of $H$, where $R$ is the amenable radical and $H^{ss}$ 
is a semisimple connected Lie group with finite center and no compact factors. 
Then it follows from Corollary~\ref{cor:2.3}
that the inclusion $H^{ss} \to G$ is tight. 
Let 
\bqn
H^{ss} =  H_1\cdot \dots\cdot H_n
\eqn 
be the decomposition of $H^{ss}$ into almost simple factors, 
and let $H_1, \dots, H_l$, $l\leq n$ be the almost simple factors of $H^{ss}$ 
for which the restriction $\kgb|_{H_i} \in \hcb^2(H_i, \RR)$ is nonzero. 
Then 
\bqn
 \|\kgb\|
=\|\kgb|_{H^{ss}}\|
= \sum_{i=1}^l\|\kgb|_{H_i}\|
\eqn
and the 
inclusion 
\bqn
H_1  \dots  H_l \to G_{nc}
\eqn 
is
tight. 
Let $\Xx$ be the symmetric space associated to $G_{nc}$, 
$\Yy_i$ the symmetric space associated to $H_i$, $1\leq i\leq l$ and  
$\Yy_1 \times \dots \times \Yy_l \to \Xx$ the corresponding tight embedding.

Then, by Theorem~\ref{thm:tight_embedding_one}, 
the centralizer $\Zz_{G_{nc}}(H_1  \dots  H_l)$ is compact, 
which implies first that $\ell=n$, that is
$\Zz_{G_{nc}}(H^{ss})$ is compact and hence that $\Zz_G(H)$ is compact.
Now, it $H$ were not reductive, it would be contained in a proper
parabolic subgroup of $G$ and hence $\Zz_G(H^{ss})$ would be noncompact.
Hence $H$ is reductive and, since $\Zz_G(H^{ss})$ is compact
and $\Yy_1\times\dots\Yy_n$ is Hermitian symmetric, 
the group $H$ is of type (RH).  

Finally, (4) follows from Lemma~\ref{lem:makepositive}.
\end{proof}

%Then, by Theorem~\ref{thm:tight_embedding_one}, 
%the centralizer $\Zz_{G_{nc}}(H_1  \dots  H_l)$ is compact, hence 
%$\Zz_{G}(H_1  \dots  H_l)$ is compact. 
%This proves (2) since $\Zz_G(H) \subset \Zz_{G}(H_1  \dots  H_l)$.
%
%Furthermore this implies that $l=n$ and that $R$ is compact, 
%and hence that $H$ and $\hH$ are reductive. 
%Since $l=n$, every almost simple factor $H_i$ of the noncompact part $H_{nc} = H^{ss}$ 
%is of Hermitian type, hence $H$ is of type (RH) 
%and $\Yy = \Yy_1 \times \dots \times \Yy_n$ is Hermitian symmetric. 
%Thus we have (1) and (3). 
%
%By Lemma~\ref{lem:makepositive} there exists a unique complex structure $\Jj$ on $\Yy$ 
%which makes the tight embedding $\Yy \to \Xx$ positive, hence we have (4).
%\end{proof} 

From Theorem~\ref{thm:structure_hom} we can now deduce the following 
\begin{thm}\label{thm:map}
Let $\G$ be a countable discrete group with probability measure $\theta$
and let $\gG$ be a semisimple real algebraic group such that
$G:=\gG(\RR)^\circ$ is of type (RH).
If $(B,\nu)$ is a Poisson boundary for $(\G,\theta)$
and $\rho:\G\to G$ is a tight homomorphism,
then there exists a $\rho$-equivariant measurable map 
\bqn
\varphi:B\to\cs_\Xx\,.
\eqn
\end{thm}

\begin{proof} Let $\hH$ be the Zariski closure of $\rho(\G)$.
By Theorem~\ref{thm:structure_hom} the symmetric space $\Yy$ associated to 
$H:=\hH(\RR)^\circ$ is Hermitian symmetric and we fix a complex
structure such that the embedding $\Yy\to\Xx$ is tight and positive. 
Theorem~\ref{thm:tight_embedding_one} gives the existence of a $\rho$-equivariant 
map $\check f$ between the corresponding Shilov boundaries
\bq\label{eq:shilov}
\check f:\cs_\Yy\to\cs_\Xx\,.
\eq
Let $\qQ_H<\hH$ be a maximal parabolic subgroup defined over $\RR$ 
such that $\cs_\Yy\cong\hH(\RR)/\qQ_H(\RR)$, 
and let $\pP_H<\qQ_H$ be a minimal parabolic subgroup
defined over $\RR$ contained in $\qQ_H$, 
so that we have an equivariant map 
\bq\label{eq:surj-par}
\hH(\RR)/\pP_H(\RR)\twoheadrightarrow\hH(\RR)/\qQ_H(\RR)\cong\cs_\Yy\,.
\eq

Since $\rho:\G\to H$ has Zariski dense image
\cite[Theorem~4.7]{Burger_Iozzi_Wienhard_kahler} implies 
the existence of a $\rho$-equivariant measurable boundary map 
\bqn
\varphi_0:B\to\hH(\RR)/\pP_H(\RR)
\eqn
which composed with the maps in (\ref{eq:shilov}) and (\ref{eq:surj-par})
provides the $\rho$-equivariant map $\phi:B \to \cs$.
\end{proof}

%% file: criterion.tex
\section{A Criterion for Tightness}\label{sec_criterion}
To get a simple criterion when a totally geodesic embedding is tight 
let us recall the relation between totally
geodesic embeddings and Lie algebra homomorphisms. 

Let $G_1, G_2$ be connected semisimple Lie groups with finite center
and no compact factors and $\Xx_1,
\Xx_2$ be the corresponding symmetric spaces. Fix two base points $x_j \in \Xx_j$,
$j=1,2$, and let $\frakg_j = \frakk_j \oplus \frakp_j$ be the
corresponding Cartan decompositions of $\frakg_j = \mathrm{Lie}(G_j)$. 
Then every totally geodesic embedding $f: \Xx_1 \to \Xx_2$ with $f(x_1)
= x_2$ induces a Lie algebra homomorphism $\rho: \frakg_1 \to
\frakg_2$ which respects the Cartan decompositions. Conversely any Lie
algebra homomorphism $\rho: \frakg_1 \to \frakg_2$ respecting the
Cartan decompositions gives rise to a totally geodesic embedding $f:
\Xx_1 \to \Xx_2$ with $f(x_1) = x_2$. 

Let $\Xx$ be a Hermitian symmetric space with a fixed complex
structure $\Jj_\Xx$ and let $Z_{\Jj_\Xx} \in \Zz(\frakk)$ be the
element in the center of $\frakk$ such that
$\ad(Z_{\Jj_\Xx})_{|_\frakp}$ induces the complex structure $\Jj_\Xx$
on $T_{x_0}\Xx \cong \frakp$. 
The restriction of the Killing form $\frakB$ on $\frakg$ to $\frakk$ is
a negative definite symmetric bilinear form. Let 
\bqn
\frakk = \RR Z_{\Jj_\Xx} \oplus \RR Z_{\Jj_\Xx}^\perp
\eqn 
be the orthogonal decomposition of $\frakk$ with respect to $\frakB_{|_\frakk}$.
We identify $ \RR Z_{\Jj_\Xx}$ with $\RR$ by sending $ Z_{\Jj_\Xx}$
to $i$. Then the orthogonal projection onto  $\RR Z_{\Jj_\Xx} $
defines a homomorphism $\lambda_{Z_{\Jj_\Xx}} \in \hom(\frakk, i\RR)$.

To relate tightness of a totally geodesic embedding with properties of
the corresponding Lie algebra homomorphism we make use of the isomorphism 
\bq\label{eq_homcom}
\hom(\frakk, i\RR) \to \Omega^2(\Xx)^{G_\Xx} \to \hcb^2(G_\Xx,\RR)\,,
\eq 
where the first map associates to a homomorphism $\lambda \in
\hom(\frakk, i\RR)$ the unique $G_\Xx$-invariant differential form on
$\Xx$ whose value at $x_0$ is 
\bqn (\omega_\lambda)_{x_0}(X,Y) :=
\frac{1}{4\pi i} \lambda\big([X,Y]\big)\,, 
\eqn 
for $X,Y \in \frakp \cong T_{x_0}\Xx$.

Let $\Xx_1, \Xx_2$ be two Hermitian symmetric spaces with complex
structures $\Jj_{\Xx_j}$ given by $Z_j = Z_{\Jj_{\Xx_j}}\in \Zz(\frakk_j)$, $j= 1,2$. Let 
$f: \Xx_1 \to \Xx_2$ be a totally geodesic embedding, $f(x_1) =
x_2$, and $\rho: \frakg_1 \to \frakg_2$ the corresponding Lie algebra
homomorphism. 
Let $\DD$ be the Poincar\'e disk and $Z_\DD \in \fraks\frako(2)
\subset \fraks\frakl(2,\RR)$ the element which induces the standard
complex structure on $\DD$. 
Let $d_j: \DD \to \Xx_j$, $j = 1,2$ be diagonal disks with $d_j(0) = x_j$ and 
$\rho_j:\fraks\frakl(2,\RR) \to \frakg_j$ the corresponding Lie algebra homomorphisms.

\begin{lemma}\label{lem:criterion}
The embedding $f: \Xx_1 \to \Xx_2$ is tight and positive if and only
if 
\bqn
\lambda_{Z_2}\big(\rho\big(\rho_1(Z_\DD)\big)\big) =
\lambda_{Z_2}\big(\rho_2(Z_\DD)\big)\,.
\eqn
\end{lemma}
\begin{proof}
We start by showing that since the embedding $d_1: \DD \to \Xx_1$ is tight and positive,
then $f:\Xx_1 \to \Xx_2$ is tight and positive if
and only if $h:=f\circ d_1 :\DD \to \Xx_2$ is tight and positive. 
In fact, let $\Xx_1=\Xx_{1,1}\times\dots\times\Xx_{1,n}$ be the decomposition of $\Xx$
into irreducible subspaces and let 
\bqn
f^\ast(\omega_{\Xx_2})=\sum_{i=1}^n\lambda_i\omega_{\Xx_1,i}\,.
\eqn
Then 
\bqn
 h^\ast(\omega_{\Xx_2})
=\sum_{i=1}^n\lambda_i d_1^\ast(\omega_{\Xx_1,i})
=\sum_{i=1}^n\lambda_i \r_{\Xx_{1,i}}\omega_\DD\,,
\eqn
where the last inequality follows from the fact that $d_1$ is tight and positive.
If $h$ is tight and positive, then 
\bqn
 h^\ast(\omega_{\Xx_2})=\r_{\Xx_2}\omega_\DD
\eqn
so that 
\bq\label{eq:tipo}
\sum_{i=1}^n\lambda_i \r_{\Xx_{1,i}}=\r_{\Xx_2}\,.
\eq
Since $f$ is norm decreasing then 
\bqn
\sum_{i=1}^n|\lambda_i| \r_{\Xx_{1,i}}\leq\r_{\Xx_2}\,,
\eqn
which together with \eqref{eq:tipo} implies that $f$ is positive and,
by Corollary~\ref{cor:tight_rank}, tight.

Let $\omega_{\lambda_{Z_2}} \in \Omega^2(\Xx_2)^{G_2}$ be the
differential form corresponding to $\lambda_{Z_2}$. 
Then, since $\Omega^2(\DD)^{\PSL(2,\RR)} = \RR \omega_\DD$ is one dimensional  
\bqn
 {h}^*\omega_{\lambda_{Z_2}} 
=\frac{\lambda_{Z_2}\big(\rho\big(\rho_1(Z_\DD)\big)\big)}
      {\lambda_{Z_2}\big(\rho_2(Z_\DD)\big)}d_2^*\omega_{\lambda_{Z_2}}\,.
\eqn 
But since $d_2$ is tight and positive, $h$ is tight and positive 
if and only if the proportionality constant is equal to $1$.
\end{proof}

Lemma~\ref{lem:criterion} gives a criterion for tightness which takes
on a particular nice form when $\Xx_2$ is of tube type. Recall from 
\cite[Proposition~3.12.]{Koranyi_Wolf_65_annals}\label{lem:tubetype}
that a Hermitian symmetric space $\Xx_2$ is of tube type 
if and only if there exists a diagonal disk $d_2:\DD \to \Xx_2$
such that the corresponding Lie algebra homomorphism satisfies $\rho_2(Z_\DD) = Z_2$.

\begin{cor}\label{cor:criterion_tube}
Let $\Xx_1, \Xx_2$ be Hermitian symmetric spaces of tube type. A
totally geodesic embedding $f: \Xx_1 \to \Xx_2$ is tight and positive
if and only if  
the corresponding Lie algebra homomorphism $\rho: \frakg_1 \to \frakg_2$ satisfies 
\bqn
\lambda_{Z_2}\big(\rho(Z_1)\big) = 1\,.
\eqn
\end{cor}
\begin{proof}
Since $\Xx_1, \Xx_2$  are of tube type, we can choose the tight
holomorphic disks $\d_j: \DD \to \Xx_j$, $j=1,2$ such that 
the corresponding Lie algebra homomorphisms 
$\rho_j: \fraks\frakl(2,\RR) \to \frakg_j$  satisfy $\rho_j(Z_\DD)= Z_j$. 
Then $\lambda_{Z_2}\big(\rho_2(Z_\DD)\big) = \lambda_{Z_2}(Z_2)= 1$  and 
Lemma~\ref{lem:criterion} implies the claim. 
\end{proof}

Let us compare this criterion for tightness with the characterization of Lie
algebra homomorphisms corresponding to  
holomorphic totally geodesic embeddings. 

\begin{defi}\cite{Satake_book}
%Suppose $\frakg_1, \frakg_2$ are simple Lie algebras of Hermitian type. 
%Then a 
A homomorphism $\rho:\frakg_1\to \frakg_2$ is said to be of type 
\begin{enumerate}
\item[$(\h1)$] if $\ad\big(\rho(Z_1)\big)= \ad(Z_2)$.
\item[$(\h2)$] if $\rho(Z_1)=Z_2$.
\item[$({\h2}^\prime)$] if $\rho$ is $(\h1)$ and the induced holomorphic totally 
geodesic map $\Dd_1\to \Dd_2$ maps the Shilov boundary of $\Dd_1$
into the Shilov boundary of $\Dd_2$.
\end{enumerate}
\end{defi}

Lie algebra homomorphisms of type $(\h1)$ are in one-to-one
correspondence with holomorphic totally geodesic embeddings $\Xx_1 \to \Xx_2$. 

With these definitions, Theorem~\ref{thm:tight_embedding_one} implies:

\begin{cor}
Assume that $f:\Xx_1\to \Xx_2$ is a holomorphic tight embedding. 
Then the corresponding Lie algebra 
homomorphism $\rho: \frakg_1\to \frakg_2$ is an $(\h2^\prime)$ homomorphism.
\end{cor}

Corollary~\ref{cor:criterion_tube} together with 
\cite[Proposition~10.12]{Satake_book} implies the following

\begin{cor}\label{H2_hom}
Suppose $\Xx_1, \Xx_2$ are Hermitian symmetric spaces of tube type. 
Then $f:\Xx_1\to \Xx_2$
is a tight and holomorphic embedding if and only if the 
corresponding homomorphism of Lie algebras 
$\rho:\frakg_1\to \frakg_2$ is an $(\h2)$-homomorphism.
\end{cor}

\begin{rem}
When $\frakg_1, \frakg_2$ are not of tube type, the property of being an
$(\h2)$-Lie algebra homomorphism does not imply tightness: for 
$\fraks\fraku(1,n) = \aut(V,h)$ the representations of $\fraks\fraku(1,n)$  
on $\Lambda^k(V)$ preserving the Hermitian form induced 
by $h$ are always $(\h2)$ \cite[page 188]{Satake_book}, but we will
see below that they are tight only for $k=1$. 
\end{rem}

\subsection{Examples}\label{sec:examples}

\begin{exo}\label{ex:h2} 
By Corollary~\ref{H2_hom} $(\h2)$ Lie algebra homomorphisms 
give examples of diagonal embeddings 
if the Hermitian symmetric spaces are of tube type. 
All $(\h2)$ Lie algebra homomorphisms 
were classified by Satake \cite{Satake_annals} and Ihara \cite{Ihara_65, Ihara_66}. 

If we are actually interested in the existence of tight homomorphism 
of Lie groups of Hermitian type the problem becomes more complicated 
since the Lie algebra homomorphism might only lift to a 
Lie group homomorphism of a finite cover of the adjoint Lie group. 
Satake showed in \cite{Satake_annals} (see also  \cite[Ch.~IV]{Satake_book}) 
that some of the those $(\h2)$ Lie algebra homomorphism lift to Lie group 
  homomorphisms, for example  
\bqn
&\tau:& \SU(n,n) \to \Sp(4n, \RR)\\
&\tau:& \SO^*(4n) \to \Sp(8n,\RR)\\
&\tau:& \Spin(2,n) \to \Sp(2m, \RR), \text{ where $m$ depends on
}n\mod 8\,, 
\eqn
are tight homomorphisms.
\end{exo}

\begin{exo}%\label{ex:irreducible}
An important and interesting tight embedding that 
is not holomorphic is the embedding of the Poincar\'e disk 
into the Siegel upper half space obtained 
from the 2n-dimensional irreducible representation 
$\fraks\frakl(2,\RR) \to \fraks\frakp(2n,\RR)$. 

\begin{prop}\label{irreducible_rep}
The homomorphism $\rho:\fraks\frakl(2,\RR)\to
\fraks\frakp(2n,\RR)$ given by the 
$2n$-dimensional irreducible representation of $\fraks\frakl(2,\RR)$ is tight.
\end{prop}

\begin{proof}
Let $Z_2$ and $Z_{2n}$ be 
  generators of the center of the maximal compact Lie subalgebras in
  $\fraks\frakl(2,\RR)$ and $\fraks\frakp(2n,\RR)$ respectively. 
Let $\lambda$ be the homomorphism $\hom(\frakk_{2n}, i\RR)$ given by the 
orthogonal projection onto $\RR\cdot Z_{2n}$. 
Then we have to determine $\lambda(\rho(Z_2))$. 

Let $V = \RR_{2n-1}[x,y]$ be the vector space of homogeneous
polynomials of degree $2n-1$ 
in two variables $x,y$, with 
a basis is given by $(P_0, \dots P_m)$, $m= 2n-1$, where $P_k(x,y)= x^{m-k}y^k$. 
The $2n$-dimensional 
irreducible representation of $\fraks\frakl(2,\RR)$ is given by the following action: 
Let $X=\begin{pmatrix} a & b \\ c & -a \end{pmatrix} \in \fraks\frakl(2,\RR)$, then 
\bqn
\rho(X) P_k(x,y)= a(m-2k) P_k + b(m-k) P_{k+1} + c k P_{k-1}\,.
\eqn
This action preserves the skew symmetric bilinear form $\< .,.\>$ on $V$, 
defined by 
$\< P_k, P_{l}\> = (-1)^k {\binom{k+l}{k}}^{-1} \delta_{m-k,l}$ 
%$\< P_k, P_{m-k}\> = (-1)^k {\binom{m}{k}}^{-1}$ 
and gives rise to the irreducible representation 
\bqn
\rho:\fraks\frakl(2,\RR) \to \fraks\frakp(2n,\RR)\,,
\eqn 
into the Lie algebra of the symplectic group ${\rm Sp}(V,\langle \,.\,,\,.\,\rangle)$.
The map $J$ defined by $J P_k= (-1)^k P_{m-k}$ gives a complex
structure on $V$ and the element in the center of 
$\frakk_{2n}\subset \fraks\frakp(2n,\RR)$ which induces the complex structure 
on $\frakp_{2n}\subset \fraks\frakp(2n,\RR)$ via the adjoint action is $Z_{2n}
=\frac{1}{2} J$.

The image of the element 
\bqn
Z_2 = \frac{1}{2} \begin{pmatrix} 0 & -1 \\ 1 & 0\end{pmatrix}
\eqn
is given by 
$\rho(Z_2) P_k = \frac{1}{2}((k-m)P_{k+1} + k P_{k-1})$. 
Decomposing $\rho(Z_2)=
\lambda\big(\rho(Z_2)\big)Z_{2n} \mod Z_{2n}^\perp$, 
we get 
\bqn
Z_{2n}  \rho(Z_2) = \frac{-\lambda}{4} \id_V \mod  
Z_{2n}^\perp
\,.
\eqn
Since $\tr \big(Z_{2n}  Z_{2n}^\perp\big)= 0$, 
we have that 
$\tr\big(Z_{2n}  \rho(Z_2)\big) 
= \frac{-\lambda}{4} \mathrm{dim}(V)$. 
Now 
$Z_{2n}\rho(Z_{2 }) P_k 
= \frac{1}{4} (-1)^{k+1} \big[(m-k) P_{m-k-1}- k P_{m-k+1}\big]$. 
Thus the diagonal terms are 
\bqn
\frac{1}{4} (-1)^{n} (2n-1-n+1)\,  \text{ for }\,  k={n-1} \quad \text{ and
} \quad \frac{1}{4} (-1)^n n \,  \text{ for }\,  k=n\,.
\eqn
Hence 
\bqn
\left|\tr\big(Z_{2n}  \rho(Z_{2 })\big)\right| =
\frac{n}{2} =\frac{1}{4} \mathrm{dim}(V)=|\tr(Z_{2n}
Z_{2n})|\,, 
\eqn 
and $|\lambda|=1$.
\end{proof}
\end{exo}

\begin{exo}%\label{ex:irreducible}
\begin{prop}\label{prop:irred_not_tight} The irreducible representation $\rho: \fraks\fraku(1,2) \to
   \fraks\fraku(2,4)$ is not tight.
\end{prop}
\begin{proof}
Let $V$ be a $3$-dimensional complex vector space with Hermitian form
of signature $(1,2)$.  
The irreducible representation $\rho: \fraks\fraku(1,2) \to
   \fraks\fraku(2,4)$ is the representation given by the action of
   $\fraks\fraku(1,2)$ on $\sym^2(V)$ with the induced Hermitian form.
Let $Z_{\fraks\fraku(2,4)}$ be the generator of the center of the maximal compact subalgebra on
$\fraks\fraku(2,4)$ and $Z_{\fraks\fraku(2,2)}$ 
the generator of the center of the maximal compact subalgebra
of the tightly embedded subalgebra $\fraks\fraku(2,2)\subset
\fraks\fraku(2,4)$. Then
$\tr(Z_{\fraks\fraku(2,4)}Z_{\fraks\fraku(2,2)})=-1$, so the
representation $\rho: \fraks\fraku(1,2) \to 
   \fraks\fraku(2,4)$ is tight if and only if 
$|\tr(Z_{\fraks\fraku(2,4)} \rho(Z_{\fraks\fraku(1,1)}))| =
   1$, where $Z_{\fraks\fraku(1,1)}$ is a generator of the center of  
the maximal 
compact subalgebra of $\fraks\fraku(1,1) \subset \fraks\fraku(1,2)$.
A direct computation shows that $|\tr(Z_{\fraks\fraku(2,4)}
\rho(Z_{\fraks\fraku(1,1)}))| = \frac{1}{6}$, thus $\rho$ is not
tight. 
\end{proof}
\end{exo}

\begin{exo}\label{nonex:complex}
Considering a complex vector space $V_\CC$ of dimension $(1+n)$ 
with a Hermitian form of signature $(1,n)$ as 
real vector space $V_\RR$ of dimension $(2+n)$ with a quadratic form of signature $(2,n)$ provides a 
natural embedding ${\rm SU}(1,n)\to {\rm SO}(2,2n)$. 
The corresponding embedding $\Hh_\CC^n \to \Xx_{2,2n}$ is holomorphic
but it is not tight. 
Since $\Hh_\CC^n$ is of rank one, the totally geodesic embedding extends 
continuously to a map of the boundary of $\Hh_\CC^n$, but its image 
does not lie in the Shilov boundary of $\Xx_{2,2n}$.
\end{exo}

\begin{exo}\label{ex:twosp4}
%Tightness indeed depends 
%on the complex structure. 
In ${\rm SL}(4,\RR)$ 
there are two copies of ${\rm Sp}(4, \RR)$
\bqn
{\rm Sp}(4, \RR)_A :=\{ g \in {\rm SL}(4,\RR) \,|\, g^* J g = J\}\\
{\rm Sp}(4, \RR)_B :=\{ g \in {\rm SL}(4,\RR) \,|\, g^* \tilde{J} g = \tilde{J}\}\,, 
\eqn
where $J= \begin{pmatrix} 0 & \id\\ -\id & 0 \end{pmatrix}$ and 
$\tilde{J}=\begin{pmatrix} 0 & \Lambda \\ -\Lambda & 0\end{pmatrix}$ 
with $\Lambda=\begin{pmatrix} 0 & 1\\ 1 & 0\end{pmatrix}$,  
which are conjugate by $s= \begin{pmatrix} \id & 0 \\ 0 & \Lambda\end{pmatrix}$. 
%satisfies $ \tilde{J}= sJs^{-1}$ 
%and thus conjugates $s^{-1} {\rm Sp}(4, \RR)_B s={\rm Sp}(4, \RR)_A$.
%
The two embeddings $i_{A,B}: {\rm SL}(2,\RR) \to {\rm SL}(4,\RR)$  
\bqn
i_A\left( \begin{pmatrix} a & b\\ c & d\end{pmatrix}\right) = \begin{pmatrix} a\,\id
  & b\,\id \\ c\,\id & d\,\id \end{pmatrix}
\eqn

\bqn
i_B\left( \begin{pmatrix} a & b\\ c & d\end{pmatrix}\right) = \begin{pmatrix} a\,\id
  & b\,\Lambda \\ c\,\Lambda & d\,\id \end{pmatrix}\,.
\eqn
are also conjugate by $s$. 
The images of ${\rm SL}(2,\RR)$ under these two embeddings are contained in  
${\rm Sp}(4, \RR)_A \cap {\rm Sp}(4, \RR)_B$. 
The embedding $i_A$ is tight and positive with respect to ${\rm Sp}(4, \RR)_A$ 
but totally real with respect to ${\rm Sp}(4, \RR)_B$.

The boundary $\partial\DD$ of $\DD$ is mapped under both embeddings 
into the Shilov boundaries 
$\cs_A$ respectively $\cs_B$. 
The totally real embedding 
extends to an embedding of ${\rm SL}(2,\CC)$ 
whereas the tight embedding extend to an embedding of ${\rm SO}(2,2)$ into 
${\rm Sp}(4, \RR)$ .
\end{exo}
%%% Local Variables: 
%%% mode: latex
%%% TeX-master: "tight"
%%% End: 

%% file: classification.tex
\section{Classification of Tight Embeddings of the Poincar\'e Disk}\label{sec:classification}
In this section we classify all tight embeddings $f: \DD \to \Xx$,
where $\Xx$ is any Hermitian symmetric space.

\begin{defi}\label{def:hull}
Let $\Xx$ be a Hermitian symmetric space of noncompact type. 
Let $V \subset \Xx$ be a subset. The {\em Hermitian hull} $\Hh(V)$ of $V$ is the smallest Hermitian 
symmetric subspace $\Hh(V)\subset \Xx$, such that
$V\subset \Hh(V)$.

If $\Xx_1, \Xx_2$ are Hermitian symmetric spaces and 
$f: \Xx_1 \to \Xx_2$ is a totally geodesic embedding, 
we denote by $\Hh(f) = \Hh(f(\Xx_1))$ the Hermitian hull of
$f(\Xx_1)\subset \Xx_2$.
\end{defi}

\begin{rem}
We make some observations.
\begin{enumerate}
\item Let $\Xx_1$ be irreducible, then $f: \Xx_1 \to \Xx_2$ is (anti)-holomorphic 
if and only if $\Hh(f) = f(\Xx_1)$.
\item If  $f: \Xx_1 \to \Xx_2$ is tight, then  $f: \Xx_1 \to \Hh(f)$
  is tight and $\Hh(f) \to \Xx_2$ is tight and holomorphic.
\item If  $f: \Xx_1 \to \Xx_2$ is tight, then $\Hh(f)$ is of tube type
  if and only if $\Xx_1$ is of tube type (Theorem~\ref{thm:tight_embedding_tubetype}).
\end{enumerate} 
\end{rem}

\begin{prop}\label{prop:hull}
Let $\Xx_1, \Xx_2$ be Hermitian symmetric spaces of tube type and $f:
\Xx_1 \to \Xx_2$ a tight embedding with corresponding Lie algebra
homomorphism $\rho: \frakg_1 \to \frakg_2$. Let $\Hh(f) \subset \Xx_2$
be the Hermitian hull and $\frakh \subset \frakg_2$ the Lie subalgebra
corresponding to the subgroup of Hermitian type $G_{\Hh(f)}$
determined by $\Hh(f)$.

Then $\frakh$ is the subalgebra generated by $\rho(\frakg_1)$ and
$Z_2$, where $Z_2\in \frakk_2$ is the element defining the complex structure on $\Xx_2$. 
\end{prop}
\begin{proof}
The Hermitian symmetric space 
$\Hh(f)$ is of tube type and the embedding $\Hh(f)\to \Xx_2$ is tight and holomorphic, 
therefore the 
corresponding Lie algebra homomorphism is an $(\h2)$ homomorphism
(Lemma~\ref{H2_hom}). In particular, 
the element $Z_\frakh$ in the center of the maximal compact subalgebra
of $\frakh$ defining the complex structure on $\Hh(f)$ equals $Z_2$.
Thus $\big\<\rho(\frakg_1),Z_2\big\> \subset \frakh$, and equality follows from
the minimality of $\Hh(f)$.
\end{proof}

Proposition~\ref{prop:hull} allows us to define in the above context 
($\Xx_1$ and $\Xx_2$ of tube type) the Hermitian hull 
of the Lie algebra homomorphism $\rho:\frakg_1 \to \frakg_2$ as 
\bqn
\Hh(\rho):= \big\<\rho(\frakg_1), Z_2\big\>\,.
\eqn

\begin{rem}
A similar characterization of the Hermitian hull is not 
true if $\Xx$ is not of tube type. 
Consider for example the canonical embedding of 
$\fraks\fraku (p,p) \to \fraks\fraku(p,q)$. It is holomorphic and tight, 
but the central element $Z_{p,q}$ of the maximal compact Lie algebra 
defining the complex structure on the symmetric space associated to $\SU(p,q)$ 
is not contained in $\fraks\fraku(p,p)$ if $p\neq q$.
\end{rem}

\begin{lemma}\label{lem:irreducible}
Let $\Xx$ be an irreducible Hermitian symmetric space and $f:\DD \to
\Xx$ a tight embedding with corresponding Lie algebra homomorphism 
$\rho:\fraks\frakl(2,\RR)\to \frakg$.

If $\Hh(\rho) = \frakg$, then $\frakg\cong\fraks\frakp(2n,\RR)$ and 
$\rho:\fraks\frakl(2,\RR)\to\frakg$ is the $2n$-dimensional irreducible representation.
\end{lemma}

\begin{proof}
Since $\DD$ is of tube type, $\Hh(\rho) = \frakg$ implies that
necessarily $\Xx$ is of tube type. 
Thus Proposition~\ref{prop:hull} gives that 
$\frakg = \big\<\rho(\fraks\frakl(2,\RR)), Z_{\Jj_\Xx}\big\>_\RR$. 
Let $\frakg_\CC$ be the complexification of $\frakg$ and $\rho_\CC:
\fraks\frakl(2,\CC) \to \frakg_\CC$ the complexification of
$\rho$, then $\frakg_\CC = \big\<\rho_\CC(\fraks\frakl(2,\CC)), Z_{\Jj_\Xx}\big\>_\CC$. 

By tightness $Z_{\Jj_\Xx}$ cannot lie in the centralizer of $\rho(\fraks\frakl(2,\RR))$ in
$\frakg$ and so $Z_{\Jj_\Xx}$ cannot lie in the centralizer
of $\rho_\CC(\fraks\frakl(2,\CC))$ in $\frakg_\CC$. Hence 
the centralizer $\Zz_{\frakg_\CC} \big(\rho_\CC\big((\fraks\frakl(2,\CC)\big)\big)$ is
trivial. 

This means that $\rho_\CC(\fraks\frakl(2,\CC))$ is a semiprincipal
three-dimensional simple subalgebra of $\frakg_\CC$. 
Semiprincipal subalgebras were classified by Dynkin, and we refer the reader 
to \cite{Dynkin, Onishchik_Vinberg_III} for more details. Using the classification by 
Dynkin (see \cite{Dynkin, Onishchik_Vinberg_III}), we consider all
possible cases of semiprincipal
three-dimensional simple subalgebras in $\frakg_\CC$ which are
complexifications of three-dimensional simple subalgebras of the
specific real form $\frakg$ of $\frakg_\CC$.
This case by case study gives the following result: 
\begin{enumerate}
\item  When $\frakg= \fraks\frakp(2n,\RR)$, $\frakg_\CC= \fraks\frakp(2n,\CC)$ 
the semiprincipal subalgebra 
$\frakh$ is given by the image of the  irreducible representation of  
$\fraks\frakl(2,\CC) \to \fraks\frakp(2n,\CC)$, hence $\rho:
\fraks\frakl(2,\RR) \to\fraks \frakp (2n,\RR)$ is the irreducible representation. 
\item When $\frakg= \fraks\fraku(n,n)$, $\frakg_\CC= \fraks\frakl(2n, \CC)$ 
the semiprincipal subalgebra is also given by 
the irreducible representation of $\fraks\frakl(2,\CC)$, which in dimension $2n$ 
is always contained 
in $\fraks\frakp(2n,\CC)$. Thus 
%with   
%$Z_{\fraks\fraku(n,n)} = Z_{\fraks\frakp(2n, \RR)}$, 
we have $ \Hh(\rho)=\<\rho(\fraks\frakl(2,\RR)), Z_\frakg\>=
\fraks\frakp(2n,\RR) \subset \fraks\fraku(n,n)$.
\item When $\frakg=\fraks\frako(2,2n-1)$, $\frakg_\CC= \fraks\frako(2n+1,\CC)$ the semiprincipal 
subalgebra $\frakh$ is the image of the irreducible representation 
of $\fraks\frakl(2,\CC)$. But any real irreducible representation 
of $\fraks\frakl(2,\RR)$ into $\fraks\frako(2,2n-1)$ is 
contained either in $\fraks\frako(2,3)\cong \fraks\frakp(4,\RR)$ or 
$\fraks\frako(2,1)\cong \fraks\frakp(2,\RR)$.
\item In the remaining cases $\frakg = \fraks\frako(2,2n),
  \fraks\frako^*(2n)$ or $\frake_{VII}$ 
there are no semiprincipal three dimensional subalgebras in $\frakg_\CC$ which are
complexifications of a real three dimensional simple subalgebra in $\frakg$.
\end{enumerate}

Summarizing we get the result: $\frakg\cong\fraks\frakp(2n,\RR)$ and
$\rho:\fraks\frakl(2,\RR)\to \frakg$ is given by the irreducible
representation of $\fraks\frakl(2,\RR)$.
\end{proof}

\begin{cor}\label{cor:hull_disk}
Let $f: \DD \to \Xx$ be a positive tight embedding with corresponding
Lie algebra homomorphism $\rho: \fraks\frakl(2,\RR) \to \frakg$. Then 
\bqn
\Hh(\rho) = \oplus_{i=1}^k \fraks\frakp(2n_i,\RR) \subset \frakg
\eqn
with $\sum_{i=1}^k n_i \leq \r_\Xx$ and $\rho_i: \fraks\frakl(2,\RR)
\to \fraks\frakp(2n_i, \RR)$ being the irreducible representation, and 
\bqn
\Hh(f) = \Yy_1\times \cdots \times \Yy_k\,,
\eqn
where $\Yy_k$ is a symmetric space associated to $\Sp(2n_i, \RR)$.
\end{cor}
\begin{proof}
The subalgebra $\Hh(\rho)$ is a semisimple Lie algebra of Hermitian
 type, so $\Hh(\rho)= \oplus_{i=1}^k \frakh_i$, where all simple factors $\frakh_i$ are of tube type. 
The representations  $\rho_i:\fraks\frakl(2,\RR)\to \frakh_i$ correspond again to tight embeddings 
with $\frakh_i= \Hh(\rho_i)$. Hence Lemma~\ref{lem:irreducible} implies
 the claim.
\end{proof}

%%% Local Variables: 
%%% mode: latex
%%% TeX-master: "tight"
%%% End: 

%% file: app.tex
\section{\ }
\subsection{The Norm of the Bounded K\"ahler Class}
We give here a proof of the following

\begin{thm}\label{thm:gromov_norm}
Let $M$ be a connected simple Lie group with finite center
and assume that its associated symmetric space $\Mm$ is Hermitian.
Let $\kappa_M^\mathrm{b}\in\hcb^2(M,\RR)$ be the continuous bounded class
given by the K\"ahler form associated to the Hermitian metric of
holomorphic sectional curvature -1.  Then
\bqn
\|\kappa_M^\mathrm{b}\|=\frac{\r_\Mm}{2}\,.
\eqn
\end{thm}

Since $\kappa_\Mm^\mathrm{b}$ is defined by the cocycle 
$c_{\omega_\Mm}$ which, according to Theorem~\ref{thm:dtco} is bounded by 
$\frac{\r_\Mm}{2}$, the inequality 
\bqn
\|\kappa_M^\mathrm{b}\|\leq\frac{\r_\Mm}{2}
\eqn
follows.  Observe that the opposite inequality cannot be immediately
deduced from the statement that 
\bqn
\|c_{\omega_\Mm}\|_\infty=\frac{\r_\Mm}{2}\,,
\eqn
since the norm $\|\kappa_\Mm^\mathrm{b}\|$ is the infimum of the supremum norms
over all bounded cocycles on $\Mm$ representing $\kappa_\Mm^\mathrm{b}$.

We shall proceed as follows: let $d:\DD\to\Mm$
be a diagonal disk (see Definition~\ref{defi:maxi-polydisk})
and $\rho:L\to M$ the corresponding homomorphism,
where $L$ is some finite covering of $\mathrm{SU}(1,1)$.
Then 
\bqn
d^\ast(\omega_\Mm)=\r_\Mm\omega_\DD
\eqn
and hence it follows from Lemma~\ref{lem:diagram} that 
\bqn
\rho_\mathrm{b}^\ast(\kappa_M^\mathrm{b})=\r_\Mm\,\kappa_L^\mathrm{b}\,.
\eqn
Since the pullback in continuous bounded cohomology is norm decreasing,
we have
\bqn
    \|\kappa_M^\mathrm{b}\|
\geq\|\rho_\mathrm{b}^\ast(\kappa_M^\mathrm{b})\|
   =\r_\Mm\|\kappa_L^\mathrm{b}\|\,,
\eqn
and it suffices to determine the values of $\|\kappa_L^\mathrm{b}\|$.
Thus the theorem will follow from the following:

\begin{prop} With the above notation we have that 
\bqn
\|\kappa_L^\mathrm{b}\|=\frac12\,.
\eqn
\end{prop}

\begin{proof}  Let $e:(\partial\DD)^3\to\{-1,0,1\}$
be the orientation cocycle on the circle $\partial\DD$.
we use the fact that the space of $L$-invariant alternating bounded measurable
cocycles on $\partial\DD$ is isometrically isomorphic to 
$\hcb^2(L,\RR)$ and that, under this isomorphism,
$e$ corresponds to $2\kappa_L^\mathrm{b}$, \cite{Iozzi_ern}.  
Thus, since $\|e\|_\infty=1$, we deduce that $\|\kappa_L^\mathrm{b}\|=\frac12$.
\end{proof}

\subsection{Surjection onto Lattices}

\begin{prop} Let $\G$ be a countable discrete group, $G$ a Lie group of Hermitian
type  and $\rho: \G \to G$ a homomorphism such that the image $\rho(\Gamma)$ is Zariski 
dense and the action of $\G$ on the Shilov boundary of the associated symmetric space is minimal.  
Then $\rho$ is tight.
\end{prop}

The proof of this proposition relies on functoriality properties of bounded cohomology. 
We use that the bounded continuous cohomology $\hb^2(L,\RR)$ in degree two 
of a locally compact group $L$
can be realized isometrically as the space $\Zz\la\big(B^3,\RR\big)^L$
of $L$-invariant bounded alternating $\linfty$ cocycles on any
space $(B,\nu)$ on which the $L$-action is amenable and mixing.
In particular if $G$ is a group of Hermitian type, then 
\bq\label{eq:fp}
\hb^2(G,\RR)\cong \Zz\la\big((G/P)^3,\RR\big)^G\,,
\eq
where $P<G$ is a minimal parabolic subgroup.  Likewise, 
we use that if $\G$ is the countable discrete group with a probability measure $\theta$ 
then a Poisson boundary $(B,\nu)$ for $(\G,\theta)$ always exists and then
\bqn
\hb^2(\G,\RR) \cong \Zz\la\big(B^3,\RR\big)^\G\,.
\eqn
%as well as of a lattice $\Lambda< G$, can be realized isometrically 
%as the cohomology of the complex of invariants of bounded alternating measurable functions 
%$\left( (\la\big((G/P)^\bullet,\RR\big)^G, d_\bullet\right)$, 
%respectively $\left( (\la\big((G/P)^\bullet,\RR\big)^\Lambda, d_\bullet\right)$, 
%on Cartesian products of $G/P$, where $P<G$ is a minimal parabolic subgroup.
%In particular, 
%Moreover, we use that if $\G$ be the countable discrete group with a probability measure $\theta$ 
%and $(B,\nu)$ is a Poisson boundary for $(\G,\theta)$, 
%then $\hb^2(\G,\RR) \cong \Zz\la\big(B^3,\RR\big)^\G$. 
For more details and proofs of the precise functoriality properties 
we refer the reader to  \cite[\S~4]{Burger_Iozzi_Wienhard_kahler}
and to the references therein.

\begin{proof}
We realize the Shilov boundary of the bounded domain realization $\Dd$ 
of the symmetric space associated to $G$ as $\cs = G/Q$. 
We fix a minimal parabolic subgroup $P<Q<G$ and 
denote by ${\rm pr}: G/P \to G/Q$ the canonical projection. 
If $\b_\Dd: \cs^3 \to \RR$ is the normalized Bergmann cocycle, 
then ${\rm pr}^* \b_\Dd: (G/P)^3 \to \RR$ is a cocycle in $\la\big((G/P)^3,\RR\big)^G$
representing the class $ \kgb\in \hcb^2(G,\RR)$. 
It follows from \eqref{eq:fp} that
\bqn
  \|\kgb\| =  \esssup_{x_1, x_2, x_3 \in G/P} {\rm  pr}^*\b_\Dd (x_1, x_2, x_3) \,.
\eqn
%where $r: \hcb^2(G, \RR) \to \hb^2(\Lambda, \RR)$ is the restriction
%map. 

Since the image of $\rho$ is Zariski dense, 
there exists a $\rho$-equivariant measurable boundary map $\phi:(B,\nu) \to G/P$ 
\cite[Theorem~4.7]{Burger_Iozzi_Wienhard_kahler}, and moreover $\phi^* {\rm pr}^*\b_\Dd$ 
represents $\rho_{\rm b}^*\kgb$ \cite[Proposition~4.6]{Burger_Iozzi_Wienhard_kahler}.

The essential image of ${\rm pr}\circ\phi:(B,\nu)\to G/Q$ is defined as
the support of the push-forward measure $({\rm pr}\circ\phi)^\ast(\nu)$
and is hence a closed $\rho(\Gamma)$-invariant subset which, 
by minimality of the $\Gamma$-action, must be equal to $G/Q$.
But then this implies that
\bqn
 \esssup_{x_1,x_2,x_3\in G/P}\big|{\rm pr}^\ast\beta_\Dd(x_1,x_2,x_3)\big|
=\esssup_{b_1,b_2,b_3\in G/P}\big|{\rm pr}^\ast\beta_\Dd\big(\phi(b_1),\phi(b_2),\phi(b_3)\big)\big|\,,
\eqn
and hence
\bqn
  \|\kgb\|
 =\|\rho_{\rm b}^*\kgb\|\,.
\eqn
\end{proof}

From the above proposition we obtain immediately the following

\begin{cor}\label{cor:latticetight}
Let $\G$ be a countable discrete group, $G$ a Lie group of Hermitian
type  and $\rho: \G \to G$ a homomorphism. 
If $\rho(\G)$ contains a lattice $\Lambda < G$, then $\rho$ is tight.
\end{cor}

\begin{cor}
Let $\Mod_g$ be the mapping class group of a closed oriented surface of genus $g$. 
Then the natural homomorphism $\rho: \Mod_g \to \Sp(2g,\RR)$ is tight. 
In particular if $\kgb\in \hcb^2(\Sp(2g,\RR))$ is the bounded K\"ahler class 
associated to the normalized K\"ahler form, 
then the norm of $\rho^*_{\rm b} \kgb \in \hb^2(\Mod_g, \RR)$ is $\frac{g}{2}$.
\end{cor}
\begin{proof}
The natural homomorphism $\rho: \Mod_g \to \Sp(2g,\RR)$ surjects onto 
$\Sp(2g,\ZZ)$, so $\|\rho_{\rm b}^*\kgb\| = \|\kgb\|$ 
which equals $\frac{g}{2}$ by Theorem~\ref{thm:gromov_norm}.
\end{proof}

%%% Local Variables: 
%%% mode: latex
%%% TeX-master: "tight"
%%% End: 